\documentclass[11pt]{article}

\usepackage{amsmath,
            amsfonts,
            amssymb,
            amsthm,
            amscd,
            latexsym,
            tangle-s,
	    }

\theoremstyle{plain}
\newtheorem{theorem}{Theorem}[section]
\newtheorem{proposition}[theorem]{Proposition}
\newtheorem{corollary}[theorem]{Corollary}
\newtheorem{lemma}[theorem]{Lemma}

\theoremstyle{definition}
\newtheorem{definition}{Definition}[section]
\newtheorem{example}{Example}[section]
\newtheorem*{acknowledgement}{Acknowledgements}

\theoremstyle{remark}
\newtheorem{remark}{Remark}[section]

\numberwithin{equation}{section}

\newcommand\be{\begin{equation}}
\newcommand\ee{\end{equation}}
\newcommand\bea{\begin{eqnarray}}
\newcommand\eea{\end{eqnarray}}
\def\barr#1{\begin{array}{#1}}
\def\earr{\end{array}}
\newcommand\ba{\begin{array}}
\newcommand\ea{\end{array}}

\newcommand{\field}[1]{\mathbb{#1}}
\newcommand\RR{\field{R}}

\newcommand\QQ{\field{Q}}
\newcommand\ZZ{\field{Z}}
\newcommand\NN{\field{N}}
\newcommand\Zn{\mathbb{Z}/n}

\newcommand\1{{1\mkern-5mu {\mathrm I}}}
\newcommand\2{{\underline2}}
\newcommand{\0}{{\underline0}}

\newcommand\DY[1]{\mathcal{DY}\left({#1}\right)}

\newcommand\cO[1]{{\mathcal{O}(#1)}}

\DeclareMathOperator\Aut{Aut}
\DeclareMathOperator\C{{\mathcal {C}}}
\DeclareMathOperator\Obj{Obj}
\DeclareMathOperator\im{im}
\DeclareMathOperator\coim{coim}
\DeclareMathOperator\tr{tr}

\DeclareMathOperator\Hom{Hom}
\DeclareMathOperator\End{End}

\DeclareMathOperator\Int{Int}
\DeclareMathOperator\Inv{Inv}

\newcommand\op{\mathrm{op}}
\newcommand{\pti}{^\vee}
\newcommand{\lpti}{{}^\vee}

\newcommand\ev{\mathrm{ev}}
\newcommand\coev{\mathrm{coev}}
\newcommand\id{\mathrm{id}}
\newcommand\ad{\mathrm{ad}}
\newcommand{\TC}{TC}
\newcommand{\RT}{RT}
\newcommand{\DRT}{DoubleRT}

\newcommand\intl[1]{{\textstyle\int_{#1}}}
\newcommand\cointl[1]{{\textstyle\int^{#1}}}
\newcommand\intr[1]{{\sideset{_{#1}}{}{\textstyle\int}}}
\newcommand\cointr[1]{{\vphantom{\textstyle\int}}^{#1}\!\!{\textstyle\int}}
\newcommand\miniint{\hbox{\scriptsize$\int$}}
\newcommand\cll{c^\ell_\ell}
\newcommand\clr{c^\ell_r}
\newcommand\crl{c^r_\ell}
\newcommand\crr{c^r_r}
\newcommand\ququ{\rho}
\newcommand\psipsi{\sigma}
\newcommand\monodromy[2]{\Omega^{#1}_{#2}}
\newcommand{\rib}{\nu}

\def\n-{\nobreakdash-\hspace{0pt}}

\let\tens\otimes
\let\xra\xrightarrow

\def\Atop#1#2{\genfrac{}{}{0pt}{}{#1}{#2}}
\def\I{{\mathrm I}}
\def\d{{\mathrm d}}
\def\m{{\mathrm m}}
\def\S{{[S]}}
\def\Op{{[\mathrm {op}]}}
\def\Vee{{[\vee]}}
\def\PPi{\overset\vee\Pi{}}

\author{
{Yuri Bespalov
 }
\and
{Thomas Kerler
 }
\and
{Volodymyr Lyubashenko
  \thanks{ Research was supported in part by NSF grant 530666.}
 }
\and
{Vladimir Turaev
}}
\title{\mbox{\huge \bf Integrals for braided Hopf algebras}}

\date{q-alg/9709020, September 12, 1997}

\begin{document}


\bibliographystyle{amsplain}
\thispagestyle{empty}
\maketitle

\begin{abstract}
Let $H$ be a Hopf algebra in a rigid braided monoidal category with
split idempotents.  We prove the existence of  integrals on (in) $H$
characterized by the universal property, employing results about Hopf
modules, and show that their common target (source) object $\Int H$ is
invertible.  The fully braided version of Radford's formula for the
fourth power of the antipode is obtained.  Connections of integration
with cross-product and transmutation are studied.
\end{abstract}

\begin{center}
\parbox{11 cm}{\small \tableofcontents}
\end{center}

\allowdisplaybreaks[1]

\bigskip

\section{Introduction}
\label{sec-intro}

The notion of integrals for finite-dimensional Hopf algebras was first
introduced by Larson and Sweedler in an attempt to generalize the
notion of  Haar measures on groups.  In their seminal paper~\cite{LS}
they prove that integrals always  exist for finite dimensional Hopf
algebras and give a variety of interesting applications, among them the
Maschke theorem for Hopf algebras.  Infinite-dimensional Hopf algebras
were also considered by Sweedler~\cite{Sweedler1}.

 By definition, e.g., a left  integral \emph{in} a Hopf algebra
$\,(H,m,\eta,\Delta,\epsilon,S)\,$
is an element, $\,l\in H\,$,  such that $\,hl=\epsilon(h)l\,$
 holds for any element, $\,h\in H\,$. The defining equation for a
right integral, $\,r\in H\,$, is analogously $\,rh=\epsilon(h)r\,$.
        From the uniqueness of  integrals, as proven in
\cite{LS}, one readily infers that the opposite regular actions also
leave the integrals invariant, only that now the one-dimensional
representation is described by a non-trivial character
$\alpha\in H^*\,$. I.e., we have   $\,lh=\alpha(h)l\,$ for all
$\,h\in H\,$.  The character $\alpha\in H^*\,$ is called the {\em
 module} on $\,H\,$, and plays a r\^ole similar to that of a modular
function for an invariant measure.  For the analogous dual definition
of integrals in $\,H^*\,$ (or {\em on} $\,H\,$) and the respective
modular element $a\in H\,$ the reader is referred to  Sweedler's
book~\cite{Sweedler}.
\medskip

 The theory of integrals, that was subsequently developed, turned into
a powerful instrument in the study of finite-dimensional Hopf algebras,
 and eventually revealed very rigid   structures inherent to
Hopf-algebras.  In particular, using  it, Radford proved that the order
of the antipode of a finite-dimensional Hopf algebra, $\,H\,$, is
finite \cite{Radford:antipode}.\nocite{Radford:trace} His approach is
based on a formula of his, which expresses the fourth power of the
antipode $S:H\to H$ in terms of the moduli $a\in H$, and $\alpha\in
H^*\,$.  More precisely, we have the following:
\[ S^{-4} = \ad^\alpha \circ \ad_a . \]
Here $\ad^\alpha\,$ and $\,\ad_a:H\to H\,$ are the usual Hopf algebra
automorphisms defined on an element,
$\,h\in H\,$, as $\ad_a.h= aha^{-1}\,$ and
$\,\ad^\alpha.h=
(\alpha\tens\id\tens\alpha^{-1}) (\Delta\tens\id) \Delta (h)\,$,
respectively.  Knowing that moduli are always group like, the only
thing that remains to be verified in order to infer now finiteness of
the order of $\,S\,$ is that group like elements in finite dimensional
Hopf algebras are of finite order, see \cite{Radford:antipode}.
\medskip

The classical theory of Larson, Sweedler, Radford, and others deals
only with algebras in the ordinary sense, meaning, the Hopf algebra
$\,H\,$ is always assumed to be a linear space over a given field
$\Bbbk$, and the defining operations are given as $\Bbbk$-linear maps.
For example, the multiplication, $m\,$, is a linear homomorphism of the
form $\,m\,\in\,Hom_\Bbbk\bigl(H\otimes_\Bbbk H,\,H\bigr)\,$.

In this paper we shall abandon the concept  of $\,H\,$ being a linear
space. Instead we shall assume that  $\,H\,$ is merely an abstract
object in a monoidal category, $\,{\cal C}\,$, with no further
structure of its own.  The content of the Hopf algebra structure thus
lies entirely in the operations of $\,H\,$, which are now expressed in
terms of morphisms in $\,{\cal C}\,$.  For example the multiplication
will be given by a morphism,
$\,m\,\in\,Hom_{\cal C}\bigl(H\otimes_{\cal C}H,\,H\bigr)\,$, where
$\,\otimes_{\cal C}\,$ stems from the monoidal structure of
$\,{\cal C}\,$.

The classical notion is recovered in the special case, where the
category $\,{\cal C}\,$ is abelian and admits a tensor fiber functor
$\,{\cal C}\,\longrightarrow\,Vect(\Bbbk)\,$, i.e., when $\,{\cal C}\,$
is Tannakian. If the category $\,{\cal C}\,$ is not Tannakian but still
abelian, a special Hopf algebra in $\,{\cal C}\,$ can often be found as
a special {\em coend} in $\,{\cal C}\,$.

We will, however, encounter numerous generic examples, in which
$\,{\cal C}\,$ is neither Tannakian nor abelian and not even additive,
but still Hopf algebra objects with interesting interpretations can be
extracted.  Some of these categories are defined combinatorially or
 purely topologically.
\medskip

 Our main goal in this article is to  extend  the above outlined
classical results for Hopf algebras to the framework of categories.
Specifically, we will show that for any braided Hopf algebra, $\,H\,$,
in a braided, rigid, monoidal category $\,{\cal C}\,$ with split
idempotents, the analogues of integrals and moduli are defined and,
similarly, satisfy existence and uniqueness assertions. Equipped with
these tools we shall then derive the respective generalization of
Radford's formula for these types of categories, which in
Theorem~\ref{thm-Radford-formula} will take on the form
\[ S^4\circ u^0_{-2} = (\ad^\alpha)^{-1}\circ \ad_a^{-1} \circ
\monodromy{(\Int H)}H\quad\,\in\;\Aut_{\cal C}(H)\,. \]
Here the morphism $\,u^0_{-2}\,$ is defined in Figure~\ref{Fig-u} and,
the $\,\monodromy{(\Int H)}H\,$ is constructed  from the square of the
braiding of $\,H\,$ with $\,\Int H\,$ -- the object of the integrals,
see  Lemma~\ref{lem-monodromy}. We shall also see that all three
 multiplicands in the right hand side are Hopf algebra automorphisms
commute with each other.
\medskip

Hopf algebras in tensor categories and their properties have been
investigated  by a number of other people.  Let us briefly review some
of the previous contributions to the subject, which we build up on and
generalize, and explain how our results are threaded into this
development.
\medskip

 The idea of an algebra in a tensor category, see~\cite{MacLane}, is
straightforwardly extended to that of a Hopf algebra, if that category
is also symmetric.  Most proofs of analogous assertions for categorical
  Hopf algebras in the symmetric case are formally identical to those
in the classical situation of a linear Hopf algebra over $\Bbbk$.  The
first proofs of existence and uniqueness of integrals in this setting
are due to Drinfel'd.  It is immediate from definitions that in a
symmetric category the elements $u^0_{-2}$ and
$\monodromy{(\Int H)}H$ are  trivial (i.e. units)  so that also the
categorical version of Radford's formula will differ from the original
one only in the interpretations of where the automorphisms live.

A nice way of illustrating the calculations, leading to the uniqueness
of integrals and the Radford formula, has been given by Kuperberg
\cite{Kuperberg}.  He uses diagrammatic techniques that are essential
for the construction of his 3-manifold invariants. In this language the
generalizations  to the symmetric categories  follow immediately
without much further explanation.

 The notion of  Hopf modules for ordinary Hopf algebras, which combines
an action with a compatible coaction on  the same space, can be
similarly generalized without much difficulty to  a Hopf algebra living
in an abelian, symmetric, monoidal category.  The classical results on
Hopf modules can be rederived again by ``imitation'' in the categorical
framework, where, similarly, the modules appear as abstract objects.
See  \cite{Sch:hopf}.
\medskip

In the theory  of Hopf algebras in braided categories a variety of
complications arise due to the fact that the isomorphism  between the
tensor product of two objects and the transposed product is no longer
canonical, i.e., also no longer coherent.  However, a  transposition
isomorphism is needed, e.g., to formulate the compatibility axiom
between multiplication and comultiplication for bialgebras.

If for the latter the braid isomorphism is employed one is left with
the definition of a  {\em braided} Hopf algebra, as it was given and
first studied by  Majid~\cite{Ma:bg}.

In the case where the braided tensor category is also rigid and
 abelian, integrals for such Hopf algebras were investigated by
 Lyubashenko \cite{Lyu:mod}, leading to general proofs for the
existence of integrals as well as for the  invertibility of the object
of integrals for a braided Hopf algebra. The results in \cite{Lyu:mod}
follow from the study of Hopf modules in abelian braided monoidal
categories.  This strategy generalizes the approach to integrals via
Hopf modules as proposed by Sweedler in \cite{Sweedler1}.
\medskip

As for symmetric categories we shall  prove here analogously the
existence and uniqueness of integrals for braided Hopf algebras in
braided monoidal categories, but in addition we shall drastically
weaken  the condition of abelianness. More precisely, we shall no
longer assume that the category has kernels or direct sums. As has
already been pointed out in remarks in  \cite{Bes:cross,BesDra:hopf} it
will suffice that our category has  {\em split idempotents}. This means
roughly that the category contains  for a given projection, $P\,$, in
an $End$-set the image of $\,P\,$ as an object, and that $\,P\,$
factors through this object.

 The property of split idempotent is central in   Karoubi's definition
of pseudo-abelian categories \cite{karoubi1}, which will be further
discussed later on.  It is easy to satisfy this condition taking the
{\em Karoubian envelope} of a category \cite{karoubi1}.
\medskip

Besides existence and uniqueness of integrals we will also prove under
our very general assumptions the invertibility of the object of
integrals $\Int H\,$, which, in the case of abelian, symmetric
categories, has already been pointed out  by Drinfel'd.

Furthermore, we will generalize the results on Hopf modules for
abelian, symmetric, monoidal category as in \cite{Sch:hopf} to any
braided monoidal category with  split idempotents.

Many of the algorithmic parts of the proofs will be done in a similar
diagrammatic language as in \cite{Kuperberg}. Only now the diagrams are
no longer plat graphs, but projections of graphs that distinguish
between over and under crossings.  As a result we will often encounter
additional special elements resulting from non trivial full twists,
such as $u^0_{-2}$ and   $\monodromy{(\Int H)}H$, which  enter the
Radford formula.

It is worth mentioning that our results apply to ordinary Hopf algebras
and Hopf modules over a commutative ring $\,R\,$, if they  are
projective $R$-modules of finite rank (cf. \cite{DelMil}).
\medskip

 Our interest in  braided Hopf algebras in non-abelian, braided
categories has  primarily been triggered  by realizing their crucial
r\^ole in the most recent discoveries in three dimensional topology
related to quantum physics.  In particular, integrals turned out to be
 the algebraic objects that are stringently  associated to  elementary
surgery data in the  construction of invariants of surgically presented
3-manifolds. Although the point of view of integrals was not used in
the more computational approach in  \cite{Tur:q3}, it is inevitable in
the construction of the non-semisimple analogues of invariants of
3-manifolds as in \cite{Lyu:3inv} and, more generally, 3-dimensional
topological quantum field theories as in  \cite{KerLyu}.

The fact that braided Hopf algebras and their integrals are inseparable
from 3-dimensional topology in this approach, naturally leads one to
identifying the torus with one hole as such a Hopf algebra in the
braided category of 3-dimensional cobordisms between one holed
 surfaces.  As a non trivial example of our generalized theory, we will
precisely identify the integrals and all ingredients to the Radford
formula for this braided Hopf algebra and category, via their explicit
presentations in a tangle category as in \cite{Ker:bridge}.

\bigskip

\subsection*{Summary of Contents:}
In  Section \ref{sec-prelim} we review the necessary preliminaries and
standard notations pertaining to rigid braided categories and braided
Hopf algebras.  We shall further define  categories with split
idempotents and discuss  basic properties of invertible objects.
  In Section \ref{sec-basic} results about Hopf modules are applied to
prove  the existence of ($\,\Int H$)-valued(-based) integrals on (in)
$\,H\,$. There the object $\Int H$ of $\,\Int H\,$ is  characterized by
a universal property implying uniqueness, and will turn out to be
invertible. The fully
  braided version of  Radford's formula for the fourth power of the
antipode is derived.  Section~\ref{sec-proofs}  starts with an
exposition of several results about Hopf modules, that are useful for
our purposes.  We then continue to prove the  main results in this
section, notably the  invertibility of the object of integrals, the
relationships of integrals with modular group-like elements and the
antipode, and, eventually, the generalized Radford formula.  In
Section~\ref{sec-exam} we shall develop the example of a braided Hopf
 algebra in a category of tangles, and explain its applications to
topological field theory.  This algebra can be functorially mapped to
representatives of another important class of braided Hopf algebras,
namely coends in abelian rigid braided categories.  Section
\ref{sec-cross} is devoted to connections of integration with
cross-products and transmutations.  In  Section~\ref{sec-appl-int} we
shall attempt to define external Hopf algebras, and discuss  duality
properties for Hopf bimodules. We present an explicit equivalence
between the categories of Hopf $H$\n-bimodules and Hopf
$H\pti$-bimodules given by the tensor product with $\Int H$.  From
these duality considerations we are able to obtain   two more proofs of
the generalized Radford formula.

\begin{acknowledgement}
One of us (V.L.) wishes to thank mathematicians at Institut de
Recherche Math\'ematique Avanc\'ee for warm hospitality during a visit,
which initiated the work on this paper. We also like to thank David
Radford for encouraging discussion.
\end{acknowledgement}

\section{Preliminaries}
\label{sec-prelim}
Throughout this paper the symbol $\C=(\C,\otimes,\1)$ denotes a strict
rigid monoidal category with braiding $\Psi$. For the convenience of
the reader we shall begin with the definition of strict rigid monoidal
categories.  In the following paragraphs the notions of categories with
split idempotents, braided Hopf algebras, and invertible objects are
introduced.

\subsection*{Rigid braided categories}
A (strict) {\em rigid category}, $\,\C\,$, is a (strict) monoidal
category, in which for every object, $\,X\in\C\,$, one can find a pair
of {\em  dual} objects, $\,X\pti\,$ and $\,\lpti X \in\C$,  as well as
morphisms of evaluation and coevaluation denoted as follows.
\begin{alignat*}3
\ev &: X\tens X\pti \to \1\  &=
\makebox[20mm][l]{
\raisebox{-3mm}[5mm][5mm]{
\unitlength 0.800mm
\linethickness{0.4pt}
\begin{picture}(16,14)
\put(10,7){\oval(10,14)[b]}
\put(4,5){\makebox(0,0)[rb]{$X$}}
\put(16,5){\makebox(0,0)[lb]{$X\pti$}}
\end{picture}
}} \ ,
 &\qquad \ev &: \lpti X\tens X \to \1\ &=
\makebox[20mm][l]{
\raisebox{-3mm}[5mm][5mm]{
\unitlength 0.800mm
\linethickness{0.4pt}
\begin{picture}(16,14)
\put(10,7){\oval(10,14)[b]}
\put(4,5){\makebox(0,0)[rb]{$\lpti X$}}
\put(16,5){\makebox(0,0)[lb]{$X$}}
\end{picture}
}} \ ,   \\
\coev &: \1\to X\pti\tens X\ &=
\makebox[20mm][l]{
\raisebox{-3mm}[5mm][5mm]{
\unitlength 0.80mm
\linethickness{0.4pt}
\begin{picture}(16,7)
\put(10,0){\oval(10,14)[t]}
\put(4,0){\makebox(0,0)[rb]{$X\pti$}}
\put(16,0){\makebox(0,0)[lb]{$X$}}
\end{picture}
}} \ ,
 &\qquad \coev &: \1\to X\tens\lpti X\ &=
\makebox[20mm][l]{
\raisebox{-3mm}[5mm][5mm]{
\unitlength 0.80mm
\linethickness{0.4pt}
\begin{picture}(16,7)
\put(10,0){\oval(10,14)[t]}
\put(4,0){\makebox(0,0)[rb]{$X$}}
\put(16,0){\makebox(0,0)[lb]{$\lpti X$}}
\end{picture}
}} \ \ .
\end{alignat*}
They are subject to the condition that
 the following  compositions between evaluations and coevaluations
are all  equal the identity morphism in $\,\End(X)\,$:
\begin{gather*}
X = X\tens \1 \xra{1\tens\coev} X\tens(X\pti\tens X) =
(X\tens X\pti)\tens X \xra{\ev\tens1} \1\tens X = X ,
\\
X = \1\tens X \xra{\coev\tens1} (X\tens\lpti X)\tens X
= X\tens(\lpti X\tens X) \xra{1\tens\ev} X\tens \1 = X ,
\\
X\pti = \1\tens X\pti \xra{\coev\tens1} (X\pti\tens X)\tens X\pti
= X\pti\tens(X\tens X\pti) \xra{1\tens\ev} X\pti\tens \1 = X\pti,
\\
\lpti X = \lpti X\tens \1 \xra{1\tens\coev}
\lpti X\tens(X\tens\lpti X) =
(\lpti X\tens X)\tens\lpti X \xra{\ev\tens1} \1\tens\lpti X = \lpti X
\end{gather*}

This data allows us to introduce the {\em transposes },
 $\,f^t:Y^\vee \to X^\vee\,$ and $\,{}^tf:\lpti Y \to \lpti X\,$,
 of a given morphism, $\,f:X\to Y\,$ in $\,\C\,$.

The  braiding $\,\{\Psi_{X,Y}:\,X\otimes Y\to Y\otimes X\}\,$
has to satisfy axioms, such as naturality and the
hexagonal equation. For details seek, e.g., \cite{JS:geom,Tur:q3}.

For a morphism, $\,f\,$, in a rigid braided monoidal category its trace
$\tr_8 f$ is defined using Figure~\ref{Fig-u}, the subscript being
motivated by the  resemblance with the digit 8.
The {\em  dimension} associated to this trace  is
$\,\dim_8(X):=\tr_8(\id_X)\;\in {\rm End}(\1)\,$.
We will also need the natural automorphism $\,u^0_{-2}: X \to X\,$
defined via the diagram in Figure~\ref{Fig-u}. The notation is borrowed
from \cite{Lyu:tan}.

\begin{figure}
\[
\tr_8 f:=\enspace
\begin{tangle}\hh\coev \\
       \O{f}\step\id \\
       \hxx \\
       \hh\ev\end{tangle}
\qquad\qquad
u^0_{-2}:=\enspace
\begin{tangle}\hh\Step\id\step\coev \\
       \hcoev\step\hx\step\id \\
       \id\step\hx\step\hev \\
       \hh\ev\step\id\end{tangle}
\]
\caption{}
\label{Fig-u}
\end{figure}

\subsection*{Categories with split idempotents}
M.~Karoubi introduced in \cite{karoubi1} a class of Banach categories,
which he calls {\em pseudo-abelian}. If we drop, irrelevant for us,
Banach and additive structures his definition reduces  to the
following:

  An idempotent, $\,e=e^2:X\to X\,$, in a category, $\,\mathcal D\,$,
 is said to be {\em split}, if there exists an object, $X_e\,$, and
morphisms, $\, i_e:X_e\to X\,$ and $\,p_e:X\to X_e\,$, such that
$\,e =i_e\circ p_e\,$ and $\,\id_{X_e}=p_e\circ i_e\,$.  If every
idempotent in $\,\mathcal{D}\,$ is split then we say that
$\,\mathcal{D}\,$ \emph{admits split idempotents}.

For a given category, $\,\C\,$, there exists a universal embedding,
$\,{\C} \xra i \widehat{\C}\,$, such that the category $\,\widehat{\C}\,$
admits split idempotents. Moreover, $\,\widehat{\C}\,$ can be chosen
universal in sense that for any category, $\,\mathcal{D}\,$, with split
idempotents every functor, $\,F:{\C}\to\mathcal{D}\,$, factors in the
form $\,F = \bigl(\C \xra i \widehat{\C} \xra G \mathcal{D}\bigr)\,$,
such that the functor $\,G\,$ is unique up to an isomorphism of functors.
The category $\widehat{\C}$ is called the \emph{Karoubi enveloping
category of} $\C$. According to Karoubi~\cite{karoubi1} it is realized
as the category with objects, $\,X_e=(X,e)\,$, where $\,X\,$ is an
object in $\,\C\,$ and $\,e:\,X\rightarrow X\,$ is an idempotent in
$\,\C\,$.  The morphisms in $\,\widehat{\C}\,$ are defined by
$\,{\widehat{\C}}(X_e,Y_f):= \{t\in{\C}(X,Y)\mid fte=t\}\,$. The
functor $\,i\,$, defined by $\,i(X)=X_{{\rm id}_X}\,$ and $i(f)=f\,$,
is a full embedding, that is we have
$\,\C(X,Y) = \widehat{\C}(i(X),i(Y))\,$.

It is noted in \cite{Lyu:mod} that if $\,\C\,$ is a (braided) monoidal
category then the category $\,\widehat\C\,$ can be equipped with a
(braided) monoidal structure:
\[
\1 := (\1,\id_{\1})\,,
\qquad
X_e\otimes Y_f :=(X\otimes Y)_{e\otimes f}\,,
\qquad
\Psi_{X_e,Y_f} := (f\otimes e)\circ\Psi_{X,Y}\,.
\]
In this case $\,i\,$ is a (braided) monoidal functor. Furthermore, if
the category $\,\C\,$ is rigid, so is $\,\widehat\C\,$, and the dual
objects of $\,(X,e)\,$ are $\,(X\pti,e^t)\,$ and $\,(\lpti X,{}^te)\,$.

        From now on we assume that our braided rigid category $\,\C\,$
admits split idempotents.

\subsection*{Hopf algebras}
Recall that a Hopf algebra, $\,H\in \C\,$,  \cite{Ma:bg} is an object,
$\,H\in\Obj\C\,$, together with  an associative multiplication,
$\,m:H\tens H \to H\,$, and an associative comultiplication,
$\,\Delta:H\to H\tens H\,$, obeying the bialgebra axiom
\begin{multline*}
\bigl(H\tens H \xra m H \xra\Delta H\tens H\bigr) \\
= \bigl(H\tens H \xra{\Delta\tens\Delta} H\tens H\tens H\tens H
\xra{H\tens\Psi\tens H}  H\tens H\tens H\tens H \xra{m\tens m}
H\tens H\bigr) .
\end{multline*}
Moreover, $\,H\,$ shall have  a unit, $\,\eta: \1\to H\,$, a counit,
$\,\varepsilon:H\to\1\,$, an antipode, $\,S:H\to H\,$, and an inverse
antipode, $\,S^{-1}:H\to H\,$, which shall satisfy axioms analogous to
the classical case.

A left (resp. right) module over an algebra, $\,H\,$, is an object,
$\,M\in\C\,$, equipped with an associative action,
$\,\mu_\ell: H\tens M\to M\,$ (resp. $\mu_r: M\tens H \to M$).  The
category of left (resp.  right) $\,H$\n-modules will be denoted by
$\,{}_H\C\,$ (resp.  $\,\C_H\,$). The morphisms in $\,{}_H\C\,$ are
those from $\,\C\,$ that are equivariant with respect to
$\,\mu_\ell\,$.  A left (resp. right) comodule over a coalgebra,
 $\,H\,$, is an object, $\,M\in\C\,$, equipped with an associative
coaction, $\,\Delta_\ell:M\to H\tens M\,$ (resp.
$\Delta_r:M\to M\tens H$).  The category of left (resp. right)
$H$\n-comodules will be denoted by $\,{}^H\C\,$ (resp.  $\,\C^H\,$).

\begin{figure}
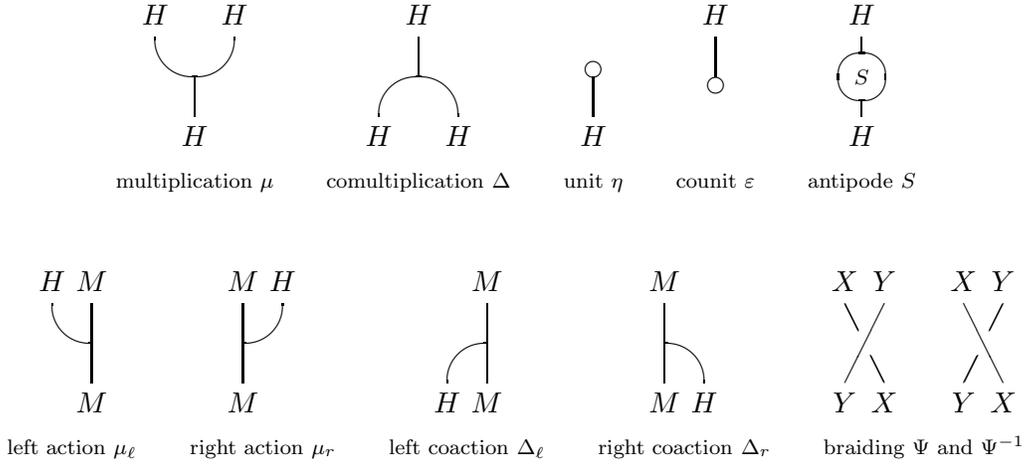

\[
\ba{ccccccccc}
\object H\Step\object H && \object H &&&& \object H && \object H \\
\begin{tangle}\cu\end{tangle}
&&
\begin{tangle}\cd\end{tangle}
&&
\begin{tangle}\unit\end{tangle}
&&
\begin{tangle}\counit\end{tangle}
&&
\begin{tangle}\S\end{tangle}
\\
\object H && \object H\Step \object H && \object H &&&& \object H \\
   \hbox{\scriptsize multiplication $\mu$}&&
   \hbox{\scriptsize comultiplication $\Delta$}&&
   \hbox{\scriptsize unit $\eta$}&&
   \hbox{\scriptsize counit $\varepsilon$}&&
   \hbox{\scriptsize antipode $S$}
\ea
\]
\[
\ba{ccccccccc}
\object H\step\object M && \object M\step\object H &&
\step\object M && \object M\step &&
\object X\step\object Y \Step \object X\step\object Y \\
\begin{tangle}\lu\end{tangle}
&&
\begin{tangle}\ru\end{tangle}
&&
\begin{tangle}\ld\end{tangle}
&&
\begin{tangle}\rd\end{tangle}
&&
\begin{tangle}\hxx\Step\hx\end{tangle}
\\
\step\object M && \object M\step &&
\object H\step\object M && \object M\step\object H &&
\object Y\step\object X \Step \object Y\step\object X \\
   \hbox{\scriptsize left action $\mu_\ell$}&&
   \hbox{\scriptsize right action $\mu_r$}&&
   \hbox{\scriptsize left coaction $\Delta_\ell$}&&
   \hbox{\scriptsize right coaction $\Delta_r$}&&
   \hbox{\scriptsize braiding $\Psi$ and $\Psi^{-1}$}
\ea
\]
\caption{Graphical notations}
\label{Fig-notations}
\end{figure}

If $\,({\C},\otimes,\1,\Psi)\,$ is a braided monoidal category we shall
denote by $\,\overline{\C}=({\C},\otimes,\1, \overline\Psi)\,$ the same
monoidal category with the mirror-reversed braiding
$\,\overline{\Psi}_{X,Y}:={\Psi_{Y,X}}^{-1}\,$. For a Hopf algebra,
$\,H\,$, in $\,\C\,$ we further denote by $\,H^\op\,$ (resp.
$\,H_\op\,$) the same coalgebra (resp. algebra) with opposite
multiplication $\,\mu^\op\,$ (resp.  opposite comultiplication
$\,\Delta^\op\,$) defined as follows:
\begin{equation}\label{opp}
\mu^\op:=\mu\circ\Psi_{H\,H}{}^{-1}\qquad\left({\rm resp.}\qquad
\Delta^\op:= \Psi_{H\,H}{}^{-1}\circ\Delta\right)\,.
\end{equation}
It is easy to see that $\,H^\op\,$ and $\,H_\op\,$ are Hopf algebras in
$\,\overline \C\,$ with antipode $\,S^{-1}\,$. We will always consider
$\,H^\op\,$ and $\,H_\op\,$ as objects of the category
$\,\overline\C\,$.  In what follows we often use graphical notations
for morphisms in monoidal categories, see
\cite{Bes:cross,JS:geom,Lyu:tan,Ma:introm,Tur:q3}.  The graphics and
notations  for (co\n-)multiplication, (co\n-)unit, antipode, left and
right (co\n-)action, and braiding are given in
Figure~\ref{Fig-notations}, where $\,H\,$ is a Hopf algebra and $\,M\,$
is an $\,H$\n-module ($\,H\,$\n-comodule).

\subsection*{Invertible objects}
\label{Append-invert}
An object, $\,K\,$, of a braided, monoidal category, $\,\C\,$, is
 called {\em invertible} if there exists an object,
$\,K^{-1}\in\Obj(\C)\,$, such that $\,K\otimes K^{-1}\simeq\1\,$ (and,
 hence $\,K^{-1}\otimes K\simeq\1\,$).  Properties of invertible
objects in a monoidal category are summarized in the following lemmas.
Most of the proofs are straightforward and left to the reader.

\begin{lemma}
\label{lemma-invert-obj}
Let $\,K\in\Obj(\C)\,$ be invertible. Then the following hold:
\begin{itemize}
\item
The morphisms
$\,\ev_K\,$ and $\,\coev_K\,$ are invertible;
\item
The map $\,\End(\1)\longrightarrow\End(K)\,:\;
c\mapsto c\otimes\id_K\,$
is an isomorphism of commutative monoids;
\item
$\dim_8(K)=\dim_8(K^{-1})\,$ is invertible and
\[
\Psi_{K,K}\;=\;\left(\dim_8K\right)^{-1}\cdot\id_{(K\otimes K)}\,.
\]
\end{itemize}
\end{lemma}

Only the last claim is not entirely obvious. The dimension
$\,\dim_8K\in\End(\1)\,$ is an isomorphism, because it is a composition
of isomorphisms, $\,\coev_K\,$, $\,\Psi_{K,K}\,$, and $\,\ev_K\,$. We
know that $\,\Psi_{K,K}=\psipsi\cdot\id_{(K\otimes K)}\,$ for some
$\,\psipsi\in\End(\1)\,$.  Finally
$\,\psipsi=\left(\dim_8K\right)^{-1}\,$ is implied by the diagrammatic
calculation below:

\[
 \begin{tangle}\hcoev \\
        \hxx \\
        \hev\end{tangle}
\enspace = \enspace
 \begin{tangle}\hh\hcoev \\
        \hxx\step\hcoev \\
        \id\step\hxx\step\id \\
        \hev\step\hxx \\
        \hh\Step\hev\end{tangle}
\enspace = \enspace \psipsi\
 \begin{tangle}\hh\hcoev \\
        \hxx\step\hcoev \\
        \id\step\id\step\id\step\id \\
        \hev\step\hxx \\
        \hh\Step\hev\end{tangle}
\quad.
\]

For a Hopf algebra, $\,H\,$, a morphism, $\,a:\1\to H\,$, is called a
{\em group-like element}, if $\,\Delta\circ a=a\otimes a\,$. A
morphism, $\,\alpha:H\to\1\,$, is called a {\em multiplicative
functional} if $\,\alpha\circ\mu=\alpha\otimes\alpha\,$.  For such
morphisms we set $\,a^{-1}:=S\circ a=S^{-1}\circ a\,$ and
$\,\alpha^{-1}:=\alpha\circ S=\alpha\circ S^{-1}\,$.

Invertible objects can be interpreted as the categorical analogues of
one-dimensional modules and  one-dimensional comodules of a linear Hopf
algebra.  The following two lemmas express properties that are obvious
for one-dimensional (co)modules  in the framework of a braided,
monoidal category. The proofs follow easily using the identities in
Lemma~\ref{lemma-invert-obj}.

\begin{lemma}
Suppose $\,K\,$ is an invertible object in a rigid, monoidal category,
$\,\C\,$.  Then any (co)module structure on $\,K\,$ over an
(co)algebra, $\,A\,$, in $\,\cal C\,$ is given by a multiplicative
functional  $\,\alpha\,$ (group-like element $\,a\,$):
\[
\mu_r=\id_K\otimes\alpha\,,\qquad
\Delta_r=\id_K\otimes a\,.
\]
\end{lemma}

\begin{lemma}
Let $\,K\,$ be an invertible object in a rigid, monoidal category,
$\,\C\,$, and $X,Y\in\Obj(\C)$. Then the map
$\,\Hom(X,Y) \to \Hom(X\tens K, Y\tens K)\,:\;f\mapsto f\tens\id_K\,$
is bijective.
\end{lemma}

The previous lemma, applied to the isomorphism
$\,g=\Psi_{K,X}\Psi_{X,K} : X\tens K \to X\tens K\,$ for an invertible
object $\,K\,$, implies that there exists a natural isomorphisms
$\,\monodromy{K}{X}:X\to X\,$ such that:
\[
\Psi_{K,X}\Psi_{X,K}=\monodromy{K}{X}\otimes\id_K
\qquad(\text{or, equivalently,}\quad
\Psi_{X,K}\Psi_{K,X}=\id_K\otimes\monodromy{K}{X}
)\,,
\]
We shall call $\,\monodromy{K}{X}\,$  the \emph{ monodromy}.  It may be
thought of as acting on $\,\cal C\,$ or Hopf algebras therein in the
following way:

\begin{lemma}\label{lem-monodromy}
The set of isomorphisms $\bigl\{\monodromy{K}{X}\bigr\}_X$ constitutes
an automorphism of the monoidal identity functor
\[
\monodromy{K}{\ }\,:\;(\id_{\C},\id_\tens)\longrightarrow
(\id_{\C},\id_\tens) \;:\; (\C,\tens) \to (\C,\tens)\qquad.
\]
Specifically, this means that
$\,I\circ \monodromy{K}{X}=\monodromy{K}{X}\circ I\,$ for any
$\,I\in {\rm Hom}(X,Y)\,$,  and, further, that
$\monodromy{K}{X\tens Y} = \monodromy{K}X\tens\monodromy{K}Y
: X\tens Y \to X\tens Y$.
\end{lemma}

\begin{corollary}\label{monodr-corol}
For any invertible object, $\,K\,$, and any Hopf algebra, $\,H\in\C\,$,
the morphism $\,\monodromy{K}H : H \to H\,$ is a Hopf algebra
automorphism.
\end{corollary}

\bigskip

\section{Integrals and the generalized Radford formula}
\label{sec-basic}
This section contains the detailed statement  of the main result of our
paper, namely the braided, categorical version of Radford's formula, as
well as the precise definitions of left and right integrals in a
category and of the respective (co)moduli entering this formula. We
shall defer all proofs to Section~4.  In the first paragraph on
integrals we will introduce a set of canonical projections,
$\,\PPi^\bullet_\bullet\,\in\,{\rm End}_{\cal C}(H)\,$, which will
guarantee the existence of integrals for a braided Hopf algebra,
$\,H\,$, in a category with split idempotents. We shall also illustrate
the generalizations in Radford's formula in the example of the category
of $\Zn$-graded vector spaces equipped with a braiding given by
Heisenberg-type relations.  Here the moduli will turn out to be units,
whereas the additional elements $\,u^0_{-2}\,$ and
$\,\monodromy{(\Int H)}H\,$ are going to be non-trivial.

\subsection*{Integrals for Hopf algebras}
The following definition  directly generalizes the classical one by
rewriting the defining formula into a relation of operations given by
morphisms.

\begin{definition}
Let $\,H\,$ be a bialgebra in a braided, monoidal category,
$\,\C\,$. A \emph{left $\,X$-valued integral on}
$H$ is a morphism $\,f: H\to X\,$ such that
\[ \bigl(H \xra\Delta H\tens H \xra{\id_H\tens f} H\tens X\bigr) =
\bigl(H \xra f X \xra{\eta_H\otimes\id_X} H\tens X\bigr) .\]
Right $X$\n-valued integrals on $H$ are defined analogously.

A \emph{left $X$\n-based integral in} $H$ is a morphism,
$\,f: X\to H\,$, such that
\[ \bigl(H\tens X \xra{\id_H\tens f} H\tens H \xra{\mu_H} H\bigr) =
\bigl(H\tens X \xra{\varepsilon_H\tens\id_X} X \xra f H\bigr) .\]
Right $\,X$\n-based integrals in $\,H\,$ are defined similarly.
\end{definition}

In other words, a left (resp. right) $X$\n-valued integral is a
homomorphism from the left (resp. right) regular $\,H$\n-comodule to
the trivial $\,H$\n-comodule $\,X\,$. A left (resp. right)
$\,X$\n-based integral is a homomorphism from the trivial $H$\n-module
$\,X\,$ to the left (resp.  right) regular $\,H$\n-module.

We shall consider only bialgebras \emph{with invertible antipode} and
call them Hopf algebras.

The next proposition not only asserts that a Hopf algebra in a category
with split idempotents always admits integrals, but also that the
objects of the integrals are invertible, and can be chosen the same for
all cases.

\begin{proposition}\label{pro-exist-int}
Assume $\,H\,$ is a Hopf algebra  with an invertible antipode, $\,S\,$,
in a braided, monoidal category, $\,\C\,$, with split idempotents.
Then there exist an invertible object $\Int\, H$ in $\,\C\,$, for which
the following hold simultaneously:

\begin{enumerate}
\item There exist left and right $\,\Int\, H$-valued integrals,
$\,\intl H\,$ and $\intr H :H\to \Int\,H\,$, such that any left (resp.
right) $\,X$\n-valued integral on $\,H\,$ admits a unique factorization
of the form $\,H \xra{\intl H} \Int H \xra g X\,$ (resp.
$\,H \xra{\intr H} \Int H \xra g X\,$).

\item There exist left and right $\,\Int\, H$-based integrals,
$\,\cointl H\,$ and $\,\cointr H :\Int H \to H\,$, such that any left
(resp.  right) $\,X$\n-based integral in $\,H\,$ admits a unique
factorization of the form $\,X \xra h \Int H \xra{\cointl H} H\,$
(resp.  $\,X \xra h \Int H \xra{\cointr H} H\,$).
\end{enumerate}
\end{proposition}

For proofs of this and the following results, see  Section~4.

It follows immediately from the universality of integrals that the
object $\,\Int_H\,$ is defined uniquely up to a unique isomorphism.

\begin{figure}
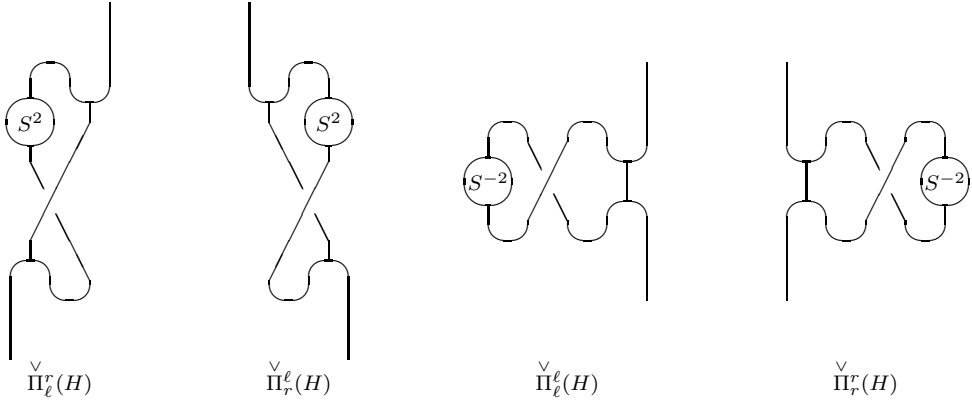

\[
\ba{ccccccc}
\begin{tangle}\hstep\hcoev\step\id \\
       \hstep\O{S^2}\step\hddcu \\
       \hstep\hxx \\
       \hh\cd\hstep\d \\
       \id\step\hev\end{tangle}
&\qquad\quad&
\begin{tangle}\id\step\hcoev \\
       \hdcu\step\O{S^2} \\
       \step\hxx \\
       \hh\hstep\dd\hstep\cd \\
       \hstep\hev\step\id\end{tangle}
&\qquad\quad&
\begin{tangle}\hcoev\step\hcoev\step\id \\
       \O{S^{-2}}\step\hxx\step\cucd \\
       \hev\step\hev\step\id\end{tangle}
&\qquad\quad&
\begin{tangle}\id\step\hcoev\step\hcoev \\
       \cucd\step\hxx\step\O{S^{-2}} \\
       \id\step\hev\step\hev\end{tangle}\\
\hbox{\scriptsize $\PPi^r_\ell(H)$}&&
\hbox{\scriptsize $\PPi^\ell_r(H)$}&&
\hbox{\scriptsize $\PPi^\ell_\ell(H)$}&&
\hbox{\scriptsize $\PPi^r_r(H)$}
\ea
\]
\caption{Four projections in $H$}
\label{Fig-4-Pi}
\end{figure}

Notice that $\,\intl H\,$ is a coequalizer for a pair
of left $\,{}^\vee H$-actions on $\,H\,$, given by the formulae:
\[
\mu' = (\ev\otimes\id_H)\circ(\id_{({}^\vee H)}\otimes\Delta_H)=
\ba{c}
\object{{}^\vee H}\step\hstep\object{H}\hstep\\
\begin{tangle}\hh\id\step\cd \\
       \hh\hev\step\id\end{tangle}
\\
\Step\object{H}
\ea\,,\enspace\quad
\mu'' = \varepsilon_{({}^\vee H)}\otimes\id_H=
\ba{c}
\object{{}^\vee H}\step\object{H}\\
\begin{tangle}\counit\step\id\end{tangle}\\
\step\object{H}
\ea
\]
Concretely, this means that $\,\intl H\,$ is a universal morphism for
which
\[
\bigl(\lpti H\tens H \xra{\mu'} H \xra{\intl H} \Int_H\bigr) =
\bigl(\lpti H\tens H \xra{\mu''} H \xra{\intl H} \Int_H\bigr).
\]
Also the right integral $\intr H$ on $\,H\,$ can be identified as a
 coequalizer for a similar pair of actions.  The integrals in $\,H\,$
are equalizers in a dual fashion. In the case of abelian categories
this property served as the definition of integrals, see
\cite{Lyu:mod}. In order to prove the existence of integrals in the
case of a category with split idempotents
(Proposition~\ref{pro-exist-int}), we construct the following four
idempotents (the next lemma will be derived as a special case of
Corollary~\ref{idemp-corollary}):

\begin{lemma}\label{Prop-4-idems}
The endomorphisms
$\PPi^r_\ell(H), \PPi^\ell_r(H), \PPi^\ell_\ell(H), \PPi^r_r(H) :H \to H$
given in Figure~\ref{Fig-4-Pi} are idempotents in $\End_{\C} H$.
\end{lemma}

The morphisms which split these idempotents are precisely the integrals
we are looking for. Some of their properties are described in the
following theorem.

\begin{theorem}\label{thm-intint=Pi}
The morphisms
$\intl H\circ\cointl H\,,\enspace
 \intl H\circ\cointr H\,,\enspace
 \intr H\circ\cointl H\,,\enspace
 \intr H\circ\cointr H\,\in\End(\Int\,H)$
are invertible, and, hence, by Lemma~\ref{lemma-invert-obj}
\begin{eqnarray*}
\intl H\circ\cointl H=\cll\cdot\id_{(\Int\,H)}\,,&\qquad&
\intl H\circ\cointr H=\clr\cdot\id_{(\Int\,H)}\,,\\
\intr H\circ\cointl H=\crl\cdot\id_{(\Int\,H)}\,,&\qquad&
\intr H\circ\cointr H=\crr\cdot\id_{(\Int\,H)}
\end{eqnarray*}
for certain invertible
$\cll,\,\clr,\,\crl,\,\crr\,\in\End(\1)$.

The idempotents from Lemma~\ref{Prop-4-idems} are $H$-valued and
$H$-based integrals. They are split by $\Int H$ as follows:
\begin{align*}
(\crl)^{-1}\cdot\cointl H\circ\intr H=\PPi^r_\ell(H) &\quad&
\vcenter{\hbox{is a right $H$-valued integral on $H$}
\hbox{and a left $H$-based integral in $H$,}} \\[2mm]
(\clr)^{-1}\cdot\cointr H\circ\intl H=\PPi^\ell_r(H) &\quad&
\vcenter{\hbox{is a left $H$-valued integral on $H$}
\hbox{and a right $H$-based integral in $H$,}} \\[2mm]
(\cll)^{-1}\cdot\cointl H\circ\intl H=\PPi^\ell_\ell(H) &\quad&
\vcenter{\hbox{is a left $H$-valued integral on $H$}
\hbox{and a left $H$-based integral in $H$,}} \\[2mm]
(\crr)^{-1}\cdot\cointr H\circ\intr H=\PPi^r_r(H) &\quad&
\vcenter{\hbox{is a right $H$-valued integral on $H$}
\hbox{and a right $H$-based integral in $H$.}}
\end{align*}
For any of the idempotents $\PPi^\bullet_\bullet(H)$ we have
\[
\tr_8\left(\PPi^\bullet_\bullet(H)\right)=\dim_8(\Int\,H)\,.
\]
\end{theorem}

The following lemma is a useful tool for recognizing the object
$\,\Int H\,$ of the integrals as the object of non-degenerate pairings
of a special form.

\begin{lemma}[\cite{KerLyu}]
Suppose there is an invertible object, $\,K\,$, and a morphism,
$\,t:  H\to K\,$, such that the pairing
$\,\phi: H\tens H \xra m H \xra t K\,$ is side-invertible (the induced
morphisms $\,H \to K\tens\lpti H\,$ and $\,H\to H\pti\tens K\,$ are
invertible). Then $\,\Int H \simeq K\,$.
\end{lemma}

\subsection*{The generalized Radford formula}
In the following lemma  group-like elements, $\,a\,$ and $\,\alpha\,$,
are extracted, which characterize the, e.g., right integral as a
homomorphism with respect to a left regular action. These special
elements are essential ingredients in  the generalized Radford formula.

\begin{lemma}\label{prop-group-like}
There exists a unique group-like element, $\,a:\1\rightarrow H\,$, and a
unique multiplicative functional, $\,\alpha :H\to\1\,$, such that the
following identities hold:
\begin{alignat}2
\left(\intl H\otimes\id_{H}\right)\circ\Delta&=\intl H\otimes a
&\enspace :\enspace& H\to\Int H\otimes H\,, \label{a-def}
\\ \vspace*{.3cm}
\left(\id_{H}\otimes\intr H\right)\circ\Delta&=a^{-1}\otimes\intr H
&\enspace :\enspace& H\to H\otimes\Int H\,, \label{a-1-def}
\\
\mu\circ\left(\cointl H\otimes\id_{H}\right)&=\cointl H\otimes\alpha
&\enspace :\enspace& \Int H\otimes H\to H\,,
\\
\mu\circ\left(\id_{H}\otimes\cointr H\right)&=\alpha^{-1}\otimes\cointr H
&\enspace :\enspace& H\otimes\Int H\to H\,. \label{alpha-1-def}
\end{alignat}
I.e. integrals are (co)module morphisms between $\,H\,$ with the
regular structure and $\,\Int H\,$ with the structure determined by
the group-like element $a$ (multiplicative functional $\alpha$).
\end{lemma}

\noindent
We shall use the following notations for morphisms in $\,\End(H)\,$:
\[
\ad_a:=\mu^{(3)}\circ(a\otimes\id_H\otimes a^{-1})\,,\qquad
\ad^\alpha:=(\alpha\otimes\id_H\otimes\alpha^{-1})\circ\Delta^{(3)}\,,
\]

where
\[
\mu^{(3)}:=\mu(\mu\otimes\id_H)=\mu(\id_H\otimes\mu)\,,\qquad
\Delta^{(3)}:=(\Delta\otimes\id_H)\Delta=(\id_H\otimes\Delta)\Delta\,.
\]

Given this the precise form of our main result is as follows:

\begin{theorem}[The generalized Radford formula]
\label{thm-Radford-formula}
\[ S^4\circ u^0_{-2} = (\ad^\alpha)^{-1}\circ \ad_a^{-1} \circ
\monodromy{(\Int H)}H\quad\in\Aut(H)\,, \]
where the monodromy action $\,\monodromy{K}{}\,$ of an invertible
object, $\,K\,$, is defined in Lemma~\ref{lem-monodromy}, and the
morphism $\,u^0_{-2}\,$ is defined on Figure~\ref{Fig-u}.
\end{theorem}

\begin{remark}
The components in the generalized Radford formula, namely,
$\,S^4\circ u^0_{-2}\,$, $\,\ad^\alpha\,$, $\,\ad_a\,$, and
$\,\monodromy{(\Int H)}H\,$, are Hopf algebra automorphisms, commuting
with each other. See also Corollary~\ref{monodr-corol}.
\end{remark}

The first proof of Theorem~\ref{thm-Radford-formula} is given in
Section~\ref{sec-proofs}. It follows the scheme of the original proof
given by Radford~\cite{Radford:antipode} for usual Hopf algebras, and
adapted by Kuperberg~\cite{Kuperberg} to the case of symmetric rigid
monoidal categories. We use several presentations of the antipode and
the inverse antipode via integrals.

The second proof of Theorem~\ref{thm-Radford-formula}, given in
Section~\ref{sec-appl-int}, is based on the properties of Hopf
bimodules.  The generalized Radford formula follows from the
equivalence of two presentations of the Hopf $\,H\pti$\n-bimodule
structure on $\,X\tens_H H\pti\,$, where $\,X\,$ is a Hopf
$\,H$\n-bimodule.

In the case of usual Hopf algebras the object of the integrals is a
one-dimensional vector space. In the case of $\ZZ/2$\n-graded Hopf
algebras the object of the integrals can be odd or even. In braided
categories there are more possibilities as the following example shows.

\begin{example}[Non-trivial object of integrals]
\label{non-triv-int}
Let $\,\C = \Zn$-grad-Vect($\Bbbk$) be the category of $\,\Zn$-graded
$\,\Bbbk$\n-vector spaces with the braiding
$\,\psi : V\tens W \to W\tens V\,$, given by
$\,\psi(v_i\tens w_j) = q^{ij} w_j\tens v_i\,$ for homogeneous vectors,
$\,v_i\in V^i\,$ and $\,w_j\in W^j\,$.  Here $\,q\in k\,$ is a
primitive $n^{\text{th}}$ root of unit (a root of the
$\,n^{\text{th}}\,$ cyclotomic polynomial $\,\phi_n\,$).  Notice that
the category $\,\C\,$  is also ribbon with respect to the ribbon twist
$\,\nu: V \to V\,$, where $\,\nu(v_i) = q^i v_i\,$, for
$\,v_i\in V^i\,$.

Denote by $\,X = X^1\,$ a selected one-dimensional vector space
concentrated  in degree 1. The tensor algebra $\,T(X)\,$, equipped with
a comultiplication induced by $\,\Delta (x) = x\tens 1 + 1\tens x\,$
for a basis vector, $\,x\in X^1\,$, constitutes a {\em braided} Hopf
algebra in the braided category $\,\C\,$.

Applying the bialgebra axiom, which contains the braiding $\,\psi\,$,
we find iteratively the comultiplication and the antipode as given
next:
\begin{gather*}
\Delta(x^m) = \sum_{k=0}^m \binom{m}k_q x^k\tens x^{m-k} , \\
S(x^m) = (-1)^m q^{m(m-1)/2} x^m .
\end{gather*}
For any $\,0<k<n\,$ the $q$-binomial coefficient $\,\tbinom nk_q\,$ is
divisible by $\,\phi_n(q)\,$. Therefore,
$\,\Delta x^n = x^n\tens 1 + 1\tens x^n\,$, and the span of $\,x^n\,$
is a coideal. Hence, $\,H = T(X)/(X^{\tens n})\,$ is a rigid Hopf
algebra in $\,\C\,$.

Using the definition of integrals one easily verifies that
$\,\intr H = \intl H : H \to X^{\tens n-1}\;:
\;x^k \mapsto \delta^k_{n-1} x^{\tens n-1}$,
is a two-sided $\,X^{\tens n-1}$-valued integral on $\,H\,$, and
$\,\cointr H = \cointl H : X^{\tens n-1} \to H\;:\;
x^{\tens n-1} \mapsto x^{n-1}$, is a right and left
$\,X^{\tens n-1}$-based integral in $\,H\,$. Since $\,X^{\tens n-1}\,$ is
one-dimensional, the above integrals are universal and
$\,\Int H = X^{\tens n-1} \simeq X\pti\,$. In particular,
$\,(\Int H)^{\tens n} \simeq \1\,$. One concludes that all idempotents
$\,\PPi^\bullet_\bullet\,$ coincide and are given by
$\,x^k \mapsto \delta^k_{n-1} x^{n-1}\,$.

The element $\,a\,$ is the unit, the functional $\,\alpha\,$ is the
counit.  The monodromy $\,\monodromy{(\Int H)}V\,$ determined by
$\,\Int H\,$ on $\,V\,$ is defined on homogeneous vectors by
$\,\monodromy{(\Int H)}V(v_k) = q^{-2k} v_k\,$. Since
$\,(u^0_{-2})_{V^k} = \nu_{V^k}^{-2} = q^{-2k^2}\,$,  the map on the
right-hand side of the generalized Radford formula turns out to be
$\,\monodromy{(\Int H)}H (x^k) = q^{-2k} x^k\,$, and on the left-hand
side the same map results from
$\,(S^4 \circ u^0_{-2})(x^k) = q^{2k(k-1)} \cdot q^{-2k^2} x^k\,$.
\end{example}

\section{Proofs of the main results}\label{sec-proofs}

The proof of Theorem~\ref{thm-intint=Pi} we will follow in this section
relies heavily on the theory of Hopf modules, which will be discussed
in the first paragraph. The projections   $\,\PPi^\bullet_\bullet\,$
are explicitly constructed, and their images identified as integrals.
In the following paragraph we prove the invertibility of the universal
object of the integrals, permitting the  conclusion of the proof of
Theorem~\ref{thm-intint=Pi}. In subsequent parts of this section we
discuss the ingredients of  the generalized Radford formula, namely the
special group-like elements, and the action of the antipode on the
integrals. With this preparation the proof of the formula for the
forth order  of the antipode is finally only a matter of combining
several identities in the right  way.

\subsection*{Hopf modules}

        From \cite{Lyu:mod} we know that the Structure Theorem of Hopf
modules \cite{Sweedler} also holds for Hopf algebras in braided,
monoidal categories, if there exist (co\n-)equalizers. In
\cite{BesDra:hopf,Bes:cross} it is shown that in order to prove the
Structure Theorem in a braided monoidal category, $\,\C\,$, it suffices
to assume that $\,\C\,$ admits split idempotents.

In a monoidal category, $\,\mathcal{D}\,$, which admits split
idempotents the subcategory of (co\n-)modules over a (co\n-)algebra
admits split idempotents, too.  Similar facts hold for the categories
of crossed modules and categories of Hopf (bi-)modules over a bialgebra
in a braided, monoidal category with split idempotents.

\begin{figure}
\[
\ba{c}
\ba{c}
\object{H}\step\object{X}\\
\begin{tangle}\lu \\
       \hh\ld\end{tangle} \\
\object{H}\step\object{X}
\ea
\, =\;
\begin{tangle}\hh\cd\step\hld \\
       \id\step\hxx\hstep\id \\
       \hh\cu\step\hlu\end{tangle}
\qquad
\ba{c}
\step\object{X}\step\object{H}\\
\begin{tangle}\step\ru \\
       \hh\ld\end{tangle}\\
\object{H}\step\object{X}\step
\ea
\!\!\!=\;
\begin{tangle}\ld\step\hcd \\
       \id\step\hxx\step\id \\
       \hcu\step\ru\end{tangle}
\qquad
\ba{c}
\object{H}\step\object{X}\step\\
\begin{tangle}\hh\lu \\
       \step\rd\end{tangle}\\
\step\object{X}\step\object{H}
\ea
\!\!\!=\;
\begin{tangle}\hcd\step\rd \\
       \id\step\hxx\step\id \\
       \lu\step\hcu\end{tangle}
\qquad
\ba{c}
\object{X}\step\object{H}\\
\begin{tangle}\ru \\
       \hh\rd\end{tangle} \\
\object{X}\step\object{H}
\ea
\!=\;
\begin{tangle}\hh\hrd\step\cd \\
       \id\hstep\hxx\step\id \\
       \hh\hru\step\cu\end{tangle}
\\
\hbox{\scriptsize a) Hopf module axioms}
\ea
\]
\[
\ba{c}
\Pi^\ell_\ell(X):=\quad
\begin{tangle}\ld \\
       \S\step\id \\
       \lu\end{tangle}
\qquad
\Pi_\ell^r(X):=\quad
\begin{tangle}\ld \\
       \SS\step\id \\
       \hx \\
       \ru\end{tangle}
\qquad
\Pi_r^\ell(X):=\quad
\begin{tangle}\rd \\
       \id\step\SS \\
       \hx \\
       \lu\end{tangle}
\qquad
\Pi^r_r(X):=\enspace
\begin{tangle}\rd \\
       \id\step\S \\
       \ru\end{tangle}
\\
    \hbox{\scriptsize b) Idempotents $\Pi^\bullet_\bullet(X)$}
\ea
\]
\[
\ba{ccccccc}
\begin{tangle}\hh\coev\step\id \\
       \O{S^2}\step\lu \\
       \xx \\
       \rd\step\id \\
       \hh\id\step\ev\end{tangle}
&\qquad\quad&
\begin{tangle}\hh\id\step\hcoev \\
       \ru\step\O{S^2} \\
       \xx \\
       \id\step\ld \\
       \hh\ev\step\id\end{tangle}
&\qquad\quad&
\begin{tangle}\hcoev\step\hcoev\step\id \\
       \O{S^{-2}}\step\hxx\step\kk \\
       \hev\step\hev\step\id\end{tangle}
&\qquad\quad&
\begin{tangle}\id\step\hcoev\step\hcoev \\
       \k\step\hxx\step\O{S^{-2}} \\
       \id\step\hev\step\hev\end{tangle}\\
\hbox{\scriptsize $\PPi^r_\ell(X)$}&&
\hbox{\scriptsize $\PPi^\ell_r(X)$}&&
\hbox{\scriptsize $\PPi^\ell_\ell(X)$}&&
\hbox{\scriptsize $\PPi^r_r(X)$}
\ea
\]
\caption{Hopf modules and projections}
\label{Fig-Hopf}
\end{figure}

\begin{definition}\

\begin{enumerate}
\item A {\em left Hopf module}, $\,X\,$, (resp. right-left Hopf module)
over a bialgebra, $\,H\,$, in $\,\C\,$ is a left (resp. right)
$\,H\,$\n-module and a left $\,H\,$\n-comodule, such that the action is
a comodule morphism.  Here the (co\n-)actions of $\,H\,$ on
$\,X\otimes H\,$ and $\,H\otimes X\,$ are determined by the braided,
 diagonal tensor (co\n-)action, given that $\,H\,$ is an
$\,H\,$-(co\n-)module via the regular (co\n-)action.  As illustrated in
Figure \ref{Fig-Hopf}a) the axiom implied by this condition can be
thought of as a ``polarized'' version of the bialgebra, which has to
hold  besides the comodule and module axioms for $\,X\,$.  We shall
denote by ${}_H^H{\C}$ (resp. ${}^H{\C}_H$) the subcategory of Hopf
modules, whose objects are the left (resp. right-left) Hopf modules,
and whose morphisms are the left-(resp.
right-)$H$-module-left-$H$-comodule homomorphisms.  The  definitions
for ${\C}_H^H$ and ${}_H{\C}^H$ are analogous.

\item A {\em two-fold  Hopf module},
$\,X =( X , \mu_\ell, \mu_r,\Delta_\ell)\,$,  is an object, which is a
$\,H$\n-bimodule in the category of left $H$\n-comodules, or in the
language of Hopf modules: $\,X\in\Obj({}_H^H\C)\,$ and
$\,X\in\Obj({}^H\C_H)\,$.  The two-fold Hopf modules together with the
$\,H$-bimodule-left-$H$-comodule morphisms form the category of
two-fold Hopf modules $\,{}_H^H\C_H\,$.  The remaining three types of
two-fold Hopf modules are defined in similar ways.

\item An object, $\,(X,\mu_r,\mu_\ell,\Delta_r,\Delta_\ell)\,$, is
called an $\,H$-Hopf bimodule if $\,(X,\mu_r,\mu_\ell)\,$ is a
$H$\n-bimodule, and $\,(X,\Delta_r,\Delta_\ell)\,$ is a
$\,H$\n-bicomodule in the category of $\,H$\n-bimodules, where the
regular (co\n-)action on $H$ and the diagonal (co\n-)action on tensor
products of modules are used.  Hopf bimodules together with the
$H$-bimodule-$H$-bicomodule morphisms form the category which will be
denoted by ${}_H^H\C^H_H$.
\end{enumerate}
\end{definition}

It follows directly from the definition that any bialgebra, $\,H\,$, is
a Hopf bimodule over itself with the regular actions and the regular
coactions. Next we give the construction of the morphism
$\,\Pi^\ell_\ell(X):X\to X\,$. Idempotency follows from straightforward
application of the Hopf algebra and Hopf module axioms as depicted in
Figure~\ref{Fig-Hopf}. In the same way the object, through which
$\,\Pi^\ell_\ell(X):X\to X\,$ factors, is identified as a
(co)invariance.

\begin{lemma}\label{idem-pi}
Let $\,H\,$ be a Hopf algebra in $\,\C\,$ and $\,(X,\mu_l,\Delta_l)\,$
be a left $\,H$-Hopf module.

Then the following hold:
\begin{itemize}
\item The morphism $\Pi^\ell_\ell(X):X\to X$ defined via
\begin{equation}
\label{Pi}
\Pi^\ell_\ell(X):=\mu_l\circ(S\otimes\id_X)\circ\Delta_l
\end{equation}
(see also Figure \ref{Fig-Hopf}b))
is an idempotent in ${\rm End}_{\C}( X)$.

\item Let $X\xrightarrow{{}_Xp}{}_HX\xrightarrow{{}_Xi}X$ be the
morphisms, which split the idempotent $\Pi^\ell_\ell(X)$, i.e.
${}_Xi\circ{}_Xp=\Pi^\ell_\ell(X)$ and
${}_Xp\circ {}_Xi={\rm id}_{{}_HX}$.  Then
\begin{equation}
\label{equalizer}
{}_HX \xrightarrow{{}_X i}  X
\genfrac{}{}{0pt}{2}{\xrightarrow[\hphantom{\eta\otimes\id_X}]{\Delta_l}}
                   {\xrightarrow[\eta\otimes\id_X]{}}
H\otimes X
\qquad\hbox{and}\qquad
 H\otimes X
\genfrac{}{}{0pt}{2}{\xrightarrow[\hphantom{\varepsilon\otimes\id_X}]{\mu_l}}
                   {\xrightarrow[\varepsilon\otimes\id_X]{}}
X\xrightarrow{{}_Xp}  {}_HX
\end{equation}
are equalizer and coequalizer respectively.  Hence $\,{}_HX\,$ is at
the same time an object of invariants and coinvariants of $X$.
\end{itemize}
\end{lemma}

The following lemma embodies  the categorical generalization of the
Fundamental Theorem of Hopf modules, which can be thought of as the
main reason for the structural rigidity of Hopf algebras.  In the
classical case it  essentially states that all Hopf modules are free.
The version below can be  derived by translating the classical proof
into the diagrammatic language, and generalizing from there to the
braided situation.

\begin{lemma}\label{equi-Hopfmod-lem}
The following functors define an equivalence of categories,
$\,\C
\genfrac{}{}{0pt}{2}{\xrightarrow{H\ltimes(\_)}}
                   {\xleftarrow[{}_H(\_)]{\hphantom{H\ltimes(\_)}}}
\sideset{^H_H}{}\C\,:
$
\begin{itemize}
\item $\,H\ltimes X\,$ is an object $H\otimes X$ with Hopf module
structure given by $\,\Delta_\ell:=\Delta\otimes\id_{X}\,$ and
 $\,\mu_\ell:=\mu\otimes\id_{X}\,$, called {\em standard Hopf module},
and $\,H\ltimes f:=\id_H\otimes f\,$ as morphism in $\,\C\,$;

\item $\,{}_HX\,$ splits the idempotent $\,\Pi^\ell_\ell(X)\,$ from the
previous lemma, and $\,{}_Hf:={}_Yp\circ f\circ {}_Xi\,$ for a Hopf
module morphism $\,f:X\to Y\,$;

\item For
$\,(X,\mu_\ell,\Delta_\ell)\in\Obj\left( \sideset{^H_H}{}\C \right)\,$
\[
H\otimes ({}_HX)
\genfrac{}{}{0pt}{2}{\xrightarrow
                      [\hphantom{(\id_H\otimes {}_Xp)\circ\Delta_\ell}]
                      {\mu_\ell\circ(\id_H\otimes{}_Xi)}}
                     {\xleftarrow[(\id_H\otimes {}_Xp)\circ\Delta_\ell]
                      {\hphantom{\mu_\ell\circ(\id_H\otimes{}_Xi)}}}
X
\]
are mutually inverse Hopf module morphisms.
\end{itemize}
\end{lemma}

Our strategy for extending the  result from the previous two lemmas for
$\,\Pi^\ell_\ell(X)\,$ to the other projections is to find appropriate
functors that turn left actions into right actions and map
$\,\Pi^\ell_\ell(X)\,$ to the respective $\,\Pi^\bullet_\bullet(X)\,$.
To this end suppose $\,X\,$ is an object in $\,\C\,$ equipped with some
or all of the $H$-module and $H$-comodule structures from above, namely
$\,\mu_\ell\,$, $\,\mu_r\,$, $\,\Delta_\ell\,$, and $\,\Delta_r\,$. We
shall introduce the following notations for the same object $\,X\,$
with modified (co)actions:

\begin{itemize}
\item $\,X^\Op\,$ denotes the underlying (bi)comodule with opposite
$\,H^\op\,$-actions\newline
\ \hspace*{2cm}$\,\mu^\Op_\ell:=\mu_r\circ\Psi^{-1}\,$,  \quad and\quad
$\,\mu^\Op_r:=\mu_\ell\circ\Psi^{-1}\,$.  \newline
$\,X_\Op\,$ denotes the underlying (bi)module with opposite
$\,H_\op$-coactions\newline
\ \hspace*{2cm}$\,\Delta^\Op_\ell:=\Psi^{-1}\circ\Delta_r\,$,
\quad and\quad $\,\Delta^\Op_r:=\Psi^{-1}\circ\Delta_\ell\,$;

\item $\,X^\S\,$ denotes the underlying (bi)comodule with
$\,H^\op$-actions\newline
\ \hspace*{2cm}$\,\mu^\S_\ell:=\mu_\ell\circ(S\otimes\id)\,$,
\quad and\quad $\,\mu^\S_r:=\mu_r\circ(\id\otimes S)\,$. \newline
$\,X_\S\,$ denotes the underlying (bi)module with
$\,H_\op$-coactions\newline
\ \hspace*{2cm}
$\,\Delta^\S_\ell:=(S\otimes\id)\circ\Delta_\ell\,$, \quad and\quad
$\,\Delta^\S_r:=(\id\otimes S)\circ\Delta_r\,$.

\item For $\,(X,\mu_\ell,\Delta_r)\,$ (resp.
$\,(X,\mu_r,\Delta_\ell)\,$\ ) we put
$\,X^\Vee:=(X,\mu^\Vee_r,\Delta^\Vee_\ell)\,$
(resp. $\,{}^\Vee X:=(X,{}^\Vee\mu_\ell,{}^\Vee\Delta_r)\,$\ ) \
to be the object equipped with the $\,H^\vee$- (resp. $\,{}^\vee H$-)
(co)actions as depicted below:
\begin{equation}
\label{Dual-act}
{}^\Vee\mu_\ell=
\begin{tangle}\hh\id\step\hld \\
       \hh\ev\hstep\id\end{tangle}\,,
\qquad
\mu^\Vee_r=
\begin{tangle}\hh\hrd\step\id \\
       \hh\id\hstep\ev\end{tangle}\,,
\qquad
\Delta^\Vee_\ell=
\begin{tangle}\hh\coev\hstep\id \\
       \hh\id\step\hlu\end{tangle}\,,
\qquad
{}^\Vee\Delta_r=
\begin{tangle}\hh\id\hstep\coev \\
       \hh\hru\step\id\end{tangle}
\end{equation}
\end{itemize}

Hence the so defined operations
$\,(\_)^\Op,\,(\_)_\Op,\,(\_)^\S,\,(\_)_\S,\,(\_)^\Vee\,$, and
$\,{}^\Vee(\_),\,$ identify, e.g. for $\,(\_)^\Op\,$, a module-object
$\,X\,\in\,{\rm Obj}(\C)\,$ for a Hopf algebra $\,H\,$ to be also a
comodule in $\,X\,\in\, {\rm Obj}(\overline\C)\,$ for $\,H^\Op\,$,
where $\,\overline\C\,$ is  the same monoidal category as before but
 with mirrored braiding.
It is easy to see that, similarly, an $\,H-$module
morphism in $\,\C\,$ is at the same time an $\,H^\Op-$comodule morphism
in $\,\overline\C\,$ with the modified coaction as above. Thus in this
 way  we obtain functors between the (co)module subcategories of
$\,\C\,$ as specified next. As functors of only the monoidal category
$\,\C\,$ they act like identity on objects and morphisms.


\begin{lemma}
\label{cat-equiv-lem}
The following functors define equivalences of categories of Hopf
(bi)modules:
\begin{eqnarray*}
(\_)^\Op:\sideset{^H_H}{^H_H}{\C}\to
\sideset{^{H^\op}_{H^\op}}{^{H^\op}_{H^\op}}{\mathop{\overline\C}}\,,
\quad&&\enspace
(\_)_\Op:\sideset{^H_H}{^H_H}{\C}\to
\sideset{^{H_\op}_{H_\op}}{^{H_\op}_{H_\op}}{\mathop{\overline\C}}\,,
\\
\bigl((\_)^\S\bigr)_\S:\sideset{^H_H}{^H_H}{\C}\to
\sideset{^{(H^\op)_\op}_{(H^\op)_\op}}{^{(H^\op)_\op}_{(H^\op)_\op}}%
{\C}\,,
\quad&&\enspace
\bigl((\_)_\S\bigr)^\S:\sideset{^H_H}{^H_H}{\C}\to
\sideset{^{(H_\op)^\op}_{(H_\op)^\op}}{^{(H_\op)^\op}_{(H_\op)^\op}}%
{\C}\,,
\\
\bigl((\_)^\Vee\bigr)^\S:\sideset{_H}{^H}{\C}\to
\sideset{^{H^{\vee\;\op}}}{_{H^{\vee\;\op}}}{\mathop{\overline\C}}\,,
\quad&&\enspace
\bigl((\_)^\S\bigr)^\Vee:\sideset{_H}{^H}{\C}\to
\sideset{^{H^{\op\;\vee}}}{_{H^{\op\;\vee}}}{\mathop{\overline\C}}\,,
\\
\bigl({}^\Vee(\_)\bigr)^\S:\sideset{^H}{_H}{\C}\to
\sideset{_{({}^\vee H)^\op}}{^{({}^\vee H)^\op}}{\mathop{\overline\C}}\,,
\quad&&\enspace
{}^\Vee\bigl((\_)^\S\bigr):\sideset{^H}{_H}{\C}\to
\sideset{_{{}^\vee(H^\op)}}{^{{}^\vee(H^\op)}}{\mathop{\overline\C}}\,,
\end{eqnarray*}
\end{lemma}

For any functor from the above list with the source category
$\,\sideset{^H_H}{^H_H}{\C}\,$ we keep the same notation for other
functors, which act by the same rule on categories of Hopf modules or
two-fold Hopf modules. We now find  the desired results for the
morphisms $\,\PPi^\bullet_\bullet(X)\,$:

\begin{corollary}
\label{idemp-corollary}\

For a Hopf bimodule, $\,X\,$,  the endomorphisms
$\,\Pi^\bullet_\bullet(X)\,$ and  $\,\PPi^\bullet_\bullet(X)\,$, as
given in Figure~\ref{Fig-Hopf}, are idempotents.

The objects that split idempotents $\Pi^\bullet_\bullet(X)$ (resp.
$\PPi^\bullet_\bullet(X)$) are all isomorphic.
A representative of the corresponding class of isomorphic objects
shall be denoted by $\,\Inv X\,$ (resp. by \ $\,\Int X\,$).
\end{corollary}

\begin{proof}
All assertions about the endomorphisms can be  reduced to the ones for
the idempotent from Lemma~\ref{idem-pi} using the identities below:
\begin{eqnarray*}
\Pi^r_\ell(X)=\Pi^\ell_\ell(X^\Op)\,,\quad&\quad
\Pi^\ell_r(X)=\Pi^\ell_\ell(X_\Op)\,,\\
\PPi^r_\ell(X)=\Pi^\ell_\ell(X^{\Vee\S\Op})\,,\quad&\quad
\PPi^\ell_r(X)=\Pi^\ell_\ell(({}^\Vee X)_{\S\Op})\,.
\end{eqnarray*}
Let $\,X\xra{p}Y\xra{i}X\,$ (resp.
$\,X\xra{p^\prime}Y^\prime\xra{i^\prime}X\,$) be morphisms, which split
the idempotent $\,\Pi^\ell_\ell(X)\,$ (resp. $\,\Pi^\ell_r(X)\,$).
Then by Lemma~\ref{idem-pi} both morphisms $\,i\,$ and $\,i^\prime\,$
are equalizers of the pair
$\,\Delta_\ell,\eta\otimes\id_X:X\to H\otimes X\,$.  Thus by the
universal property of equalizer we have $\,Y\simeq Y^\prime\,$.
\end{proof}

\begin{figure}
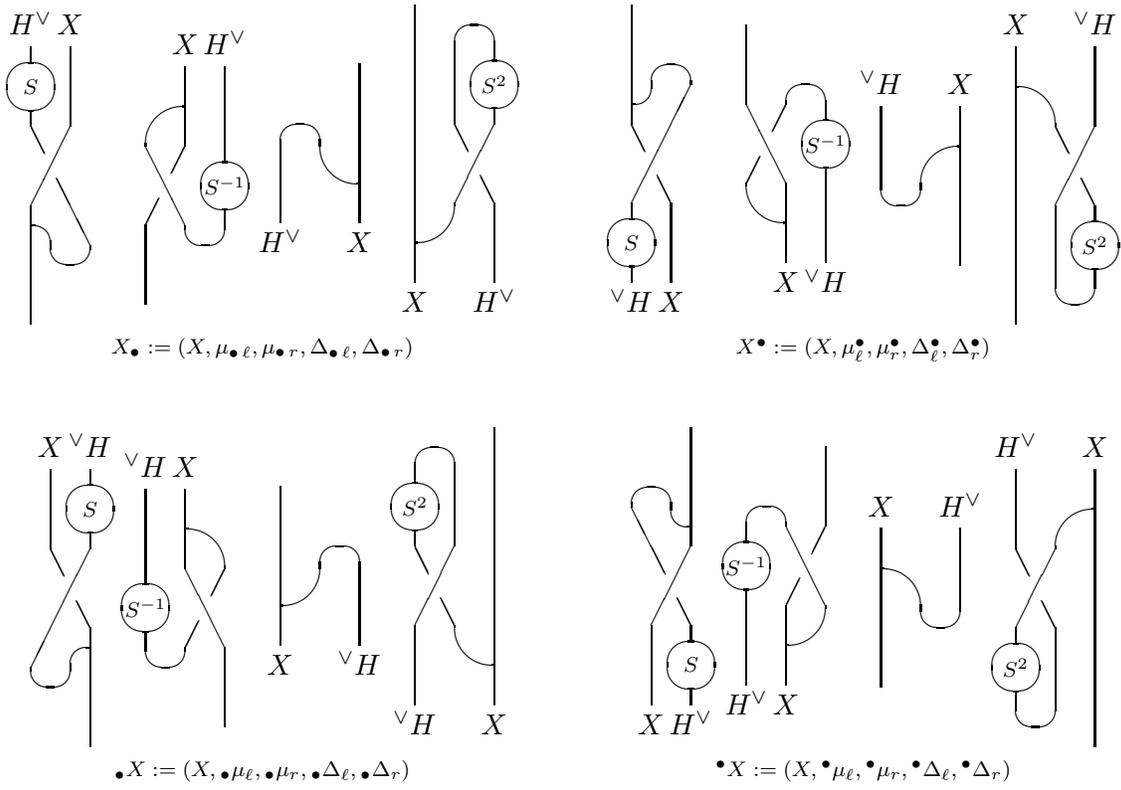

\[
\ba{ccc}
\ba{l}
\object{H^\vee}\step\object{X}\\
\begin{tangle}\S\step\id \\
       \hxx \\
       \hh\hrd\hstep\d \\
       \id\hstep\hev\end{tangle}
\ea
\quad
\ba{r}
\object{X}\step\object{H^\vee}\\
\begin{tangle}\ld\step\id \\
       \hx\step\SS \\
       \id\step\hev\end{tangle}
\ea
\quad
\ba{c}
\begin{tangle}\hcoev\step\id \\
       \id\step\lu\end{tangle}\\
\object{H^\vee}\Step\object{X}
\ea
\quad
\ba{r}
\begin{tangle}\hh\id\step\coev \\
       \id\step\id\step\O{S^2} \\
       \id\step\hxx \\
       \ru\step\id\end{tangle}\\
\object{X}\Step\object{H^\vee}
\ea
&\qquad&
\ba{l}
\begin{tangle}\id\hstep\hcoev \\
       \hh\hru\hstep\dd \\
       \hxx \\
       \S\step\id\end{tangle}\\
\object{{}^\vee H}\step\object{X}
\ea
\quad
\ba{r}
\begin{tangle}\id\step\hcoev \\
       \hx\step\SS \\
       \lu\step\id\end{tangle}\\
\object{X}\step\object{{}^\vee H}
\ea
\quad
\ba{c}
\object{{}^\vee H}\Step\object{X}\\
\begin{tangle}\id\step\ld \\
       \hev\step\id\end{tangle}
\ea
\quad
\ba{r}
\object{X}\Step\object{{}^\vee H}\\
\begin{tangle}\rd\step\id \\
       \id\step\hxx \\
       \id\step\id\step\O{S^2} \\
       \hh\id\step\ev\end{tangle}
\ea
\\
\hbox{\scriptsize $X_\bullet:=
(X,\mu_{\bullet\,\ell},\mu_{\bullet\,r},
\Delta_{\bullet\,\ell},\Delta_{\bullet\,r})$}&&
\hbox{\scriptsize $X^\bullet:=
(X,\mu^\bullet_\ell,\mu^\bullet_r,
\Delta^\bullet_\ell,\Delta^\bullet_r)$}
\ea
\]
\[
\ba{ccc}
\ba{r}
\object{X}\step\object{{}^\vee H}\\
\begin{tangle}\hstep\id\step\S \\
       \hstep\hxx \\
       \hh\dd\hstep\hld \\
       \hev\hstep\id\end{tangle}
\ea
\quad
\ba{l}
\object{{}^\vee H}\step\object{X}\\
\begin{tangle}\id\step\rd \\
       \SS\step\hx \\
       \hev\step\id\end{tangle}
\ea
\quad
\ba{c}
\begin{tangle}\id\step\hcoev \\
       \ru\step\id\end{tangle}\\
\object{X}\Step\object{{}^\vee H}
\ea
\quad
\ba{r}
\begin{tangle}\hh\coev\step\id \\
       \O{S^2}\step\id\step\id \\
       \hxx\step\id \\
       \id\step\lu\end{tangle}\\
\object{{}^\vee H}\Step\object{X}
\ea
&\qquad&
\ba{r}
\begin{tangle}\hcoev\hstep\id \\
       \hh\d\hstep\hlu \\
       \hstep\hxx \\
       \hstep\id\step\S\end{tangle}\\
\object{X}\step\object{H^\vee}
\ea
\quad
\ba{l}
\begin{tangle}\hcoev\step\id \\
       \SS\step\hx \\
       \id\step\ru\end{tangle}\\
\object{H^\vee}\step\object{X}
\ea
\quad
\ba{c}
\object{X}\Step\object{H^\vee}\\
\begin{tangle}\rd\step\id \\
       \id\step\hev\end{tangle}
\ea
\quad
\ba{r}
\object{H^\vee}\Step\object{X}\\
\begin{tangle}\id\step\ld \\
       \hxx\step\id \\
       \O{S^2}\step\id\step\id \\
       \hh\ev\step\id\end{tangle}
\ea
\\
\hbox{\scriptsize ${}_\bullet X:=
(X,{}_\bullet\mu_\ell,{}_\bullet\mu_r,
{}_\bullet\Delta_\ell,{}_\bullet\Delta_r)$}
&&
\hbox{\scriptsize ${}^\bullet X:=
(X,{}^\bullet\mu_\ell,{}^\bullet\mu_r,
{}^\bullet\Delta_\ell,{}^\bullet\Delta_r)$}
\ea
\]
\caption{Hopf module structures over the dual Hopf algebra}
\label{Fig-over-dual}
\end{figure}

Other functors can be obtained as compositions of the functors from
Lemma~\ref{cat-equiv-lem}. For a Hopf $\,H$-bimodule $\,X\,$ we
consider the underlying object equipped with (co)module structures over
the dual Hopf algebra $\,H^\vee\,$ or $\,{}^\vee H\,$ as shown in
Figure~\ref{Fig-over-dual}.

\begin{corollary}
\label{two-fold-Hopf-fun}
The following functors, whose actions on objects are as in
Figure~\ref{Fig-over-dual}, and which are identity on morphisms, define
the following equivalences of categories:
\begin{eqnarray*}
(\_)_\bullet:\sideset{^H_H}{^H}\C\rightarrow
\sideset{^{(H^\vee)}_{(H^\vee)}}{_{(H^\vee)}}\C\,,
\quad&\quad
(\_)_\bullet:\sideset{^H_H}{_H}\C\rightarrow
\sideset{^{(H^\vee)}_{(H^\vee)}}{^{(H^\vee)}}\C\,,
\\
(\_)^\bullet:\sideset{^H_H}{^H}\C\rightarrow
\sideset{^{{}^\vee H}_{{}^\vee H}}{_{{}^\vee H}}\C\,,
\quad&\quad
(\_)^\bullet:\sideset{^H_H}{_H}\C\rightarrow
\sideset{^{{}^\vee H}_{{}^\vee H}}{^{{}^\vee H}}\C\,,
\\
{}_\bullet(\_):\sideset{^H}{^H_H}\C\rightarrow
\sideset{_{({}^\vee H)}}{^{({}^\vee H)}_{({}^\vee H)}}\C\,,
\quad&\quad
{}_\bullet(\_):\sideset{_H}{^H_H}\C\rightarrow
\sideset{^{({}^\vee H)}}{^{({}^\vee H)}_{({}^\vee H)}}\C\,,
\\
{}^\bullet(\_):\sideset{_H}{^H_H}\C\rightarrow
\sideset{^{H^\vee}}{^{H^\vee}_{H^\vee}}\C\,,
\quad&\quad
{}^\bullet(\_):\sideset{^H}{^H_H}\C\rightarrow
\sideset{_{H^\vee}}{^{H^\vee}_{H^\vee}}\C\,,
\end{eqnarray*}
\end{corollary}

\begin{proof}
The functors $(\_)_\bullet$ between two-fold Hopf modules are ``glued''
from the following composite functors between Hopf modules:
\begin{eqnarray*}
&&\sideset{_H^H}{}\C
\xra{(\_)_\Op}
\sideset{_{(H_\op)}}{^{(H_\op)}}{\mathop{\overline\C}}
\xra{(\_)^{\Vee\Op}}
\sideset{^{(H^\vee)}}{_{(H^\vee)}}\C
\\
&&\sideset{_H}{^H}\C
\xra{(\_)^{\Vee\Op}}
\sideset{^{(H^\vee)^\op}}{_{(H^\vee)^\op}}{\mathop{\overline\C}}
\xra{(\_)_\Op}
\sideset{^{(H^\vee)}_{(H^\vee)}}{}\C\,,
\\
\sideset{^H}{_H}\C
\xra{\bigl((\_)^\S\bigr)_\S}
&&\sideset{^{(H^\op)_\op}}{_{(H^\op)_\op}}\C
\xra{\bigl((\_)_\Op\bigr)^\Op}
\sideset{_H}{^H}\C
\xra{\bigl((\_)^\Vee\bigr)_\S}
\sideset{^{(H^\vee)_\op}}{_{(H^\vee)_\op}}{\mathop{\overline\C}}
\xra{(\_)_\Op}
\sideset{}{^{H^\vee}_{H_\vee}}\C
\end{eqnarray*}
Other functors are reversed forms of $(\_)_\bullet$.
\end{proof}

Note that $\,X_\bullet\,$ is not a Hopf bimodule. The left-right Hopf
module axiom is not satisfied. In Section~\ref{sec-appl-int} we will
introduce a construction, which turns $\,X_\bullet\,$ into a Hopf
$H^\vee$-bimodule.

\subsection*{Invertibility of the object of integrals}

For the  Hopf bimodule $\,H\,$, equipped with the regular action and
coaction, we consider the underlying left Hopf $\,H^\vee$-module of
$\,H_\bullet\,$ (the explicit structures are shown in
Figure~\ref{Fig-proof-invert}a)). By Lemma~\ref{idem-pi} there exists
an isomorphism, $\,{\mathcal F}^H:H_\bullet\to H^\vee\otimes\Int H\,$,
into a standard Hopf $\,H^\vee$-module presented in
Figure~\ref{Fig-Fourier}. Similarly, we can construct Hopf module
isomorphisms corresponding to three other types of integrals.  The
inverse morphisms are given by the formulae
\begin{eqnarray*}
\crl({\mathcal F}^H)^{-1}=\overline{\mathcal F}^H&:=&
 {}_H{\mathcal F}\circ\Psi_{H^\vee,\Int H}\circ
 (S_{H^\vee}\otimes\id_{\Int H})
\\
\clr({}^H{\mathcal F})^{-1}={}^H\overline{\mathcal F}&:=&
 {\mathcal F}_H\circ\Psi_{\Int H,{}^\vee H}\circ
 (\id_{\Int H}\otimes S_{{}^\vee H})
\\
\clr({\mathcal F}_H)^{-1}=\overline{\mathcal F}_H&:=&
 (S_{{}^\vee H} \otimes \id_{\Int H})
\circ\Psi_{\Int H,{}^\vee H}\circ  {}^H{\mathcal F}
\\
\crl({}_H{\mathcal F})^{-1}={}_H\overline{\mathcal F}&:=&
 (\id_{\Int H}\otimes S_{H^\vee})
\circ\Psi_{H^\vee,\Int H}\circ  {\mathcal F}^H
\end{eqnarray*}
In the case of self-dual Hopf algebras and {\em bosonic} integrals,
meaning $\,\Int H\simeq\1\,$, these isomorphisms turn into Fourier
transforms \cite{Lyu:mod}.  The braided Fourier transform
$\,H\to H\otimes{}^\vee H\,$ defined in \cite{Cryss} (where the
$\,H$-valued integral is used) factorizes over our ${}^H\mathcal F$.
Given these isomorphisms we shall now prove
Proposition~\ref{pro-exist-int} and Theorem~\ref{thm-intint=Pi} :

\begin{figure}
\[
\ba{cccc}
\ba{c}
\Step\object H\\
\begin{tangle}\hh\coev\step\id \\
       \hh\id\step\cu \\
       \id\step[1.5]\Ointr H\end{tangle}\\
\object{H^\vee}\step[1.5]\object{\Int H}\hstep
\ea
&
\ba{c}
\object H\Step\\
\begin{tangle}\hh\id\step\coev \\
       \hh\cu\step\id \\
       \hstep\Ointl H\step[1.5]\id\end{tangle}\\
\hstep\object{\Int H}\step[1.5]\object{{}^\vee H}
\ea
&
\ba{c}
\object{{}^\vee H}\step[1.5]\object{\Int H}\hstep\\
\begin{tangle}\id\step[1.5]\Ocointr H \\
       \hh\id\step\cd \\
       \hh\ev\step\id\end{tangle}\\
\Step\object H
\ea
&
\ba{c}
\hstep\object{\Int H}\step[1.5]\object{H^\vee}\\
\begin{tangle}\hstep\Ocointl H\step[1.5]\id \\
       \hh\cd\step\id \\
       \hh\id\step\ev\end{tangle}\\
\object H\Step
\ea
\\
\hbox{\scriptsize ${\mathcal F}^H:H_\bullet\to H^\vee\otimes\Int H\;$}&
\hbox{\scriptsize
   ${}^H{\mathcal F}:{}_\bullet H\to\Int H \otimes{}^\vee H\;$}&
\hbox{\scriptsize ${\mathcal F}_H:{}^\vee H\otimes\Int H \to H^\bullet\;$}&
\hbox{\scriptsize ${}_H{\mathcal F}:\Int H\otimes H^\vee \to {}^\bullet H$}
\ea
\]
\caption{}
\label{Fig-Fourier}
\end{figure}

\begin{proof}[Proof of Proposition~\ref{pro-exist-int}]
In order to prove invertibility of $\,\Int H\,$, we consider the
Hopf module isomorphisms
$\,{\mathcal F}^H:H_\bullet\to H^\vee\otimes\Int H\,$
in $\,\sideset{^{H^\vee}_{H^\vee}}{}\C\,$ and
$\,\overline{\mathcal F}_{(H_\op)^\vee}:
((H_\op)^\vee)^\bullet\to (H_\op)\otimes\Int H^\vee\,$
in $\,\sideset{^{H_\op}_{H_\op}}{}{\mathop{\overline{\C}}}\,$.
The explicit Hopf module structures on $H_\bullet$ and
$(H_\op)^\bullet$ are shown in Figure \ref{Fig-proof-invert}a),b).  Let
us consider the functor
$\,({}^\Vee((\_)^\Op))_{\S\Op}: \sideset{^{H^\vee}_{H^\vee}}{}\C\to
      \sideset{^{H_\op}_{H_\op}}{}{\mathop{\overline{\C}}}\,$. The
$\,H_\op$-(co)module structures on $({}^\Vee((X)^\Op))_{\S\Op}$ are
shown on Figure~\ref{Fig-proof-invert}c).  This functor converts
$\,H_\bullet\,$ into the regular Hopf $\,H_\op$\n-module, and the
regular Hopf module $\,H^\vee$ into $((H_\op)^\vee)^\bullet$.  Hence
the composition
\begin{equation*}
H\,\,\xra{\,{\mathcal F}^H\,}\,\,H^\vee\otimes\Int H\,\,
\xra{\,\overline{\mathcal F}_{(H_\op)^\vee}\otimes\id_{(\Int H)}\,}\,\,
H\otimes\Int H^\vee\otimes\Int H
\end{equation*}
is an isomorphism of standard Hopf modules in
$\sideset{^{H_\op}_{H_\op}}{}{\mathop{\overline\C}}$.  From this we
conclude that
\begin{multline*}
(\varepsilon_H\otimes\id_{(\Int H^\vee\otimes\Int H)})\circ
(\overline{\mathcal F}_{(H_\op)^\vee}\otimes\id_{(\Int H)})\circ
{\mathcal F}^H\circ\eta_H \\
= (\intl{H^\vee}\otimes\intr{H})\circ\coev\,\,:\,\,\,\,
\1\,\to\,\Int H^\vee\otimes\Int H\,\,\,\,\,
\end{multline*}
is an isomorphism.
\end{proof}

\begin{proof}[Proof of Theorem \ref{thm-intint=Pi}]
Let $\,\intr H:\,H\to\Int H\,$ be a right integral on $\,H\,$, and
$\,\cointl H:\,\Int H\to H\,$ be a left integral in $\,H\,$.  By the
universal property of integrals there exist isomorphisms,
$\,f:\,K\to\Int H\,$ and $\,\enspace g:\Int H\to K\,$, such that
$\,f\intr H=fp\,$ and $\,\enspace\cointl H=ig\,f$. Here $\,i\,$ and
$\,p\,$ are as before the morphisms splitting for example the
idempotent $\,\PPi^r_\ell(H)\,$ with object $\,K\,$.  Then
$\,\intr H\circ\cointl H\,=\,fpig\,=\,fg\, = \,\crl\cdot\id_{\Int H}\,$
is an isomorphism, and we find
\[
\left(\crl\right)^{-1}\cdot\cointl H\circ\intr H=
\cointl H\circ\left(\intr H\circ\cointl H\right)^{-1}\circ\intr H=
ig(fg)^{-1}fp=ip=\PPi^r_\ell(H)
\]
We also compute,
\begin{eqnarray*}
\tr_8\left(\PPi^r_\ell(H)\right)
&=&\tr_8\left(
   \cointl H\circ\left(\intr H\circ\cointl H\right)^{-1}\circ\intr H
                                                               \right)
 \\
&=&\tr_8\left(
   \intr H\circ\cointl H\circ\left(\intr H\circ\cointl H\right)^{-1}
                                                               \right)
 \\
&=&\tr_8\left(\id_{(\Int H)}\right)=\dim_8(\Int H)\,.
\end{eqnarray*}
\end{proof}

\begin{figure}
\[
\ba{ccccc}
\ba{l}
\hstep\object{H^\vee}\step\object{H}\\
\begin{tangle}\hstep\S\step\id \\
       \hstep\hxx \\
       \hh\cd\hstep\d \\
       \id\step\hev\end{tangle}         \\
\object{H}
\ea
\qquad
\ba{l}
\Step\object{H}\\
\begin{tangle}\hcoev\step\id \\
       \id\step\hcu\end{tangle}     \\
\object{H^\vee}\step\hstep\object{H}\hstep
\ea
&\qquad&
\ba{l}
{\object{H}}
\step\hstep
{\object{{H^\vee}}}\\
\begin{tangle}\id\step\hcd \\
       \hev\step\id\end{tangle}\\
\Step\object{H^\vee}
\ea\qquad
\ba{l}
\object{H^\vee}\\
\begin{tangle}\id\step\hcoev \\
       \hx\step\SS \\
       \hh\cu\hstep\dd \\
       \hstep\hx\end{tangle}\\
\hstep\object{H}\step\object{H^\vee}
\ea
&\qquad&
\ba{l}\object{H}\Step\object{X}\\
\begin{tangle}\id\step\ld \\
       \hev\step\id\end{tangle}\\
\Step\object{X}
\ea\qquad
\ba{l}\object{X}\\
\begin{tangle}\hh\id\step\hcoev \\
       \hx\step\id \\
       \lu\step\SS \\
       \step\hx\end{tangle}\\
\step\object{H}\step\object{X}
\ea
\\
\hbox{\scriptsize a) $H_\bullet$ in $\sideset{^{H^\vee}_{H^\vee}}{}\C$}&&
\hbox{\scriptsize b) $((H_\op)^\vee)^\bullet$ in
       $\sideset{^{H_\op}_{H_\op}}{}{\mathop{\overline{\C}}}$}&&
\hbox{\scriptsize c) $({}^\Vee((\_)^\Op))_{\S\Op}:
      \sideset{^{H^\vee}_{H^\vee}}{}\C\to
      \sideset{^{H_\op}_{H_\op}}{}{\mathop{\overline{\C}}}$}
\ea
\]
\caption{}
\label{Fig-proof-invert}
\end{figure}

\begin{remark}
Applying the functor $\,\text{-}^\vee\,$ to left or right integrals on
(in) $\,H\,$, we get corresponding right or left integrals in (on)
$\,H\pti\,$. Hence, the four natural pairings of the form
$\,\Int H\tens \Int H\pti\, \xra{\miniint\otimes\miniint} \,H\otimes
H\pti \,\to\, \1\,$ and the four natural copairings
$\,\1 \,\to\, H\pti\tens H\, \xra{\miniint\otimes\miniint} \,\Int H\pti
\tens \Int H\,$ are isomorphisms.
\end{remark}

\subsection*{Group-like elements}

Let us now give the derivation of the special group-like elements (or
moduli) exploiting universality properties of integrals, and discuss
their properties as needed in the Radford formula.

\begin{proof}[Proof of Lemma \ref{prop-group-like}]
The morphism $\,\left(\intl H\otimes\id_H\right)\circ\Delta\,$ is a
left $\,(\Int H\otimes H)$-valued integral on $\,H\,$. Indeed it is
easily found from basic Hopf algebra axioms, that if $\,g:H\to X\,$ is
a left integral, then the composition
$\,H \,\xra{\,\Delta\,}\, H\tens H \,\xra{\,g\tens H\,}\, X\tens H\,$
is a left integral as well.  By the universal property for integrals
there exists a unique morphism,
$\,\Delta_r:\,\Int H\rightarrow\Int H\otimes H\,$, such that the
following diagram is commutative:
\[
\begin{CD}
H@>\Delta>>H\otimes H\\
@V{\intl H}VV @VV{\intl H\otimes\id_H}V\\
\Int H@>{\Delta_r}>>\Int H\otimes H
\end{CD}
\]
For the same reason there exists a unique morphism, $\,f\,$,
which makes  the following diagram commutative:
\[
\begin{CD}
H@>{(\id_H\otimes\Delta)\Delta=(\Delta\otimes\id_H)\Delta}>>
    H\otimes H\otimes H\\
@V{\intl H}VV @VV{\intl H\otimes\id_{(H\otimes H)}}V\\
\Int H@>f>\hphantom{(\id_H\otimes\Delta)\Delta=(\Delta\otimes\id_H)\Delta}>
                                                    \Int H\otimes H\otimes H
\end{CD}
\]
Given  that both $\,\left(\id_{(\Int H)}\otimes\Delta\right)\Delta_r\,$
and $\,\left(\Delta_r\otimes\id_{(\Int H)}\right)\Delta_r\,$
fulfill this condition for $\,f\,$, we infer from universality
that $\,\left(\id_{(\Int H)}\otimes\Delta\right)\Delta_r\,=\,
      \left(\Delta_r\otimes\id_{(\Int H)}\right)\Delta_r\,$, i.e.,
that $\,\Delta_r\,$ is a coaction.  According to
Lemma~\ref{lemma-invert-obj} any coaction, $\,\Delta_r\,$, on an
invertible object, $\,\Int H\,$, has the form
$\,\Delta_r=\id_{(\Int H)}\otimes a\,$ with a group-like morphism,
$\,a:\,\1\to H\,$.
\par
Similarly, we find that
$\,\left(\id_{H}\otimes\intr H\right)\circ\Delta\,=\,b\otimes\intr H\,$
for some group-like $\,b:\1\to H\,$. It is easy to check that
$\,\intl H\circ S\,$ is a right integral. Moreover, the coequalizer
property for right integrals implies that there exists an automorphism,
$\,t: \Int H \to \Int H\,$, such that
$\,\intl H\circ S\, = \,t\circ\intr H : H \to \Int H\,$. Composing
\eqref{a-def} with the antipode one deduces \eqref{a-1-def} for
$\,b\, = \,S^{-1} a\, = \,a^{-1}\,$.
\end{proof}

Moduli and Comoduli are dual to each other in the following sense:

\begin{corollary}
$(\alpha^{-1})_{H\pti} = a_H^t : H\pti \to \1$,
$(a^{-1})_{H\pti} = \alpha_H^t : \1 \to H\pti$.
\end{corollary}

There is a natural way in which the  group-like elements act on the
spaces $\,{\Hom}_{\cal C}(H,X)\,$ and $\,{\Hom}_{\cal C}(X,H)\,$.
Specifically, for any morphisms, $\,f:H\to X\,$ and $\,g:X\to H\,$, we
shall use the following abbreviated notations for these actions:
\begin{eqnarray*}
a.f=f\circ\mu\circ(a\otimes\id_H)\,,&\quad&
f.a=f\circ\mu\circ(\id_H\otimes a)
\enspace :\enspace H\to X\,,\\
\alpha .g=(\alpha\otimes\id_H)\circ\Delta\circ g\,,&\quad&
g.\alpha=(\id_H\otimes\alpha)\circ\Delta\circ g
\enspace :\enspace X\to H\,.
\end{eqnarray*}

Specialized to the integrals we obtain the following relations:

\begin{proposition}
\begin{eqnarray}
\label{group-like-action}
(\cll)^{-1}\left(a.\intl H\right)=
(\crl)^{-1}\cdot\intr H\,,
\quad&&\quad
(\clr)^{-1}\left(\intl H.a\right)=
(\crr)^{-1}\cdot\intr H\,,
\\
(\cll)^{-1}\left(\alpha.\cointl H\right)=
(\clr)^{-1}\cdot\cointr H\,,
\quad&&\quad
(\crl)^{-1}\left(\cointl H.\alpha\right)=
(\crr)^{-1}\cdot\cointr H\,,
\\
\label{cccc-equ}
\clr\crl=\ququ\cll\crr\,,\qquad &&\text{where}\qquad
\ququ:=\alpha\circ a\in\Aut(\1)\,.
\end{eqnarray}
\end{proposition}

\begin{figure}
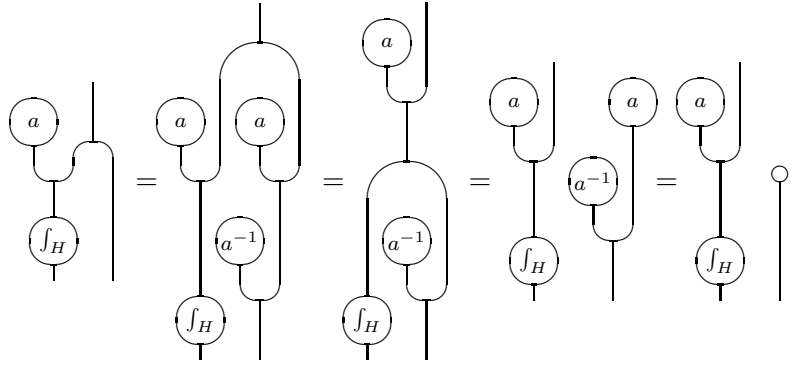

\[
\begin{tangle}\Q{a}\step\hcd \\
       \hh\cu\step\id \\
       \hstep\Ointl{H}\step\hstep\id\end{tangle}
\enspace = \enspace
\begin{tangle}\step\cd \\
       \Q{a}\step\id\step\Q{a}\step\id \\
       \hh\cu\step\cu \\
       \hstep\id\step\Q{a^{-1}}\step\id \\
       \hstep\Ointl{H}\step\hcu\end{tangle}
\enspace = \enspace
\begin{tangle}\hstep\Q{a}\step\id \\
       \hh\hstep\cu \\
       \cd \\
       \id\step\Q{a^{-1}}\step\id \\
       \Ointl{H}\step\hcu\end{tangle}
\enspace=\enspace
\begin{tangle}\Q{a}\step\id\Step\Q{a} \\
       \hcu\step\Q{a^{-1}}\step\id \\
       \hstep\Ointl{H}\hstep\step\hcu\end{tangle}
\enspace=\enspace
\begin{tangle}\Q{a}\step\id \\
       \hcu\step\unit \\
       \hstep\Ointl{H}\hstep\step\id\end{tangle}
\]
\caption{Proof of the right integral property.}
\label{Proof-int}
\end{figure}

\begin{proof}
In the case of the left equation from (\ref{group-like-action}), we can
see from the diagrammatic calculation in Figure~\ref{Proof-int}, where
the third equality follows from \eqref{a-def}, that the left hand side
is also a right integrals on $\,H\,$.  Hence $\,a.\intl H\,$ is
``proportional'' to $\,\intr H\,$ in the sense of
Proposition~\ref{pro-exist-int} by an element in in
$\,{\Aut}({\Int}H)\,$. The fact that the latter is given by the
indicated isomorphisms $\,c\in{\Aut}(\1)\,$ is seen by composing both
sides with $\,\_\circ\cointr H\,$.  The remaining three relations
between integrals are found similarly. Composition of these equalities
with $\_\circ\cointr H$ and \eqref{alpha-1-def} yield (\ref{cccc-equ}).
The scheme of this proof was taken from
\cite{Radford:antipode,Kuperberg}.
\end{proof}

For a fixed choice of $\,\Int H\,$ all integrals are defined uniquely up
to multiplication by elements of ${\rm Aut}(\1)\,$. The normalizations
can be chosen, such that   any three of the constants $\,\cll\,$,
$\,\clr\,$, $\,\crl\,$, and $\,\crr$ are unit.

\subsection*{Integrals and the antipode}

\begin{figure}
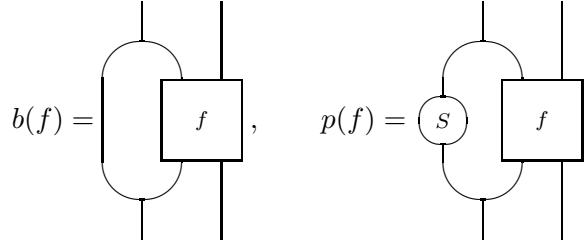

\[
b(f) =
 \begin{tangle}\cd\step\id \\
        \id\Step\Frabox{f} \\
        \cu\step\id\end{tangle}
\quad, \qquad p(f) = \quad
 \begin{tangle}\cd\step\id \\
        \S\Step\Frabox{f} \\
        \cu\step\id\end{tangle}
\]
\caption{The maps $b$, $p$.}
\label{HMHN-lem-fig}
\end{figure}

Before we consider the action of the antipode on integrals and
idempotents, we need a presentation of $\,S\,$ as
given in the following two technical but easy lemmas:

\begin{lemma}\label{HMHN-lemma}
The maps
\[ b, p: \Hom(H\tens M,H\tens N) \to \Hom(H\tens M,H\tens N) ,\]
as defined in  Figure~\ref{HMHN-lem-fig} are inverse to each other.
\end{lemma}

\begin{proof}
Straightforward.
\end{proof}

\begin{lemma}\label{S-S-lemma}
The identities for the antipode of a braided Hopf algebra, $\,H\,$,
depicted in Figure~\ref{Fig-antip-thr-int} hold true.
\end{lemma}

\begin{figure}
\[
S\otimes\id_{(\Int H)}=\left( \cll \right)^{-1}\cdot
\begin{tangle}\id\step\hstep\Ocointl{H} \\
       \hh\id\step\cd \\
       \hxx\step\id \\
       \hh\id\step\cu \\
       \id\step\hstep\Ointl{H}\end{tangle}
\qquad\qquad
\id_{(\Int H)}\otimes S=\left( \crr \right)^{-1}\cdot
\begin{tangle}\hstep\Ocointr{H}\step\hstep\id \\
       \hh\cd\step\id \\
       \id\step\hxx \\
       \hh\cu\step\id \\
       \hstep\Ointr{H}\step\hstep\id\end{tangle}
\]
\caption{}
\label{Fig-antip-thr-int}
\end{figure}

\begin{proof}
This follows from the equation derived in Figure~\ref{Fig-antip-proof},
to which we apply Lemma~\ref{HMHN-lemma}.
\end{proof}

\begin{figure}
\[
 \begin{tangle}\hstep\id\Step\Ocointl H \\
        \hh\hcd\step\hcd \\
        \id\step\hxx\step\id \\
        \hh\hcu\step\hcu \\
        \hstep\id\Step\Ointl H\end{tangle}
\enspace = \enspace
 \begin{tangle}\id\step\Ocointl H \\
        \hh\hcu \\
        \hstep\id \\
        \hh\hcd \\
        \id\step\Ointl H\end{tangle}
\enspace = \enspace
 \begin{tangle}\counit\Step\Ocointl H \\
        \Step\id \\
        \Step\id \\
        \unit\Step\Ointl H\end{tangle}
\quad = \enspace
 \begin{tangle}\hstep\id\step[1.5]\Ocointl H \\
        \hh\hcd\step\id \\
        \id\step\S\step\id \\
        \hh\hcu\step\id \\
        \hstep\id\step[1.5]\Ointl H\end{tangle}
\]
\caption{Proof of Lemma~\ref{S-S-lemma}}
\label{Fig-antip-proof}
\end{figure}

Besides the elements $\,c^{\bullet}_{\bullet}\,$ we need another
 invertible element, $\,\psipsi\in{\Aut}(\1)\,$, which determines the
self braiding of $\,\Int H\,$. It is defined by either of the following
conditions, which are equivalent by Lemma~\ref{lemma-invert-obj} and
Theorem~\ref{thm-intint=Pi}:
\begin{equation}
\Psi_{(\Int H,\Int H)}\,=\,\psipsi\cdot\id_{(\Int H\otimes \Int H)}
\,\,\qquad
\text{or}\qquad\,\,
\psipsi^{-1}\,=\,\left( \dim_8(\Int H) \right)
\,=\,\tr_8(\PPi^\bullet_\bullet(H))\,.
\end{equation}
Using cyclicity, $\,\tr_8(fg)=\tr_8(gf)\,$, one can rewrite the
latter expression for $\,\psipsi\,$ in many equivalent forms. For
example, we have  $\,\psipsi^{-1}\,=
\,\tr_8\bigl(\Delta^\op\circ\mu\circ(S^2\otimes\id)\bigr)\,$.  The
action of $\,S\,$ on the integrals is now as follows:

\begin{proposition}
\begin{eqnarray}
\label{antipode-circ-int}
S\circ\cointr H=\psipsi\cdot\clr(\cll)^{-1}\cdot\cointl H\,,
&\qquad&
S\circ\cointl H=\psipsi\cdot\crl(\crr)^{-1}\cdot\cointr H\,,
\\
\intr H\circ S=\psipsi\cdot\crl(\cll)^{-1}\cdot\intl H\,,
&\qquad&
\intl H\circ S=\psipsi\cdot\clr(\crr)^{-1}\cdot\intr H\,.
\end{eqnarray}
\end{proposition}

\begin{proof}
The composition of both sides of the first identity from Figure
\ref{Fig-antip-thr-int} with $\cointr H\otimes\id_{(\Int H)}$
gives the first identity from \eqref{antipode-circ-int}. All other
identities are found analogously.
\end{proof}

Using Theorem~\ref{thm-intint=Pi} we then find immediately the
following relations between the special idempotents and the antipode:

\begin{corollary}
\begin{eqnarray*}
S\circ\PPi^\ell_r(H) = \psipsi \PPi^\ell_\ell(H) = \PPi^r_\ell(H)\circ S\,,
&\qquad&
S\circ\PPi^r_\ell(H) = \psipsi \PPi^r_r(H) = \PPi^\ell_r(H)\circ S\,,\\
S\circ\PPi^\ell_\ell(H) = \ququ\psipsi \PPi^\ell_r(H) = \PPi^r_r(H)\circ S\,,
&\qquad&
S\circ\PPi^r_r(H) = \ququ\psipsi \PPi^r_\ell(H) = \PPi^\ell_\ell(H)\circ S\,.
\end{eqnarray*}
The idempotents $\,\PPi^\ell_\ell(H)\,$, $\,\PPi^\ell_r(H)\,$,
$\,\PPi^r_\ell(H)\,$, and $\,\PPi^r_r(H)$ commute with $\,S^2\,$.
\end{corollary}

\medskip

\subsection*{Proof of the generalized Radford formula}

\begin{figure}
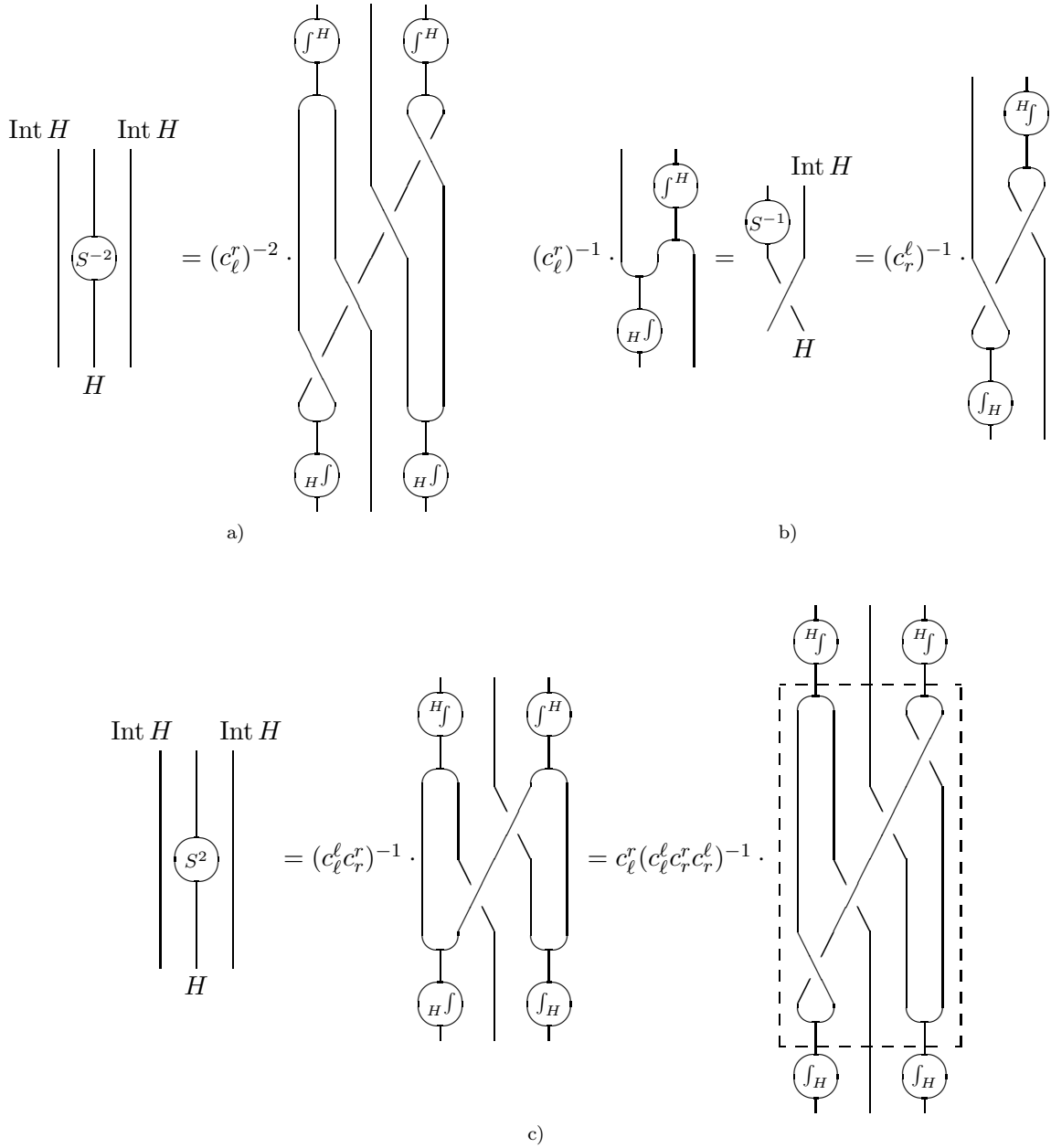

\[
\ba{ccc}
\ba{c}
\object{\Int H}\step\Step\object{\Int H}\\
\begin{tangle}\id\step\id\step\id \\
       \id\step\O{S^{-2}}\step\id \\
       \id\step\id\step\id\end{tangle}                \\
\object{H}
\ea
\enspace=
(\crl)^{-2}\cdot
\begin{tangle}\hstep\Ocointl{H}\hstep\step\id\hstep\step\Ocointl{H} \\
       \hh\cd\step\id\step\cd \\
       \id\step\id\step\id\step\hx \\
       \id\step\id\step\hx\step\id \\
       \id\step\hx\step\id\step\id \\
       \hx\step\id\step\id\step\id \\
       \hh\cu\step\id\step\cu \\
       \hstep\Ointr{H}\hstep\step\id\hstep\step\Ointr{H}\end{tangle}
&\quad\enspace&
(\crl)^{-1}\cdot
\begin{tangle}\id\step\hstep\Ocointl{H} \\
       \hh\id\step\cd \\
       \hh\cu\step\id \\
       \hstep\Ointr{H}\step\hstep\id\end{tangle}
\enspace=\enspace
\ba{l}
\step\hstep\object{\Int H}\\
\begin{tangle}\SS\step\id \\
       \hxx\end{tangle}       \\
\step\object{H}
\ea
\enspace=
(\clr)^{-1}\cdot
\begin{tangle}\id\step\hstep\Ocointr{H} \\
       \hh\id\step\cd \\
       \id\step\hxx \\
       \hx\step\id \\
       \hh\cu\step\id \\
       \hstep\Ointl{H}\step\hstep\id\end{tangle}\\
\hbox{\scriptsize a)}
&&
\hbox{\scriptsize b)}
\ea
\]
\[
\ba{c}
\ba{c}
\object{\Int H}\step\Step\object{\Int H}\\
\begin{tangle}\id\step\id\step\id \\
       \id\step\O{S^{2}}\step\id \\
       \id\step\id\step\id\end{tangle}\\
\object{H}
\ea
\enspace=
(\cll\crr)^{-1}\cdot
\begin{tangle}\hstep\Ocointr{H}\step\hstep\id\step\hstep\Ocointl{H} \\
       \hh\cd\step\id\step\cd \\
       \id\step\id\step\hxx\step\id \\
       \id\step\hxx\step\id\step\id \\
       \hh\cu\step\id\step\cu \\
       \hstep\Ointr{H}\step\hstep\id\step\hstep\Ointl{H}\end{tangle}
\enspace=
\crl(\cll\crr\clr)^{-1}\cdot\quad
\begin{tangle}\hstep\Ocointr{H}\hstep\step\id\hstep\step\Ocointr{H} \\
       \hh\cd\step\id\step\cd \\
       \id\step\id\step\id\step\hxx \\
       \id\step\id\step\hxx\step\id \\
       \id\step\hxx\step\id\step\id \\
       \hx\step\id\step\id\step\id \\
       \hh\cu\step\id\step\cu \\
       \hstep\Ointl{H}\hstep\step\id\hstep\step\Ointl{H}\end{tangle}
\unitlength\unitlens
\put(-4.5,-5){\dashbox{.3}(5,10){}}\\
\hbox{\scriptsize c)}
\ea
\]
\caption{Proof of the generalized Radford formula (part 1)}
\label{Proof-Radford1}
\end{figure}

\begin{proof}[Proof of Theorem \ref{thm-Radford-formula}]
The general scheme of our proof is taken from
\cite{Radford:antipode,Kuperberg}.  The basic steps in the proof in the
braided case are illustrated in
Figures~\ref{Proof-Radford1}--\ref{Proof-Radford3}.  The identity in
Figure~\ref{Proof-Radford1}a) is a corollary of the identities in
Figure~\ref{Fig-antip-thr-int}, rewritten for $\,H^\op\,$ and for
$\,H_\op\,$.
\par
The first identity in Figure~\ref{Proof-Radford1}b) follows from the
first or from the second identity in Figure~\ref{Fig-antip-thr-int},
rewritten for $\,H_\op\,$ or $\,H^\op\,$ respectively.
The second identity is the first one rewritten for $\,(H_\op)^\op\,$.
\par
The first identity in Figure~\ref{Proof-Radford1}c) is a result of
combining the identities in Figure~\ref{Fig-antip-thr-int}.
The second is obtained by application of the identity in
Figure~\ref{Proof-Radford1}b).
\par
In Figure \ref{Proof-Radford2} the letter $\,f\,$ denotes a morphism
in $\,\End(H^{\otimes 3})\,$, presented by the framed diagram in
Figure~\ref{Proof-Radford1}c). The first two lines of the former
figure are obtained, e.g., in a similar way as the derivation in
 Figure~\ref{Proof-int}.
\par
Denoting $\,x := \ad^\alpha\circ S^4\circ u^0_{-2}\circ\ad_a\,$, we can
rewrite the middle equation of Figure~\ref{Proof-Radford2} in a
symbolic form, using new morphisms, $\,a\,$, $\,b\,$, $\,c\,$, and
$\,d\,$, which can be visualized via Figure~\ref{Proof-Radford1}c). The
left hand side of the equation there can be transformed to the first
line in Figure~\ref{Proof-Radford3}, where $\,u^0_2\,$ denotes the
inverse to $\,u^0_{-2}\,$.  Hence, the first equation in the last line
of Figure~\ref{Proof-Radford3} holds. Using
Figure~\ref{Proof-Radford1}b) we get the last equation of
Figure~\ref{Proof-Radford3}.
\end{proof}

\begin{figure}
\[
\ba{c}
\object{\Int H}\Step\step\object{\Int H}\\
\begin{tangle}\id\step\O{S^{2}}\step\id \\
       \id\step\O{\ad^{\alpha}}\step\id\end{tangle}\\
\object{H}
\ea
\enspace
=
\crl(\cll\crr\clr)^{-1}\cdot
\begin{tangle}\hstep\kilglu
             \hbox to 0em{\put(0,0){\line(0,1){0.4}}
                          \put(0,1.6){\line(0,1){0.4}}
                          \put(0,1){\oval(2.2,1.2)}
                          \put(-0.6,0.4){\makebox(1.2,1.2)[cc]
                          {\scriptsize $\cointr H.\alpha^{-1}$}}
                                                                \hss}
             \step[1.5]\id\step[1.5]
             \kilglu
             \hbox to 0em{\put(0,0){\line(0,1){0.4}}
                          \put(0,1.6){\line(0,1){0.4}}
                          \put(0,1){\oval(2.2,1.2)}
                          \put(-0.6,0.4){\makebox(1.2,1.2)[cc]
                          {\scriptsize $\cointr H.\alpha^{-1}$}}
                                                                \hss} \\
       \hh\hstep\id\hstep\step\id\hstep\step\id \\
       \put(0,0){\dashbox{.1}(4,2){$f$}} \\
       \hh\cd\step\id\step\hstep\id \\
       \Ointl{H}\step\counit\put(.2,.2){$\alpha$}
                                      \step\id\hstep\step\Ointl{H}\end{tangle}
\qquad
=
\ququ\crr(\cll\crl\clr)^{-1}\cdot
\begin{tangle}\hstep\Ocointl{H}\hstep\step\id\hstep\step\Ocointl{H} \\
       \put(0,0){\dashbox{.1}(4,2){$f$}} \\
       \hstep\Ointl{H}\hstep\step\id\hstep\step\Ointl{H}\end{tangle}
\]
\[
\ba{c}
\object{\Int H}\Step\step\object{\Int H}\\
\begin{tangle}\id\step\O{\ad_a}\step\id \\
       \id\step\O{S^{2}}\step\id \\
       \id\step\O{\ad^{\alpha}}\step\id\end{tangle}\\
\object{H}
\ea
\enspace
=
\ququ\crr(\cll\crl\clr)^{-1}\cdot
\begin{tangle}\hstep\Ocointl{H}\hstep\step\id\step
                     \unit\put(-.8,1.3){$a^{-1}$}\step\Ocointl{H} \\
       \hh\hstep\id\hstep\step\id\step\cu \\
       \put(0,0){\dashbox{.1}(4,2){$f$}} \\
       \hstep\kilglu
             \hbox to 0em{\put(0,0){\line(0,1){0.4}}
                          \put(0,1.6){\line(0,1){0.4}}
                          \put(0,1){\oval(1.4,1.2)}
                          \put(-0.6,0.4){\makebox(1.2,1.2)[cc]
                                         {\scriptsize $a.\intl H$}}
                                                                \hss}
             \hstep\step\id\hstep\step
             \kilglu
             \hbox to 0em{\put(0,0){\line(0,1){0.4}}
                          \put(0,1.6){\line(0,1){0.4}}
                          \put(0,1){\oval(1.4,1.2)}
                          \put(-0.6,0.4){\makebox(1.2,1.2)[cc]
                                         {\scriptsize $a.\intl H$}}
                                                      \hss}\end{tangle}
\enspace
=
(\crl)^{-2}\cdot
\begin{tangle}\hstep\Ocointl{H}\hstep\step\id\hstep\step\Ocointl{H} \\
       \put(0,0){\dashbox{.1}(4,2){$f$}} \\
       \hstep\Ointr{H}\hstep\step\id\hstep\step\Ointr{H}\end{tangle}
\]
\[
 \begin{tangle}\id\Step\id\Step\id \\
        \id\step[1.5]\Frabox x\step[1.5]\id \\
        \id\Step\id\Step\id\end{tangle}
\enspace = \enspace\enspace
 \begin{tangle}\hh\hstep\id\step[1.5]\id\step[1.5]\id \\
        \hh\frabox a\step\id\step\frabox b \\
        \id\step\id\step\hxx\step\id \\
        \id\step\hxx\step\id\step\id \\
        \hh\frabox c\step\id\step\frabox d \\
        \hh\hstep\id\step[1.5]\id\step[1.5]\id\end{tangle}
\]
\caption{Proof of the generalized Radford formula (part 2)}
\label{Proof-Radford2}
\end{figure}

\begin{figure}
\[
\vcenter{\hbox to 32mm{%
\unitlength 0.50mm
\linethickness{0.4pt}
\begin{picture}(60,100)
\put(35,15){\oval(10,10)[t]}
\put(42.5,15){\oval(35,20)[lt]}
\put(45,100){\line(0,-1){7}}
\put(33,90){\line(1,0){9}}
\put(42,57.50){\oval(36,65)[r]}
\put(45,0){\line(0,1){22}}
\put(30,0){\line(0,1){15}}
\put(28,15){\oval(6,10)[lb]}
\put(32,15){\oval(16,10)[rb]}
\put(30,70){\line(0,1){17}}
\put(30,100){\line(0,-1){7}}
\put(12.25,60){\oval(25,60)[l]}
\put(33,90){\line(-1,0){15.33}}
\put(17.25,50){\oval(25,40)[rb]}
\put(23,50){\framebox(14,20)[cc]{$x$}}
\put(15,0){\line(0,1){100}}
\put(45,28){\line(0,1){60}}
\end{picture}
\hss}} \enspace = \enspace
\vcenter{\hbox to 32mm{%
\unitlength 0.50mm
\linethickness{0.4pt}
\begin{picture}(60,100)
\put(35,15){\oval(10,10)[t]}
\put(42.5,15){\oval(35,20)[lt]}
\put(10,40){\framebox(10,10)[cc]{$c$}}
\put(10,70){\framebox(10,10)[cc]{$a$}}
\put(10,70){\line(0,-1){20}}
\put(45,80){\line(0,1){7}}
\put(40,70){\framebox(10,10)[cc]{$b$}}
\put(40,40){\framebox(10,10)[cc]{$d$}}
\put(50,50){\line(0,1){20}}
\put(45,100){\line(0,-1){7}}
\put(33,90){\line(1,0){9}}
\put(20,50){\line(1,1){20}}
\put(40,50){\line(-1,2){6}}
\put(20,70){\line(1,-2){6}}
\bezier{24}(30,50)(28.67,52.83)(27.50,55)
\bezier{24}(30,70)(31.50,67.17)(32.50,65)
\put(42,57.50){\oval(36,65)[r]}
\put(45,40){\line(0,-1){12}}
\put(45,0){\line(0,1){22}}
\put(15,80){\line(0,1){20}}
\put(30,0){\line(0,1){15}}
\put(28,15){\oval(6,10)[lb]}
\put(32,15){\oval(16,10)[rb]}
\put(30,70){\line(0,1){17}}
\put(30,100){\line(0,-1){7}}
\put(12.25,60){\oval(25,60)[l]}
\put(33,90){\line(-1,0){15.33}}
\put(15,40){\line(0,-1){40}}
\put(17.25,50){\oval(25,40)[rb]}
\end{picture}
\hss}}\enspace = \enspace
\vcenter{\hbox to 22mm{%
\unitlength 0.50mm
\linethickness{0.4pt}
\begin{picture}(40,100)
\put(0,40){\framebox(10,10)[cc]{$c$}}
\put(0,70){\framebox(10,10)[cc]{$a$}}
\put(0,70){\line(0,-1){20}}
\put(30,70){\framebox(10,10)[cc]{$b$}}
\put(30,40){\framebox(10,10)[cc]{$d$}}
\put(40,50){\line(0,1){20}}
\put(10,70){\line(1,-2){10}}
\put(35,100){\line(0,-1){20}}
\put(20,70){\line(0,1){30}}
\put(35,0){\line(0,1){40}}
\put(5,0){\line(0,1){40}}
\put(5,80){\line(0,1){20}}
\put(18,58){\line(1,1){12}}
\bezier{16}(22,66)(21,68)(20,70)
\put(20,20){\oval(16,16)[]}
\put(20,0){\line(0,1){12}}
\put(20,28){\line(0,1){21.83}}
\put(20,20){\makebox(0,0)[cc]{$u^0_2$}}
\put(24.50,61){\line(1,-2){5.50}}
\bezier{28}(15.17,55.33)(12.67,52.83)(10.17,50.33)
\end{picture}
\hss}}
\]
\[
\ba{c}
\object{\Int H}\Step\step\object{H}\Step\step\object{\Int H}\\
\begin{tangle}\id\Step\step\O{\ad_a}\Step\step\id \\
       \id\Step\step\O{S^{2}}\Step\step\id \\
       \id\step\hcoev\step\O{\ad^\alpha}\Step\step\id \\
       \id\step\id\step\hx\Step\step\id \\
       \hh\id\step\ev\step\id\step\coev\step\id \\
       \id\Step\step\hx\step\id\step\id \\
       \hh\id\Step\step\id\step\ev\step\id\end{tangle}
\ea
\enspace=
(\crl)^{-2}\cdot
\begin{tangle}\hstep\Ocointl{H}\hstep\step\id\hstep\step\Ocointl{H} \\
       \hh\cd\step\id\step\cd \\
       \id\step\id\step\id\step\hxx \\
       \id\step\id\step\hxx\step\id \\
       \id\step\hx\step\id\step\id \\
       \hx\step\id\step\id\step\id \\
       \hh\cu\step\id\step\cu \\
       \hstep\Ointr{H}\hstep\step\id\hstep\step\Ointr{H}\end{tangle}
\enspace=\;
\ba{c}
\object{\Int H}\Step\step\object{\Int H}\\
\begin{tangle}\hh\id\step\id\step\id \\
       \id\step\O{S^{-1}}\step\id \\
       \id\step\hxx \\
       \id\step\hxx \\
       \id\step\id\step\id \\
       \id\step\O{S^{-1}}\step\id \\
       \hh\id\step\id\step\id\end{tangle}\\
\object{H}
\ea
\]
\caption{Proof of the generalized Radford formula (part 3)}
\label{Proof-Radford3}
\end{figure}

\begin{remark}
Under the assumption that the antipode in $\,H\,$ is an isomorphism,
Proposition~\ref{pro-exist-int} implies that the object $\Int H$ is
invertible.  Conversely, if we suppose that the object $\,\Int H\,$,
which splits the idempotent $\,\PPi^r_\ell(H)\,$, is invertible, then
$\,S^{-1}\,$ exists and can be derived from the first identity in
Figure~\ref{Proof-Radford1}b).
\end{remark}

\section{Examples of braided Hopf algebras}\label{sec-exam}

Let us furnish in this section the general theory around Radford's
formula that we have developed so far with a couple of prominent
examples of braided, monoidal categories that are not abelian. The
first is a variant of the category of surfaces and 3-cobordisms, for
which the integrals are found even without the property of split
idempotents. The second is any rigid, braided tensor category with
split idempotents, in which the {\em coend}
$\,\int^{X\in\C} X\tens X\pti\,$ exists. We shall describe in both
cases the Hopf algebras structure, integrals, and other elements
explicitly. These two categories also figure prominently in
topological quantum field theories, where the task is to functor one
into the other.

\subsection*{A topological example of a Hopf algebra}
Crane and Yetter have shown in \cite{CY} that a 3-dimensional
topological quantum field theory assigns a braided Hopf algebra to a
torus with one hole. Yetter has explained in \cite{Yet:portrait} that a
torus with one hole $\,T\,$ is a Hopf algebra in a braided category of
cobordisms between oriented connected surfaces with one hole.  The
3-dimensional manifolds which are the structure morphisms of $\,T\,$
are described in both articles \cite{CY,Yet:portrait}.

Independently Kerler~\cite{Ker:gene} recovered this structure in a very
similar category, $\,\widetilde{Cob}_3(1)\,$, whose objects are
specially selected standard surfaces, and the morphisms are
homeomorphisms classes of triples $\,(M,\psi,\sigma)\,$, where M is a
3-cobordism, $\,\psi\,$ the homeomorphisms between $\,\partial M\,$ and
the standard surfaces, and $\,\sigma\,$ the signature of a
four-manifold bounding $\,M\,$  in a standard fashion.  Using a
somewhat different topological language, not only the usual braided
tensor structure and the Hopf algebra structures associated to the
one-holed torus are found, but in addition to the elements in
\cite{CY,Yet:portrait} also the cobordisms assigned to canonical Hopf
pairings and the integrals are identified precisely in \cite{Ker:gene}.
In particular, the integrals are directly interpreted as the algebraic
elements intrinsically associated to elementary surgery.  It is also
realized in \cite{Ker:gene} that $\,\widetilde{Cob}_3(1)\,$, restricted
to connected surfaces, is generated by the morphisms of the braided
Hopf algebra structure together with the additional elements.  (The
original generators from  \cite{CY,Yet:portrait} only produce
 cobordisms embeddable into $\,\RR^3\,$). I.e., a purely algebraic
category can be defined that is freely generated by  a Hopf algebra
object, and which surjects onto $\,\widetilde{Cob}_3(1)\,$. It is
further conjectured in \cite{Ker:gene} that there is some such
definition of the algebraic category, which is in fact  isomorphic to
$\,\widetilde{Cob}_3(1)\,$ for connected surfaces.

In this  section we shall first review the explicit presentations of
the cobordisms that are used for integrals and other Hopf algebra
structures using another category, $\,\TC\,$, of a topological and
combinatorial nature, which is equivalent to
$\,\widetilde{Cob}_3(1)\,$, as proven by Kerler~\cite{Ker:bridge}. The
category $\,\TC\,$ is a subquotient of the category $\,\RT\,$ of ribbon
tangles in $\,[0,1]\times\RR^2\,$. Unlike the ordinary
category of surfaces and 3-cobordisms, its central extension
$\,\widetilde{Cob}_3(1)\,$ admits interesting representations in braided
categories built with the help of coends,  as we will see in the second
paragraph of this section.

Recall \cite{Tur:q3} that the objects of $\,\RT\,$ are non-negative
integers, and the morphisms are ribbon tangles. For our purposes we
define ribbon (or framed) tangles to be ambient isotopy classes of
smooth embeddings of rectangles $\,[0,1]\times[0,1]\,$ and annuli
$\,S^1\times[0,1]\,$ into $\,[0,1]\times\RR^2\,$, such that the edges
$\,0\times[0,1]\,$ and $\,1\times[0,1]\,$ of the rectangles are
attached to intervals on distinguished the lines
$\,0\times\RR\times0\,$ and $\,1\times\RR\times0\,$.

Furthermore, the  images of the rectangles shall be tangent to the
strip $\,[0,1]\times\RR\times0\,$ at these intervals, and the induced
isomorphism of the tangent plane to the rectangle and the tangent plane
to $\,[0,1]\times\RR\times0\,$ shall preserve the orientations.

 In the drawings of  ribbon tangles below we will represent them only
by threads, tacitly assuming the blackboard framing, and explicitly
insert powers of $\,2\pi$\n-twists $\,\rib\,$ only when necessary.  The
tensor multiplication of objects is the addition of the non-negative
integers.  The unit object is $\0$. The composition (resp. the tensor
product) of morphisms is given by  vertically stacking tangle diagrams
(resp.  horizontally juxtaposing tangle diagrams) \cite{Tur:q3}.  With
these structures it is easily seen that $\,\RT\,$ is a braided, rigid,
monoidal category.

\begin{definition}\label{def-TC}
The category $\,\TC\,$ of tangle-cobordisms is a subquotient of
 $\,\RT\,$, whose objects are {\em even}, non-negative integers and
morphisms of the form $\,\underline{2n} \,\to \,\underline{2m}\,$  are
equivalence classes of ribbon tangles obeying the following property
\dag:

\vspace*{10pt}

\noindent\ Property \dag \
\parbox{5.1in}{there are precisely $\,n+m\,$ rectangles, each of them
connects two consecutive intervals of the same distinguished line,
which is either $\,0\times\RR\times0\,$ or $\,1\times\RR\times0\,$.}

\vspace*{10pt}

A morphism in $\,\TC\,$ can be thought of as an equivalence class of
morphisms in $\,\RT\,$. The corresponding equivalence relation is
generated by the following {\em moves}. All of them are stable under
composition and tensor product and preserve Property  \dag \  so that
$\,\TC\,$ inherits these structures from $\,\RT\,$.

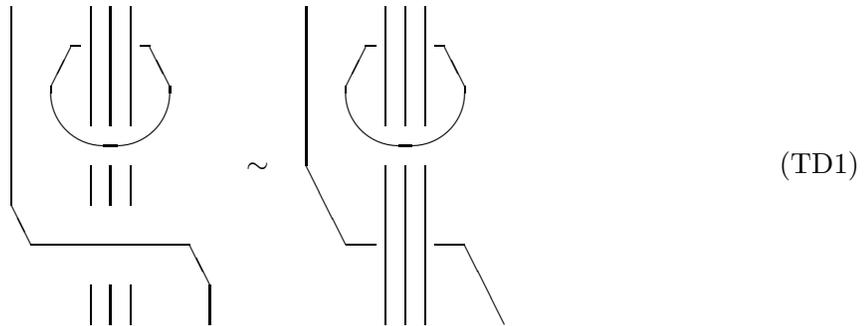
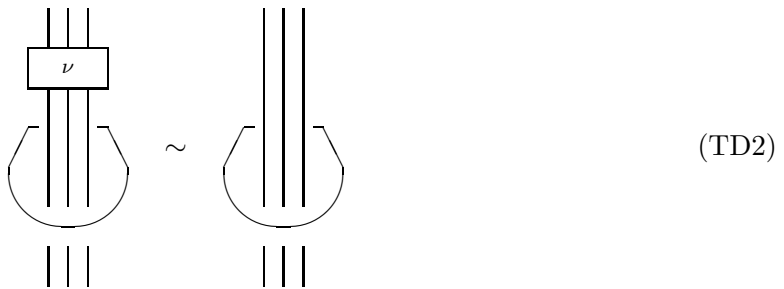
\begin{equation}
 \begin{tangle}\hh%
\id\step[1.5]\hln{.25}\step[.25]\id\hstep\id\hstep\id\step[.25]\hln{.25} \\
        \hh\id\step\dd\hstep\id\hstep\id\hstep\id\hstep\d \\
        \hh\id\step\Ev\id\hstep\id\hstep\id \\
        \hh\id \\
        \hh\id\Step\id\hstep\id\hstep\id \\
        \hh\d\hln4 \\
        \hh\step[4.5]\d \\
        \hh\Step\id\hstep\id\hstep\id\Step\id\end{tangle}
\quad\sim\quad
 \begin{tangle}\hh%
\id\step[1.5]\hln{.25}\step[.25]\id\hstep\id\hstep\id\step[.25]\hln{.25} \\
        \hh\id\step\dd\hstep\id\hstep\id\hstep\id\hstep\d \\
        \hh\id\step\Ev\id\hstep\id\hstep\id \\
        \hh\id \\
\d\hln{.75}\step[.25]\id\hstep\id\hstep\id\step[.25]\hln{.75} \\
        \Step\id\hstep\id\hstep\id\step\d\end{tangle}
\tag{TD1}\label{TD1}
\end{equation}
\begin{equation}
 \begin{tangle}\hh\step\id\hstep\id\hstep\id \\
        \hh\step\frabox\rib \\
\hh\hstep\hln{.25}\step[.25]\id\hstep\id\hstep\id\step[.25]\hln{.25} \\
        \hh\dd\hstep\id\hstep\id\hstep\id\hstep\d \\
        \hh\Ev\id\hstep\id\hstep\id \\
        \hh\vphantom{\id} \\
        \hh\step\id\hstep\id\hstep\id\end{tangle}
\quad\sim\quad
 \begin{tangle}\step\id\hstep\id\hstep\id \\
\hh\hstep\hln{.25}\step[.25]\id\hstep\id\hstep\id\step[.25]\hln{.25} \\
        \hh\dd\hstep\id\hstep\id\hstep\id\hstep\d \\
        \hh\Ev\id\hstep\id\hstep\id \\
        \hh\vphantom{\id} \\
        \hh\step\id\hstep\id\hstep\id\end{tangle}
\tag{TD2}\label{TD2}
\end{equation}
\begin{equation}
 \begin{tangle}\hh\hcoev\step\hcoev \\
        \id\step\hx\step\id \\
        \id\step\hx\step\id \\
        \hh\hev\step\hev\end{tangle}
\quad\sim\quad
\hbox{\O}
\tag{TS2}\label{TS2}
\end{equation}
\begin{equation}
\begin{array}{c}
\object{2i-1}\step[3]\object{2i} \\
 \begin{tangle}\hln5 \\
        \step\id\step[3]\id \\
        \step\id\step[3]\id \\
        \step\id\step[3]\id \\
        \step\id\step[3]\id \\
        \step\id\step[3]\id\end{tangle} \\
\dots\ \dots
\end{array}
\sim
\begin{array}{c}
\object{2i-1}\step[3]\object{2i} \\
 \begin{tangle}\hln5 \\
        \step\id\step\hcoev\step\id \\
        \step\hxx\step\hxx \\
        \hh\step\id\step\hev\step\id \\
        \hh\step\id\step\hcoev\step\id \\
        \step\hx\step\hx \\
        \step\id\step\hev\step\id\end{tangle} \\
\dots\ \dots
\end{array}
\;,\;
\begin{array}{c}
\dots\ \dots \\
 \begin{tangle}\step\id\step[3]\id \\
        \step\id\step[3]\id \\
        \step\id\step[3]\id \\
        \step\id\step[3]\id \\
        \step\id\step[3]\id \\
        \hln5\end{tangle} \\
\object{2i-1}\step[3]\object{2i}
\end{array}
\sim
\begin{array}{c}
\dots\ \dots \\
 \begin{tangle}\step\id\step\hcoev\step\id \\
        \step\hxx\step\hxx \\
        \hh\step\id\step\hev\step\id \\
        \hh\step\id\step\hcoev\step\id \\
        \step\hx\step\hx \\
        \step\id\step\hev\step\id \\
        \hln5\end{tangle} \\
\object{2i-1}\step[3]\object{2i}
\end{array}
\tag{TS3}\label{TS3}
\end{equation}
The identity morphism is the tangle
\begin{equation}
 \begin{tangle}\step[-1]\hln{14} \\
        \id\step\hcoev\step\id \\
        \hxx\step\hxx \\
        \hh\id\step\hev\step\id \\
        \hh\id\step\hcoev\step\id \\
        \hx\step\hx \\
        \id\step\hev\step\id \\
        \step[-1]\hln{14}\end{tangle}
\step[-8]\dots\Step\quad
 \begin{tangle}\id\step\hcoev\step\id \\
        \hxx\step\hxx \\
        \hh\id\step\hev\step\id \\
        \hh\id\step\hcoev\step\id \\
        \hx\step\hx \\
        \id\step\hev\step\id\end{tangle}
\quad.
\label{Id}
\end{equation}

\end{definition}

A tangle $\,\underline{2n} \,\to\, \underline{2m}\,$ represents a
cobordism between surfaces of genus $\,n\,$ and $\,m\,$.  Note also
that \eqref{TD1} is the 2-handle-slide move of any ribbon segment over
a 0\n-framed, closed, unknotted ribbon.

\begin{figure}
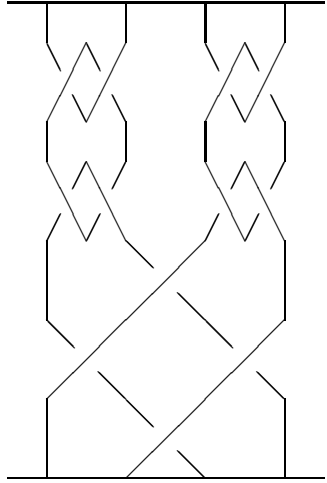

\[ \Psi: \2+\2 \to \2+\2 =\quad
 \begin{tangle}\hln8 \\
        \hh\step\id\Step\id\Step\id\Step\id \\
        \step\hxx\hxx\Step\hxx\hxx \\
        \hh\step\id\Step\id\Step\id\Step\id \\
        \step\hx\hx\Step\hx\hx \\
        \step\id\Step\xx\Step\id \\
        \step\xx\Step\xx \\
        \step\id\Step\xx\Step\id \\
        \hln8\end{tangle}
\]
\caption{Braiding in $\TC$}
\label{fig-22braid}
\end{figure}

The braiding in $\,\TC\,$ is the braiding from $\,\RT\,$ composed on
one or both sides with the identity morphism \eqref{Id} to ensure the
condition \dag. See for example Figure~\ref{fig-22braid}.

The evaluations and the coevaluations in $\,\TC\,$ are those in
$\,\RT\,$ composed with the identity morphism \eqref{Id}, as depicted
below for $\,\2\,=\,\2^{\vee}\,$:
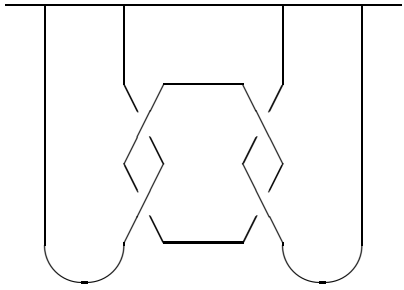
\[
\begin{array}{cccc}
 \begin{tangle}\hln{10} \\
        \step\id\Step\id\step\hln2\step\id\Step\id \\
        \step\id\Step\hxx\Step\hx\Step\id \\
        \step\id\Step\hxx\hln2\hx\Step\id \\
        \step\ev\step[4]\ev\end{tangle}
&\quad ,\quad&
 \begin{tangle}\step\coev\step\hln2\step\coev \\
        \step\id\Step\hxx\Step\hx\Step\id \\
        \step\id\Step\hxx\hln2\hx\Step\id \\
        \step\id\Step\id\step[4]\id\Step\id \\
        \hln{10}\end{tangle}
&\quad . \\
\ev: \2+\2 \to \0 && \coev: \0\to \2+\2 &
\end{array}
 \]

We assert that $\,\2\,\in\,\TC\,$ has the structure of a Hopf algebra.
Indeed, a set of structure morphisms that fulfills all the axioms of a
braided Hopf algebra is given if we choose for the multiplication and
the comultiplication the following two diagrams,
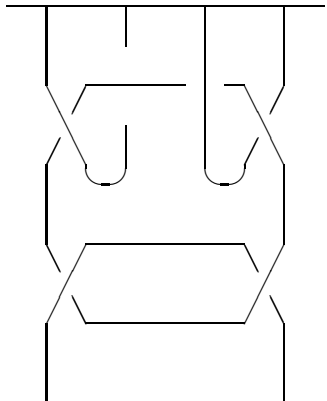
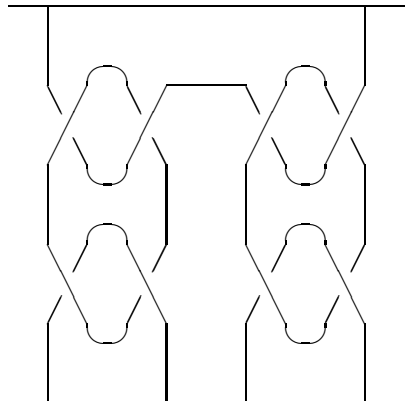
\[
\begin{array}{cccc}
 \begin{tangle}\hln8 \\
        \hh\step\id\Step\id\Step\id\Step\id \\
        \hh\step\id\step\hln{2.5}\hstep\id\hstep\hln{.5}\step\id \\
        \step\hx\step\hid\Step\id\step\hx \\
        \hh\step\id\step\hev\Step\hev\step\id \\
        \hh\step\id\step\hln4\step\id \\
        \step\hxx\hln4\hxx \\
        \step\id\step[6]\id \\
        \hln8\end{tangle}
&\quad , \quad&
 \begin{tangle}\step[-1]\hln{10} \\
        \id\step\hcoev\step\hln2\step\hcoev\step\id \\
        \hxx\step\hxx\Step\hxx\step\hxx \\
        \hh\id\step\hev\step\id\Step\id\step\hev\step\id \\
        \hh\id\step\hcoev\step\id\Step\id\step\hcoev\step\id \\
        \hx\step\hx\Step\hx\step\hx \\
        \id\step\hev\step\id\Step\id\step\hev\step\id \\
        \step[-1]\hln{10}\end{tangle}
&\quad,\\
m : \2+\2 \to \2 && \Delta : \2 \to \2+\2 &
\end{array}
\]
for the unit and the counit the ones depicted next,
\[ \eta : \0 \to \2= \qquad
 \begin{tangle}\Coev\hcoev \\
        \hxx\step\hxx \\
        \id\step\hev\step\id \\
        \step[-1]\hln5\end{tangle}
\quad  , \qquad
\varepsilon : \2 \to \0 = \quad
 \begin{tangle}\hln4 \\
        \step\ev\end{tangle}
\quad,\]
and for  the antipode and its inverse the pictures below.
\[ S = \quad
 \begin{tangle}\step[-1]\hln5 \\
        \id\step\hcoev\step\id \\
        \hxx\step\hxx \\
        \hh\id\step\hev\step\id \\
        \hh\id\step\hcoev\step\id \\
        \hxx\step\hxx \\
        \id\step\hev\step\id \\
        \step[-1]\hln5\end{tangle}
\quad, \qquad S^{-1} = \quad
 \begin{tangle}\step[-1]\hln5 \\
        \id\step\hcoev\step\id \\
        \hx\step\hx \\
        \hh\id\step\hev\step\id \\
        \hh\id\step\hcoev\step\id \\
        \hx\step\hx \\
        \id\step\hev\step\id \\
        \step[-1]\hln5\end{tangle}
\quad. \]
The bialgebra axiom is verified by employing the handle-slide
move~\eqref{TD1}. All other Hopf algebra axioms are proven
straightforwardly.

The projections of Lemma~\ref{Prop-4-idems} for this Hopf algebra are
easily found to be the following tangles:

\[ \PPi^r_l = \PPi^l_r = \PPi^l_l = \PPi^r_r = \quad
 \begin{tangle}\step[-1]\hln5 \\
        \id\step\hcoev\step\id \\
        \hx\step\hx \\
        \Ev\hev \\
        \Coev \\
        \step[-1]\hln5\end{tangle}
\quad: \2\to \2 .\]
By Theorem~\ref{thm-intint=Pi} the form of these diagrams entails that
 the object $\,\Int \2=\0\,$ of integrals is the unit object, and that
the integrals are given by the pieces
\[
 \begin{tangle}\step[-1]\hln5 \\
        \id\step\hcoev\step\id \\
        \hx\step\hx \\
        \Ev\hev\end{tangle}
\quad = \intl\2 = \intr\2 : \2\to \0\ , \quad \text{and}\qquad
 \begin{tangle}\Coev \\
        \step[-1]\hln5\end{tangle}
\quad = \cointl\2 = \cointr\2 : \0\to \2 .\]

The generalized Radford formula from Theorem~\ref{thm-Radford-formula}
yields for the Hopf algebra $\,H\,=\,\2\,$ the identity
\[
S^4 = u^0_2 = \rib^2\qquad.
\]
In fact, one can prove directly the stronger result
$\,S^2 \,= \,\rib : \,\2 \to \2\,$.

\begin{figure}
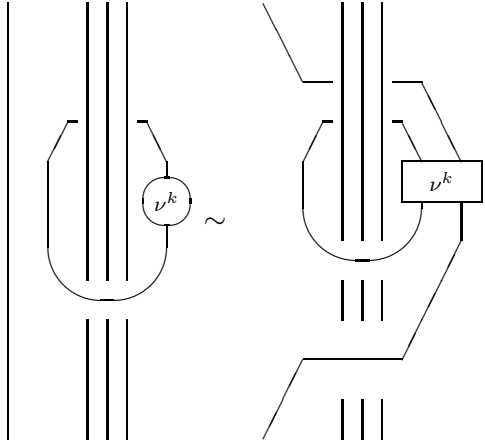

\begin{equation*}
 \begin{tangle}\id\Step\id\hstep\id\hstep\id \\
        \hh%
\id\step[1.5]\hln{.25}\step[.25]\id\hstep\id\hstep\id\step[.25]\hln{.25} \\
        \hh\id\step\dd\hstep\id\hstep\id\hstep\id\hstep\d \\
        \id\step\id\step\id\hstep\id\hstep\id\step\O{\rib^k} \\
        \hh\id\step\Ev\id\hstep\id\hstep\id \\
        \hh\id \\
        \id\Step\id\hstep\id\hstep\id \\
        \hh\id\Step\id\hstep\id\hstep\id\end{tangle}
\quad\sim\quad
\begin{tangle}\d\hln{.75}\step[.25]\id\hstep\id\hstep\id\step[.25]\hln{.75}\\
        \hh\step[1.5]\hln{.25}\step[.25]\id\hstep\id\hstep\id%
\step[.25]\hln{.25}\hstep\d \\
        \hh\step\dd\hstep\id\hstep\id\hstep\id\hstep\d\hstep\d \\
        \hh\step\id\step\id\hstep\id\hstep\id\step\frabox{\rib^k} \\
        \hh\step\Ev\id\hstep\id\hstep\id\Step\id \\
        \hh\step[4.5]\dd \\
        \hh\Step\id\hstep\id\hstep\id\step\dd \\
        \hh\step\hln{2.5}\dd \\
        \hh\hstep\dd \\
        \hh\dd\step[1.5]\id\hstep\id\hstep\id\end{tangle}
\end{equation*}
\caption{The handle-slide move with arbitrary framing}
\label{fig-arb-framing}
\end{figure}

Finally, let us mention that the  surgery calculus from
\cite{Ker:bridge} implies the Fenn Rourke and also more general
two-handle slides.  Here the case where the closed ribbon, over which
to slide, is still an unknot, but may have arbitrary framing:

\begin{proposition}\label{arb-framing}\

The handle-slide move with arbitrary framing shown on
Figure~\ref{fig-arb-framing} holds in $\,\TC\,$.
\end{proposition}

This follows from results of Kerler~\cite{Ker:bridge} and,
alternatively, can be deduced directly from relations
\eqref{TD1}--\eqref{TS3} in the category $\,\TC\,$.

\subsection*{Coends as Hopf algebras}
Let $\,(\C,\tens,\1)\,$ be a braided rigid category with split
idempotents.  Assume  that the coend $\,F=\int^{X\in\C} X\tens X\pti\,$
exists in $\,\C\,$.  See, e.g., Mac Lane's book~\cite{MacLane} for the
definition of a coend.

Whenever the coend $\,F\,$ exists in $\,\C\,$, it is a Hopf algebra
\cite{Lyu:mod,Ma:bg} with the structure morphisms given by the
following pictures or commutative diagrams. The multiplication is
determined by the condition that the diagram on the right of the
following figure commutes for any pair of objects, $\,(L,M)\,$. I.e.,
$\,m_F\,$ is the lift of the dinatural transformation defined by the
 braid on the left hand side.
\begin{equation}\label{mult3}
\unitlength=0.7mm
\raisebox{-16mm}{
\begin{picture}(46,41)
\put(1,6){\line(0,1){30}}
\put(16,6){\line(1,2){15}}
\put(31,6){\line(1,2){15}}
\put(16,36){\line(1,-1){8}}
\put(28,24){\line(1,-1){6}}
\put(38,14){\line(1,-1){8}}
\put(1,38){\makebox(0,0)[cb]{$L$}}
\put(16,38){\makebox(0,0)[cb]{$L\pti$}}
\put(31,38){\makebox(0,0)[cb]{$M$}}
\put(46,38){\makebox(0,0)[cb]{$M\pti$}}
\put(1,1){\makebox(0,0)[cb]{$L$}}
\put(16,1){\makebox(0,0)[cb]{$M$}}
\put(31,1){\makebox(0,0)[cb]{$M\pti$}}
\put(46,1){\makebox(0,0)[cb]{$L\pti$}}
\end{picture}
} \qquad \text{or} \qquad
\begin{CD}
L\tens L\pti\tens (M\tens M\pti) @>i_L\tens i_M>> F\tens F \\
@VL\tens cVV                               @V\exists Vm_FV \\
L\tens M\tens(L\tens M)\pti      @>i_{L\tens M}>> F
\end{CD}
\end{equation}
The unit is $\,\1\,=\,\1\tens\1\pti \xra{i_\1} F\,$.  The
comultiplication  $\,\Delta\,$ on $\,F\,$ is uniquely determined by the
 condition that the next equation  shall hold for any object $\,X\,$.
\begin{multline*}
\bigl(X\tens X\pti \xra{i_X} F \xra\Delta F\tens F\bigr) \\
= \bigl(X\tens X\pti = X\tens I\tens X\pti
\xra{X\tens\coev\tens X\pti} X\tens X\pti\tens X\tens X\pti
\xra{i_X\tens i_X} F\tens F\bigr)
\end{multline*}
or, pictorially,
\[ \Delta = \quad
\begin{array}{c}
\object F \\
\begin{tangle}\id\step\coev\step\id\end{tangle} \\
\object F \step[3] \object F
\end{array}
.\]
The counit $\,\varepsilon\,$ is given by the equation
\begin{equation}
\ev = \bigl(X\tens X\pti \xra{i_X} F \xra\varepsilon \1 \bigr) .
\label{couint-F}
\end{equation}
One also finds that the  antipode $\,S:F \to F\,$ exists and is
 invertible.  It is defined via the diagram
\begin{equation}\label{antgamma}
\unitlength=0.70mm
\linethickness{0.4pt}
\raisebox{-13mm}{
\begin{picture}(38,40)
\put(11,2){\line(0,1){11}}
\put(11,13){\line(6,5){17}}
\put(28,27){\line(0,1){12}}
\put(32.50,14.50){\oval(11,9)[r]}
\put(32,19){\line(-1,-4){4.33}}
\put(28,13){\line(-6,5){7}}
\put(18,21){\line(-6,5){7}}
\put(11,26.67){\line(0,1){12.33}}
\put(20,40){\makebox(0,0)[cc]{$F$}}
\put(20,1){\makebox(0,0)[cc]{$F$}}
\put(1,20){\makebox(0,0)[cc]{$S=$}}
\end{picture}
}
\end{equation}

A  pairing of Hopf algebras,  $\,\omega:F\tens F \to \1\,$, is given by
the tangle below. See also \cite{Lyu:mod}.
\begin{equation}\label{omega}
\unitlength=1mm
\linethickness{0.4pt}
\raisebox{-9mm}{
\begin{picture}(61,19)
\put(40,5){\line(0,1){4}}
\put(40,9){\line(4,5){7.33}}
\put(61,9){\line(0,1){9}}
\put(40,13){\line(-4,5){4}}
\put(28,19){\makebox(0,0)[cc]{$F$}}
\put(54,19){\makebox(0,0)[cc]{$F$}}
\put(5,10){\makebox(0,0)[cc]{$\omega\ =$}}
\put(32,9){\oval(24,18)[b]}
\put(45,10){\oval(32,20)[rb]}
\put(20,8){\line(0,1){10}}
\end{picture}
}
\end{equation}
We shall call the category $\,\C\,$ \emph{modular} in case the
pairing $\,\omega\,$  is side-invertible (i.e., non-degenerate).

Finally, it is easily seen that $\,\C\,$ is also a ribbon category.

\begin{proposition}[\cite{KerLyu}]
Let $\,\C\,$ be a modular category. Then the object of the integrals of
the Hopf algebra $\,F\,$ is isomorphic to $\,\1\,$, its integrals are
two-sided, i.e., $\,\intl F=\intr F:\, F\to\1\,$ and
$\,\cointl F\,=\,\cointr F: \,\1\to F\,$, and they are related by the
equation
\[ \intl F=\intr F =
\bigl( F = \1\tens F \xra{\cointl F\tens F}
F\tens F \xra\omega \1 \bigr)
= \bigl( F = \1\tens F \xra{\cointl F\tens F}
F\tens F \xra\omega \1 \bigr) .\]
\end{proposition}

In the case of abelian categories this was shown by
Lyubashenko~\cite{Lyu:mod}.  By Theorem~\ref{thm-intint=Pi} the number
\begin{equation}
\intl F(\cointl F) \,= \,\bigl( \1 = \1\tens\1
\xra{\cointl F\tens\cointl F} F\tens F \xra\omega \1 \bigr)\qquad \in\End\1
\label{omega-int-int}
\end{equation}
is invertible. We assume that it has a square root in $\,\End\1\,$ so
that we can rescale $\,\cointl F\,$ in a way, such that
$\,\intl F(\cointl F)=1\,$.

\medskip

The following theorem is a special case of the results proven in
\cite{KerLyu}, where abelianness was assumed although not needed for
the construction of the TQFT-functor.  We shall also outline the proof,
in which a lot of complications are avoided, since we have restricted
ourselves to the tangles associated to connected, one-holed surfaces.
For more detailed arguments the reader is referred to the original
exposition.

\begin{theorem}[see also \cite{KerLyu}]\label{thm-2-to-F}\

Let $\,\C\,$ be a modular category and assume that the number
in~\eqref{omega-int-int} is 1.  Then there is a unique monoidal functor
$\,\Phi:\,\TC\to\C\,$, for which $\,\Phi(\2)=F\,$ and
$\,\Phi(\intl\2)=\intl F\,$, compatible with the braiding and ribbon
twists, and which carries the Hopf algebra structure of $\,\2\,$ to the
Hopf algebra structure of $\,F\,$.
\end{theorem}

\begin{proof}[Outline of Proof:\ ]
(a) Denote by $\,\DRT\,$ the subcategory of $\,\RT\,$, which has even,
 non-negative integers as objects and ribbon tangles satisfying
Condition $\,\dag\,$ of Definition~\ref{def-TC} as morphisms. With no
loss of generality we may assume that $\,\C\,$ is a strict monoidal
category with $\,\lpti X=X\pti\,$, $\,X^{\vee\vee}=X\,$ and
$\,u^2_0=1\,$ (see \cite{Lyu:tan,KerLyu}).  To begin with, we construct
a functor $\,\tilde\Phi: \,\DRT \to \C\,$ following
\cite{Lyu:mod,Lyu:3inv}.

Set $\,\tilde\Phi(2n) = F^{\tens n}\,$. Consider an arbitrary tangle,
$\,T:2n \to 2m \in\DRT\,$, draw its planar diagram $\,T'\,$ with
threads, crossings and ribbon twists, $\,\rib\,$, and mark an absolute
maximum and on each of these maxima a closed thread with ends attached
to the target (bottom) line. For  an arbitrary family, $\,X_1$, \dots,
$X_n\,$, of objects of $\,\C\,$ one can then construct a morphism
\[ \phi(T'): \,\,X_1\tens X_1\pti\tens\dots\tens X_n\tens X_n\pti\, \to
\,F^{\tens m} \]
by assigning braidings to crossings, evaluations or counits to minima,
coevaluations to ordinary maxima, and the integral-element
$\,\cointl F:\,\1\to F\,$ to the special absolute maxima. By the
universal property of the coend the morphism $\,\phi(T')\,$ factorizes
through
\[
\bar\phi(T') : \,\,F^{\tens n} \simeq \int^{X_1,\dots,X_n\in\C}
X_1\tens X_1\pti\tens\dots\tens X_n\tens X_n\,\,\pti\,\,
\to F^{\tens m} .\]

In fact, this morphism depends only on $\,T\,$ and not on the choice of
maxima in $\,T'\,$, since, for example,
\[
 \begin{tangle}\hh\Step\hcoev \\
        \hh\hcoev\step\id\step\id \\
        \id\step\hxx\step\id \\
        \hxx\step\hev\end{tangle}
\ =\
 \begin{tangle}\coev \\
        \id\Step\id \\
        \id\Step\id\end{tangle}
\ = \cointl F = S\circ\cointl F = \
 \begin{tangle}\hh\hcoev \\
        \hh\id\step\id\step\hcoev \\
        \id\step\hxx\step\id \\
        \hxx\step\hev\end{tangle}
\ : \1\to F .\]
In order to see that setting $\,\tilde\Phi(T)=\bar\phi(T')\,$ gives the
desired functor $\,\tilde{\Phi}\,$, we notice that the counit of
$\,F\,$ is determined by the evaluation as in
equation~\eqref{couint-F}.

(b) We want to check that ribbon tangles, equivalent under moves
\eqref{TD1}--\eqref{TS3}, are sent to the same morphisms by
$\,\tilde{\Phi}\,$:

By the definition of integrals or by Property \eqref{a-def} of
$\,\intl F\,$ we have
\[
\begin{array}{r}
\object F\hstep \\
 \begin{tangle}\hcoev\step\id\step\id \\
        \id\step\hxx\step\id \\
        \id\step\hxx\step\id \\
        \hev\step\id\step\id\end{tangle} \\
\object F\hstep
\end{array}
\ \ \ = \ \ \
\begin{array}{r}
\object F\step[1.5] \\
 \begin{tangle}\hcoev\step\id\step\hcoev\step\id \\
        \id\step\hxx\step\id\step\id\step\id \\
        \id\step\hxx\step\id\step\id\step\id \\
        \hev\step\hev\step\id\step\id\end{tangle} \\
\object F\hstep
\end{array}
\ \ \ = \ \ \
\begin{array}{r}
\object F\step[2.5] \\
 \begin{tangle}\hcoev\step\id\step\id \\
        \id\step\hxx\step\id \\
        \id\step\hxx\step\id \\
        \hh\hev\step\hev\step\frabox\eta \\
        \hh\step[4]\id\step\id\end{tangle} \\
\object F\hstep
\end{array}
\]
If we now compose this identity  with the coaction
$\,\delta\,$  of $\,F\,$ with respect to  an arbitrary
object $\,Y\in\C\,$, pictured below,
\[ \delta = \
\begin{array}{l}
\object Y \\
 \begin{tangle}\hh\id\step\coev \\
        \hh\frabox{i_Y}\step\id \\
        \hh\id\step\id\step\id\end{tangle} \\
\hstep\object F\step[1.5]\object Y
\end{array}
\ : Y\to F\tens Y \]
we obtain
\[
\begin{array}{r}
\object Y\Step \\
 \begin{tangle}\hcoev\step\id \\
        \id\step\hxx \\
        \id\step\hxx\step\hcoev \\
        \hh\hev\step\frabox{i_Y}\step\id \\
        \hh\Step\id\step\id\step\id\end{tangle} \\
\object F\step[1.5]\object Y
\end{array}
\qquad = \qquad
\begin{array}{r}
\object Y \\
 \begin{tangle}\Step\hcoev\step\id \\
        \Step\id\step\hxx \\
        \Step\id\step\hxx \\
        \hh\frabox\eta\step\hev\step\id \\
        \hh\id\step\id\step[3]\id\end{tangle} \\
\object F\step[3.5]\object Y
\end{array}
\]
Hence,
\[
\begin{array}{l}
\step[1.5]\object X\step[1.5]\object Y \\
 \begin{tangle}\hh\step[1.5]\d\step\id \\
        \hcoev\step\hxx \\
        \id\step\hxx\step\id \\
        \id\step\hxx\step\id \\
        \hev\step\hxx \\
        \step\dd\step\d \\
        \hh\hstep\dd\step\frabox\eta\step\d \\
        \hh\dd\step[1.5]\id\step\id\step[1.5]\d\end{tangle} \\
\object X\step[2.5]\object F\step[2.5]\object Y
\end{array}
\ = \
\begin{array}{l}
\step\object X\Step\object Y \\
 \begin{tangle}\step\d\step\id \\
        \hcoev\step\hxx \\
        \id\step\hxx\step\id \\
        \id\step\hxx\step\id \\
        \hev\step\hxx \\
        \hh\step[1.5]\dd\step\d\step\hcoev \\
        \hh\step\dd\Step\frabox{i_Y}\step\id \\
        \hh\hstep\dd\step[2.5]\id\step\id\step\id\end{tangle} \\
\hstep\object X\step[3.5]\object F\step[1.5]\object Y
\end{array}
\ = \
\begin{array}{l}
\object X \\
 \begin{tangle}\id \\
        \id \\
        \id \\
        \id \\
        \hh\id\step\frabox\eta \\
        \hxx\step\id \\
        \hxx\step\id\end{tangle} \\
\object X\step[1.5]\object F
\end{array}
\step[-1.5]
\begin{array}{r}
\object Y \\
 \begin{tangle}\hcoev\step\id \\
        \id\step\hxx \\
        \id\step\hxx \\
        \hev\step\id \\
        \Step\id \\
        \hh\Step\id \\
        \Step\id\end{tangle} \\
\object Y
\end{array}
\ = \
\begin{array}{l}
\object X \\
 \begin{tangle}\id \\
        \id \\
        \id \\
        \id \\
        \hh\id\step\frabox\eta \\
        \id\step\id\step\id \\
        \id\step\id\step\id\end{tangle} \\
\object X\step[1.5]\object F
\end{array}
\step[-1.5]
\begin{array}{r}
\object Y \\
 \begin{tangle}\hcoev\step\id \\
        \id\step\hxx \\
        \id\step\hxx \\
        \hev\step\id \\
        \Step\id \\
        \hh\Step\id \\
        \Step\id\end{tangle} \\
\object Y
\end{array}
\]
for any $X,Y\in\Obj\C$. Therefore, $\,\tilde{\Phi}\,$ is stable under
Move \eqref{TD1}.

Similarly, for any $\,Y\in\Obj\C\,$ the identities
\[
\begin{array}{r}
\object Y \\
 \begin{tangle}\hcoev\step\O\rib \\
        \id\step\hxx \\
        \id\step\hxx \\
        \hev\step\id \\
        \hh\frabox\eta\step\id \\
        \hh\id\step\id\step\id\end{tangle} \\
\object F\step[1.5]\object Y
\end{array}
\ \ = \ \
\begin{array}{r}
\object Y\Step \\
 \begin{tangle}\hcoev\step\id \\
        \id\step\hxx \\
        \id\step\hxx \\
        \hev\step\O\rib\step\hcoev \\
        \hh\Step\frabox{i_Y}\step\id \\
        \hh\Step\id\step\id\step\id\end{tangle} \\
\object F\step[1.5]\object Y
\end{array}
\ \ = \ \
\begin{array}{r}
\object Y \\
 \begin{tangle}\hcoev\step\id \\
        \id\step\hxx \\
        \id\step\hxx \\
        \hh\hev\step\id \\
        \hh\frabox\eta\step\id \\
        \O\rib\step\id\step\id\end{tangle} \\
\object F\step[1.5]\object Y
\end{array}
\ \ = \ \
\begin{array}{r}
\object Y \\
 \begin{tangle}\hcoev\step\id \\
        \id\step\hxx \\
        \id\step\hxx \\
        \hh\hev\step\id \\
        \hh\frabox\eta\step\id \\
        \id\step\id\step\id\end{tangle} \\
\object F\step[1.5]\object Y
\end{array}
\]
imply Move \eqref{TD2}. Moreover, Move \eqref{TS2} holds since
$\,\intl F(\cointl F)=1\,$, see equation~\eqref{omega-int-int}.

To prove the remaining relation \eqref{TS3} notice that
\[
 \begin{tangle}\hcoev\step\hcoev\step\id \\
        \id\step\hx\step\hxx \\
        \id\step\hx\step\hxx \\
        \hev\step\hev\step\id\end{tangle}
\ = \
 \begin{tangle}\id \\
        \id \\
        \id \\
        \id\end{tangle}
\step
 \begin{tangle}\hcoev\step\hcoev \\
        \id\step\hx\step\id \\
        \id\step\hx\step\id \\
        \hev\step\hev\end{tangle}
\ = \
 \begin{tangle}\id \\
        \id \\
        \id \\
        \id\end{tangle}
\]
by the  handle-slide moves \eqref{TD1} and \eqref{TS2}. Therefore,
using \eqref{TD1} once again, we get
\[
\begin{array}{c}
\object F\step[6]\object F \\
 \begin{tangle}\id\step\id\step\hcoev\step\hcoev\step\id\step\id \\
        \id\step\hxx\step\hx\step\hxx\step\id \\
        \id\step\hxx\step\hx\step\hxx\step\id \\
        \hev\step\hev\step\hev\step\id\step\id \\
        \step[6]\hev\end{tangle}
\end{array}
\ = \
\begin{array}{r}
\object F\step[4]\object F\hstep \\
 \begin{tangle}\hh\id\step[4.5]\d\step\id\step\id \\
        \id\step\hcoev\step\hcoev\step\hxx\step\id \\
        \id\step\id\step\hx\step\hxx\step\id\step\id \\
        \id\step\id\step\hx\step\hxx\step\id\step\id \\
        \d\hev\step\hev\step\hxx\step\id \\
        \hh\step\d\hln3\dd\step\hev\end{tangle}
\end{array}
\ = \
\begin{array}{c}
\object F\Step\object F \\
 \begin{tangle}\id\step\id\step\id\step\id \\
        \id\step\hxx\step\id \\
        \id\step\hxx\step\id \\
        \hev\step\id\step\id \\
        \Step\hev\end{tangle}
\end{array}
\]
The non-degeneracy of the form $\,\omega\,$ implies now
Move~\eqref{TS3}.  This completes the construction of the  functor
$\,\Phi:\TC\to\C\,$.

(c) The graphical formulae for the structure maps for $\,\2\,$ differ
from those for $\,F\,$ by expressions as in~\eqref{Id}. Since we have
proven \eqref{TS3}, the structure morphisms of $\,\2\,$ are sent by
$\,\Phi\,$ to structure morphisms of $\,F\,$.
\end{proof}

The meaning of Theorem~\ref{thm-2-to-F} is that all calculations using
the handle-slide move with arbitrary framing (see
Proposition~\ref{arb-framing}) are valid for $\,F\,$.  As a corollary
one gets modular relations as in \cite{Lyu:mod}, representations of
mapping class groups of surfaces with one hole as in
\cite{Lyu:3inv,MatPol}, and, eventually,  topological quantum field
theories \cite{KerLyu}.

\section{Integrals and related structures}

The purpose of this section is to explain how the structures of
integrals and moduli are modified by  various procedures that can be
applied to create new braided Hopf algebras from known ones. Often new
Hopf algebras are obtained by taking some type of product between two
given ones. The first example of such a product is the {\em Heisenberg
double} $\,\mathcal H(H)=H\#{}^\vee H\,$, where a Hopf algebra,
$\,H\,$, is combined with its dual. This is of particular interest
since Hopf modules of  $\,H\,$ will turn out to be in a one-to-one
 correspondence with ordinary modules of  $\,\mathcal H(H)\,$.
We shall also discuss a very general type of cross product,
$\,A\ltimes B\,$, between two braided Hopf algebras. Finally, we will
also describe the procedure of ``transmuting'' the coproduct structure
of a Hopf algebra, an idea, which has its origin in the Tannaka Krein
theory for braided categories.

\subsection*{Heisenberg double, Hopf modules and integration.}

For a Hopf algebra, $\,H\,$, in $\,\C\,$ with a dual, $\,{}^\vee H\,$,
one can define the algebra $\,\mathcal H(H)=H\#{}^\vee H\,$, called the
{\em Heisenberg double}, in the following way:

$\,{}^\vee H\,$ becomes an algebra in the category
$\,\overline\C_{H_\op}\,$ with right $\,H_\op$\n-modules structure
given by the composition
\begin{equation*}
{}^\vee H\otimes H\;\xra{\,\Delta_{({}^\vee H)}\otimes\id_H\,}\;
 {}^\vee H\otimes{}^\vee H\otimes H\;\xra{\,\id_{({}^\vee H)}\otimes\ev\,}\;
 {}^\vee H\,.
\end{equation*}
The multiplication on the corresponding smash-product algebra
$\,H\#{}^\vee H\,$ is shown on Figure~\ref{Fig-Heisenberg}a).
Similarly, $\,H\,$ becomes an algebra in
$\,{}_{({}^\vee H)_\op}\overline\C\,$.  The corresponding smash-product
algebra $\,H\#{}^\vee H\,$ is the same.

\begin{figure}
\[
\ba{ccc}
\ba{l}
\object{H}\step\object{{}^\vee H}\step\hstep\object{H}\step\hstep
\object{{}^\vee H}\\
        \begin{tangle}\hh\id\step\id\step\cd\step\id \\
               \id\step\id\step\hx\step\id \\
               \id\step\hx\step\id\step\id \\
               \id\step\id\step\ru\dd \\
               \hh\cu\step\cu\end{tangle}
\ea
=
\ba{c}
\object{H}\step\hstep\object{{}^\vee H}\Step\object{H}\step\hstep
\object{{}^\vee H}\\
\begin{tangle}\hh\id\step\cd\step\cd\step\id \\
       \id\step\d\hev\dd\step\id \\
       \d\step\hx\step\dd \\
       \hh\step\cu\step\cu\end{tangle}
\ea
&\qquad\qquad&
f=\enspace
\ba{l}
\object{H}\Step\step\object{{}^\vee H}\\
\begin{tangle}\id\step\hcoev\step\id \\
       \hh\cu\step\cu\end{tangle}\\
\hstep\object{H}\Step\object{{}^\vee H}
\ea
\qquad
f^{-1}=\enspace
\begin{tangle}\id\step\hcoev\step\id \\
       \id\step\S\step\id\step\id \\
       \hh\cu\step\cu\end{tangle}
\\
\hbox{\scriptsize a) Heisenberg double}&&
\hbox{\scriptsize b) An isomorphism $f:H\#{}^\vee H\to H\otimes{}^\vee H$}
\ea
\]
\caption{}
\label{Fig-Heisenberg}
\end{figure}

Note that the Heisenberg double $\,{\mathcal H}(H)\,$ is isomorphic to
the ``matrix algebra'' of $\,H\,$, which is
\begin{equation*}
\bigl(\,H\otimes {}^\vee H,\,
 \mu_{(H\otimes {}^\vee H)}:=\id_H\otimes\ev\otimes\id_{{}^\vee H}\,,
 \eta_{(H\otimes {}^\vee H)}:=\coev\bigr)\,.
\end{equation*}
The existence of such an isomorphism for an ordinary Hopf algebra is a
special case of the duality theorem from \cite{BlaMon:dua}. See
\cite{Bes:duality} for the braided setting.

\begin{proposition}\

The morphism $\,f:\,H\#{}^\vee H\to H\otimes{}^\vee H\,$, shown on
Figure \ref{Fig-Heisenberg}b) is an algebra isomorphism.
\end{proposition}
\smallskip

\noindent
The category of Hopf $\,H$-modules is identified with the category of
$\,\mathcal H(H)$- modules as follows:

\begin{proposition}
\label{Hopf-Heisenberg}
The categories $\,{}^H_H\C\,$ and $\,{}_{(H\#{}^\vee H)}\C\,$ are
isomorphic.  The isomorphism is the identity on the  underlying objects
and morphisms in $\,\cal C\,$, and it turns a Hopf $\,H$-module,
$\,(X,\mu_\ell,\Delta_\ell)\,$, into an $\,(H\#{}^\vee H)$-module, for
which the action given by the composition
\begin{equation*}
H\otimes{}^\vee H\otimes X\;
\xra{\,\id_{(H\otimes{}^\vee H)}\otimes\Delta_\ell\,}
\; H\otimes{}^\vee H\otimes H\otimes X\;
\xra{\,\id_H\otimes\ev\otimes\id_X\,}\;
 H\otimes X\xra{\mu_\ell} X\,\,\,.
\end{equation*}
\end{proposition}

In \cite{Cryss} an approach to integration on braided Hopf algebras
involving Heisenberg doubles is developed.  There the following two
idempotents (``vacuum projectors'') in $\,\mathcal H(H)\,$ play a
crucial role:
\begin{equation*}
\begin{split}
\mathcal E&\,:=\,
(S_H\otimes\id_{{}^\vee H})\circ\coev\,=\,
(\id_H\otimes S_{{}^\vee H})\circ\coev\,,\\
\overline{\mathcal E}&\,:=\,\,
 \mu_{\mathcal H(H)}\circ(S^2_H\otimes\id_{{}^\vee H})\circ
  \Psi^{-1}_{H,{}^\vee H}\circ\coev\,.
\end{split}
\end{equation*}
Note that the elements $\,\mathcal E\,$ and $\,\overline{\mathcal E}\,$
act on a left $\,\mathcal H(H)$-module ($\equiv$ left Hopf $H$-module),
 $\,X\,$, in the following way:
\begin{equation}
\mu_\ell\circ({\mathcal E}\otimes\id_X)\,=\,\Pi^\ell_\ell(X)\,,\qquad\,\,\,\,
\mu_\ell\circ(\overline{\mathcal E}\otimes\id_X)\,=\,\PPi^\ell_\ell(X)\,.
\end{equation}
An idempotent, $\,\PPi^\ell_\ell(H)\,$,
 is derived in \cite{Cryss} from the formula
$\,\,i_H\circ\PPi^\ell_\ell(H)\,=\,\mu^{(3)}_{\mathcal H(H)}\circ
 (\overline{\mathcal E}\otimes i_H\otimes\mathcal E)\,\,$, where
$\,\,i_H\,=\,\id_H\otimes\eta_{{}^\vee H}\,:\,\,H\to\mathcal H(H)\,$ is
a canonical embedding of algebras.

\subsection*{Integrals, cross products, and transmutation.}
\label{sec-cross}

Cross products and transmutation are basic constructions that allow us
to obtain new braided Hopf algebras from given ones
\cite{Ma:tra,Ma:introm}.  Here we use generalizations of these
constructions from \cite{Bes:cross}, and show that integrals on cross
product Hopf algebras and on transmutation Hopf algebras are obtained
in a simple way from integrals on the initial Hopf algebras.

\begin{figure}
\[
\ba{ccc}
\ba{c}
\object{X}\hstep\step\object{A}\hstep\\
        \begin{tangle}\hh\id\step\cd \\
               \hxx\step\id \\
               \id\step\k \\
               \hxx\step\id \\
               \hh\id\step\cu\end{tangle}\\
       \object{X}\hstep\step\object{A}\hstep
\ea
\enspace ={}\enspace
\ba{c}
\object{X}\Step\hstep\object{A}\hstep\\
        \begin{tangle}\rd\step\hcd \\
               \id\step\hxx\step\id \\
               \ru\step\hcu\end{tangle}\\
\object{X}\Step\hstep\object{A}\hstep
\ea
&\qquad&
\ba{l}
\object{A}\step\object{B}\step\hstep\object{A}%
\step\hstep\object{B}\\
        \begin{tangle}\id\step\id\step\hcd\step\id \\
               \id\step\hxx\step\id\step\id \\
               \id\step\id\step\ru\step\id \\
               \hcu\step\cu\end{tangle}\\
        \hstep\object{A}\hstep\Step\object{B}
\ea
\qquad
\ba{l}
\hstep\object{A}\hstep\Step\object{B}\\
        \begin{tangle}\hcd\step\cd \\
               \id\step\id\step\rd\step\id \\
               \id\step\hxx\step\id\step\id \\
               \id\step\id\step\hcu\step\id\end{tangle}\\
\object{A}\step\object{B}\step\hstep\object{A}\step\hstep\object{B}
\ea
\qquad
\ba{l}
\object{A}\step\object{B}\\
        \begin{tangle}\id\step\rd \\
               \hxx\step\id \\
               \hh\id\step\cu \\
               \S\step\hstep\S \\
               \hh\id\step\cd \\
               \hxx\step\id \\
               \id\step\ru\end{tangle}\\
\object{A}\step\object{B}
\ea
\\
\hbox{\scriptsize a) Crossed module axiom}&&
\hbox{\scriptsize b) Cross product $A\ltimes B$}
\ea
\]
\caption{}
\label{Fig-Cross-Prod}
\end{figure}

\begin{proposition}[\cite{Bes:cross}]
Let $\,A\,$ be a Hopf algebra in $\,\C\,$, and $\,B\,$ be a Hopf
algebra in the category $\,{\DY\C}^A_A\,$ of crossed modules.  (A
crossed module is an object with both module and comodule structures
satisfying the compatibility condition presented on Figure
\ref{Fig-Cross-Prod}a)).

Then $\,A\ltimes B\,$ with underlying object $\,A\otimes B\,$,
multiplication, comultiplication, and antipode as shown on the Figure
\ref{Fig-Cross-Prod}b) is a Hopf algebra in $\,\C\,$.
\end{proposition}

$\,\Int B\,$ is an invertible object in $\,\C\,$ and, therefore, the
$\,A\,$-(co)module structures on $\,\Int B\,$ are given by the formulae
\[
\mu^B_r\,=\,\id_B\otimes \alpha_{B/A}   \,,\qquad\,\,\,\,
\Delta^B_r\,=\,\id_B\otimes a_{B/A}
\]
for a certain multiplicative functional, $\,\alpha_{B/A}\,$, and some
group-like element, $\,a_{B/A}\,$.  The crossed module axiom for
$\,\Int B\,$ then takes on the form
\[
\monodromy{\Int B}{A}\,\,=\,\,\,\ad^{\alpha_{B/A}}\circ\ad_{a_{B/A}}\,\,\,\,.
\]

\begin{proposition}
For any cross product Hopf algebra, $A\,\ltimes B\,$, in a rigid
braided monoidal category, $\,\C\,$, we have that
$\,\Int (A\ltimes B)\,\simeq\,\Int A\otimes\Int B\,$, and the following
 are valid choices for integrals in or on $A\,\ltimes B\,$:
\begin{eqnarray}
\intl{A\ltimes B}=\intl A\otimes\intl B\,,
&\qquad&
\ \vspace*{.2cm}\intr{A\ltimes B}=(\intr A).a_{B/A}\otimes\intr B\,,
\nonumber\\
\cointl{A\ltimes B}=\cointl A\otimes\cointl B\,,
&\qquad&\
\vspace*{.2cm}\cointr{A\ltimes B}=(\cointr A).\alpha_{B/A}\otimes\cointr B\,,
\label{int-cross-prod}
\\
a_{A\ltimes B}=a_A\cdot a_{B/A}\otimes a_B\,,
&\qquad&
\alpha_{A\ltimes B}=\alpha_A\cdot \alpha_{B/A}\otimes a_B\,.
\nonumber
\end{eqnarray}
\end{proposition}

The notion of braided quantum groups (quasitriangular braided Hopf
 algebras) was introduced in \cite{Ma:tra}, and the basic theory was
developed there.  We use the slightly modified definitions from
\cite{Bes:cross}, which reflect the symmetry between the two coalgebra
structures under consideration.  A pair of Hopf algebras,
$\,(A,\overline A)\,$, in $\,\C\times\overline\C\,$ and a  bialgebra
copairing,
$\,{\mathcal R}:\,\,\1\to\overline A_\op\otimes A\,$, in $\,\C\,$,
which is invertible as an algebra element, define a {\em braided
qua\-si\-tri\-angu\-lar Hopf algebra},
$\,(A,\overline A,{\mathcal R})\,$, if $\,A\,$ and $\,\overline A\,$
only differ in their comultiplication and antipode, i.e., if
$\,A=(A,\mu,\eta,\Delta,\varepsilon ,S)\,$ and
$\,\overline A=(A,\mu,\eta,\overline\Delta,\varepsilon ,\overline
S)\,$, and if the identity
$\,(\Psi_{A,A}\circ\overline\Delta)\cdot{\mathcal R}= {\mathcal
 R}\cdot\Delta\,$ holds. See Figure \ref{Fig-trans}a) for an
illustration of the latter.
Similar to the case of ordinary quantum groups it is shown in
\cite{Bes:cross} that the antipode of any quasitriangular Hopf algebra
of the form $\,(A,\overline A,{\mathcal R})\,$ is invertible.  In
particular, we have $\,\overline S^{-1}=u\cdot S\cdot u^{-1}\,,$ where
$\,u:=\mu\circ(\id_A\otimes S)\circ{\mathcal R}\,$.

For a quasitriangular bialgebra, $\,(A,\overline A,{\mathcal R})\,$, we
define the category $\,{\C}_\cO{A,\overline A}\,$ as the full
subcategory of the category $\,\C_A\,$ of $\,A$-right modules, whose
objects $\,(X, \mu_r)\,$ satisfy the identity indicated  in Figure
\ref{Fig-trans}b).
This category is identified with a full braided monoidal subcategory of
$\DY{\cal C}^A_A$ \cite{Bes:cross,Ma:tra}.  I.e., every module
$\,(X, \mu_r)\,$ in $\,\C_\cO{A,\overline A}\,$ becomes a crossed
module, $\,(X, \mu_r,\Delta_r)\,$, in $\,\DY\C^A_A\,$ through the
coaction
$\,\Delta_r:=(\mu_r\otimes\id_A)\circ(\id_X\otimes{\mathcal R})\,$.
Braiding in the category is given by the diagram in
Figure~\ref{Fig-trans}c).

\begin{definition}
\label{transmutation}
Let $\,(A,\overline A,{\mathcal R}_A)\,$ be a quantum braided group,
$\,(H,\mu_H,\Delta_H)\,$ a Hopf algebra in $\,\C\,$, and
$\,f:\, A\rightarrow H\,$ a bialgebra morphism.

We say that $\,\left((A,\overline A,{\cal R}_A),f,H\right)\,$ are
{\em transmutation data for a Hopf algebra $\,H\,$ in $\,\C\,$} if
\begin{itemize}
\item the two adjoint actions of $\,A\,$ and of $\,\overline A\,$,
defined through $\,f\,$ as presented on Figure \ref{Fig-trans}d),
 coincide, (We will denote this action by $\,\mu^H_\ad\,$ and call it
the {\em adjoint action of a quantum group}.)

\item and
$\,(H,\mu^H_\ad)\,$ is an object of $\,\C_{\cO{A,\overline A}}\,$.
\end{itemize}
The {\em  transmutation $\,\underline H\,$ of a Hopf algebra $\,H\,$}
is the underlying algebra of $\,H\,$ with a new comultiplication,
$\,\Delta_{\underline H}\,$, and a new antipode $\,S_{\underline H}\,$,
as  defined in Figures \ref{Fig-trans}e-f).  There we set
$\,{\mathcal R}_A^\sim=(\id_A\otimes S)\circ{\mathcal R}_A
 =(\overline S^{-1}\otimes\id_A)\circ{\mathcal R}_A\,$.
\end{definition}

\begin{figure}
\[
\ba{ccccc}
\begin{tangle}\hcd\put(-.2,.5){$\overline\Delta^\op$}\step\r \\
       \id\step\hxx\Step\id \\
       \hcu\step\cu\end{tangle}
\enspace=\enspace
\begin{tangle}\r\step\hcd \\
       \id\Step\hxx\step\id \\
       \cu\step\hcu\end{tangle}
&\qquad&
\ba{c}
\object{X}\step\hstep\object{A}\hstep\\
\begin{tangle}\hh\id\step\cd\hstep\put(-.25,.4){$\Delta$} \\
       \hxx\step\id \\
       \id\step\ru\end{tangle}\\
        \object{A}\step\object{X}\step
\ea
\enspace=\enspace
\ba{c}
\object{X}\step\hstep\object{A}\hstep\\
\begin{tangle}\hh\id\step\cd\hstep\put(-.25,.4){$\overline\Delta$} \\
       \hx\step\id \\
       \id\step\ru\end{tangle}\\
\object{A}\step\object{X}\step
\ea
&\qquad&
\begin{tangle}\id\step\id\hstep\ro{\mathcal R} \\
       \hh\id\step\hru\step\hstep\dd \\
       \hxx\step\dd \\
       \id\step\ru\end{tangle}
\\
\hbox{\scriptsize a)
 $\overline\Delta^{\op}\cdot{\mathcal R}={\mathcal R}\cdot\Delta$}&&
\hbox{\scriptsize b) $X\in\Obj\left(\C_\cO{A,\overline A}\right)$}
&&
\hbox{\scriptsize c) $\Psi$ in $\C_\cO{A,\overline A}$}
\ea
\]
\[
\ba{ccccccc}
\ba{l}
\object{H}\step[1.5]\object{A}\\
\begin{tangle}\id\step[1.5]\O{f} \\
       \hh\id\step\cd\put(.1,.3){$\Delta$} \\
       \hxx\step\id \\
       \S\step\hddcu \\
       \hh\cu\end{tangle}
\ea
\enspace=\enspace
\ba{l}
\object{H}\step[1.5]\object{A}\\
\begin{tangle}\id\step[1.5]\O{f} \\
       \hh\id\step\cd\put(.1,.3){$\overline\Delta$} \\
       \hx\step\id \\
       \O{\overline S}\step\hddcu \\
       \hh\cu\end{tangle}
\ea
&\qquad&
\begin{tangle}\step\cd\put(-.5,1.2){$\Delta_H$} \\
       \dd\ro{{\mathcal R}_A^\sim}\d \\
       \ru\put(0,.2){$\mu^H_\ad$}\Step\O{f}\step\id \\
       \hh\id\Step\step\cu\end{tangle}
&\qquad&
\begin{tangle}\id\step\ro{{\mathcal R}_A} \\
       \ru\put(0,.2){$\mu^H_\ad$}\Step\O{f} \\
       \S\Step\dd \\
       \cu\end{tangle}
&\qquad&
\begin{tangle}\id\step\ro{{\mathcal R}_A} \\
       \ru\put(-.2,.4){$\mu^H_\ad$}\step\dd \\
       \O{S^2_H}\step\dd \\
       \ru\put(-.2,.4){$\mu^H_\ad$}\end{tangle}
\\
\hbox{\scriptsize d) adjoint action $\mu^H_\ad$}
&&
\hbox{\scriptsize e) $\Delta_{\underline H}$}&&
\hbox{\scriptsize f) $S_{\underline H}$}&&
\hbox{\scriptsize g) $S_{\underline H}^2$}
\ea
\]
\caption{}
\label{Fig-trans}
\end{figure}

\begin{proposition}[\cite{Bes:cross}]\

The transmutation $\,\underline H\,$ of a Hopf algebra, $\,H\,$,
in $\,\C\,$ is a Hopf algebra in $\,\C_\cO{A,\overline A}\,$.
\end{proposition}

\smallskip

In addition to this observation from \cite{Bes:cross}, we shall next
determine also the integrals of $\,\underline H\,$:

\begin{proposition}
The integral object $\,\Int {\underline H}\,$ for the transmutation can
be chosen as the object $\,\Int H\,$, equipped with the $\,A$-action
$\,\id_{(\Int H)}\otimes\alpha_f\,$, where
$\,\alpha_f:=\alpha_H\circ f\,$.  Consistent with that one can further
put
\[
\cointl{\underline H}=\cointl H\,,\,\,\,
\cointr{\underline H}=\cointr H\,,\,\,\,
\intr{\underline H}=(a^{-1}_f).\intr H\,,\,\,{\rm and}\,\,\,
\intl{\underline H}=\intr{\underline H}\circ S_{\underline H}\,\,,
\]
where $\,a_f:=(\alpha_f\otimes f)\circ{\mathcal R}_A\,$.
\end{proposition}

\begin{proof}
It is verified directly that
$\,\cointl H:(\Int H,\id_{(\Int H)}\otimes\alpha_f)\to(H,\mu^H_\ad)\,$
is a morphism in $\,\C_\cO{A,\overline A}\,$.  This fact and the
universal property of $\,\cointl H\,$ imply that $\,\cointl H\,$ is a
left integral in $\,\underline H\,$.
\par
$\,(\intr{\underline H}\otimes \id_H)\circ\Delta_{\underline H}\,$ is
equal to $\,\intr{\underline H}\otimes \eta_H\,$ by the property of
integrals, which in turn is found to be equal to the composition
\[
H\;\xra{\,\Delta_H\,}\;H\otimes H\;
\xra{\,\id_H\otimes a^{-1}_f\otimes\id_H\,}\;
H\otimes H\otimes H\;
\xra{\,\intr{\underline H}\otimes\mu_H\,}\;\Int H\otimes H
\]
using the fact that $\,\intr{\underline H}\,$ is a morphism in
$\,\C_\cO{A,\overline A}\,$.  Comparing these two expressions we obtain
that $\,a_f.\intr{\underline H}\,$ is a right integral in $\,H\,$.
\end{proof}

\section{Applications of integrals}
\label{sec-appl-int}

The first part of this section is devoted to the construction of  the
 braided version of the Hopf algebra structure on the exterior algebra
$\,\bigoplus_{j} X^{\wedge j}\,$ for a given object, $\,X\,$, in a
braided category, $\,C\,$. As in the symmetric case the object of the
integral will turn out to be  the last non-vanishing power
$\,X^{\wedge n-1}\,$, and all morphisms are canonical ones.

In the second part of this section we shall explicitly construct and
 investigate the equivalence of the category of Hopf $\,H$-bimodules
and the category of Hopf $\,H^{\wedge}$-bimodules for a Hopf algebra,
$\,H\,$, in a rigid, braided category. From this we will be able to
draw two more proofs of the generalized Radford formula.

\subsection*{Integrals on external Hopf algebra}

For a given invertible solution, $\,R\in\End(V^{\otimes 2})\,$, of the
Yang-Baxter equation one can turn the tensor algebra $T(V)$ into a Hopf
algebra in the corresponding braided category, see \cite{Ma:fre}.  It
is shown in \cite{BesDra:DifCal} that, similarly, for a given object,
$\,X\,$, in an abelian braided monoidal category, $\,\mathcal C\,$, the
collection $\,\{X^{\otimes n}\}_{n\in\NN}\,$ becomes a Hopf algebra,
$\,T_{\C}(X)\,$, in the category $\,{\C}^{\NN}\,$ of graded spaces.
Moreover, because of categorical duality, there exists another Hopf
algebra, $\,{\overset\circ T}_{\C}(X)\,$, with the same underlying
object $\,\{X^{\otimes n}\}_{n\in\NN}\,$. The braided analog of the
antisymmetrizer $\,\hat A_X:T_{\C}(X)\to{\overset\circ T}_{\C}(X)\,$
is a Hopf algebra morphism, and its image
$\,T^\wedge_{\C}(X)=\{X^{\wedge n}\}_{n\in\NN}\,$ is a Hopf algebra.
Here we recall some necessary results from \cite{BesDra:DifCal} and
then restrict ourselves to the case when $\,T^\wedge_{\C}(X)\,$ has a
finite number of non-zero components and becomes an object of a rigid
category.  We determine explicitly the integrals on
$\,T^\wedge_{\C}(X)\,$ and derive from general integration theory some
properties of the objects $\,X^{\wedge n}\,$.

First let us outline the braided combinatorics, as in
 \cite{Ma:fre,BesDra:DifCal}, in an additive braided monoidal category,
$\,\C\,$.  The canonical epimorphism of the braid group $\,B_n\,$ into
the permutation group $\,S_n\,$ admits a section
$\,S_n \xrightarrow{\widehat{\;}} B_n\,$, which maps the standard
generators to the standard generators and is  uniquely determined on
other elements by the property
\begin{equation*}
\widehat{\sigma_1\sigma_2}=\widehat{\sigma_1}\widehat{\sigma_2}
\qquad\text{if}\qquad
\ell(\sigma_1\sigma_2)=\ell(\sigma_1)+\ell(\sigma_2)\,,
\end{equation*}
where $\,\ell(\sigma)\,$ is the length (of the minimal decomposition)
of $\,\sigma\,$.  Since the category $\,\C\,$ is braided, there is a
canonical mapping $\,B_n\to\End(X^{\otimes\,n})\,$.  The image of
$\,\sigma\in S_n\,$ under the composition
$\,S_n\to B_n\to\End(X^{\otimes\,n})\,$ is denoted by
$\,\sigma_{\C}(X)\,$.  Now let $\,\pi=(j_1,\ldots,j_r)\,$ be any
$\,\NN $-partition of $\,j\,$, i.e., $\,j=j_1+\cdots +j_r\,$ and
$\,j_1,\ldots j_r\in\NN \,$.  We consider the shuffle permutations
$\,S^j_{\pi}\subset S_j\,$. For every $\,l\in\{1,\ldots,r\}\,$ they are
mapping the $\,j_l\,$ elements in $\,\{1,\ldots ,j\}\,$ to
$\,\{j_1+\cdots+j_{l-1}+1,\ldots ,j_1+\cdots + j_l\}\,$ without
changing their order. The set of inverse permutations of
$\,S^j_{\pi}\,$ is denoted by $\,S_j^{\pi}\,$.  For every partition
$\,\pi=(j_1,\ldots ,j_r)\,$ of $\,j\,$, any object $\,X\in\Obj({\C})\,$
and any $\,\lambda\in\Aut(\1)\,$ we define {\em braided multinomials}
as endomorphisms of $\,X^{\otimes\,j}\,$ in $\,{\C}\,$  by the
following identities
\begin{equation}\label{multinom1}
\begin{split}
&\Big[\Atop\pi{j}\Big\vert X;\lambda\Big]\;=\;
\Big[\Atop{j_1\dots j_r}{j}\Big\vert X;\lambda\Big]\; :=\;
\sum_{\sigma\in S_j^\pi}\lambda^{\ell(\sigma)}\,\sigma_{\C}(X)\,,\\
&\Big[\Atop{j}{\pi}\Big\vert X;\lambda\Big]\;=\;
\Big[\Atop{j}{j_1\dots j_r}\Big\vert X;\lambda\Big]\; :=\;
\sum_{\sigma\in S^j_\pi}\lambda^{\ell(\sigma)}\,\sigma_{\C}(X)\,.
\end{split}
\end{equation}

\begin{proposition}\label{combi}
Let $\,\pi=(j_1,\dots, j_r)\,$ be a partition of $\,j\,$, and let
$\,\pi_k=(j_1^k,\dots ,j_{r_k}^k)\,$ be a partition of $\,j_k\,$ for
any $\,k\in\{1,\dots ,r\}\,$. Then with the notation from above the
following formulae hold true:
\begin{equation}
\label{combi-eq1}
\Big[\Atop{(\pi_1,\pi_2,\dots ,\pi_r)}{j}\Big\vert X;\lambda\Big]
=\Big[\Atop\pi{j}\Big\vert X;\lambda\Big]\circ
   \bigg(\Big[\Atop{\pi_1}{j_1}\Big\vert X;\lambda\Big]\otimes
       \Big[\Atop{\pi_2}{ j_2}\Big\vert X;\lambda\Big]\otimes\cdots
      \otimes
         \Big[\Atop{\pi_r}{j_r}\Big\vert X;\lambda\Big]\bigg)
\end{equation}
and
\begin{equation}
  \label{combi-eq2}
\Big[\Atop{j }{(\pi_1,\pi_2,\dots ,\pi_r)}\Big\vert X;\lambda\Big]
=\bigg(\Big[\Atop{j_1}{\pi_1}\Big\vert X;\lambda\Big]\otimes
      \Big[\Atop{j_2}{\pi_2}\Big\vert X;\lambda\Big]\otimes\cdots\otimes
         \Big[\Atop{j_r}{\pi_r}\Big\vert X;\lambda\Big]\bigg)\circ
         \Big[\Atop{j}{\pi}\Big\vert X;\lambda\Big]\quad,
\end{equation}
where (\ref{combi-eq2}) is the dual version of (\ref{combi-eq1}).
\end{proposition}

For any $\,0\le k\le j\,$ we shall use abbreviated
notations as follows:
$\,\big[\Atop{j}{k}\big\vert X;\lambda\big]
:=\big[\Atop{j} {(k , j-k)}\big\vert X;\lambda\big]\,$, \
$\,\big[\Atop{k}{j}\big\vert X;\lambda\big]:=
\big[\Atop{(k , j-k)}{j}\big\vert X;\lambda\big]\,$ and
$\,\big[j\big\vert X;\lambda\big]!:=
\big[\Atop{j}{1\dots 1} \big\vert X;\lambda\big]=
\big[\Atop{1\dots 1}{j}\big\vert X;\lambda\big]\,$.
In this language  we then derive by successive application of
eq.~\eqref{combi-eq2} the identities
\begin{equation}\label{combi-eq3}
\begin{split}
[j]! &= \bigg(\id_{X^{\otimes\,j-2}}\otimes\Big[\Atop{2}{1}\Big]\bigg)
  \circ\bigg(\id_{X^{\otimes\,j-3}}\otimes\Big[\Atop{3}{1}\Big]\bigg)
  \circ\dots\circ\Big[\Atop{j}{1}\Big]\,,\\
[j]! &= \bigg(\Big[\Atop{2}{1}\Big]\otimes\id_{X_{\otimes\,j-2}}\bigg)
   \circ\bigg(\Big[\Atop{3}{2}\Big]\otimes\id_{X^{\otimes\,j-3}}\bigg)
  \circ\dots\circ\Big[\Atop{j}{j-1}\Big]
\end{split}
\end{equation}
and
\begin{equation}\label{combi-eq4}
[j+k]! = ([j]!\otimes [k]!)\circ \Big[\Atop{j+k}{j}\Big]\,,
\end{equation}

\begin{definition}\label{3-6-D1}
Let $\,\I\,$ be $\,\NN\,$ or $\,\Zn\,$, considered as a discrete
 category.  We consider the functor category $\,{\C}^\I\,$ as an \/
$\,\I$-graded category over $\,{\C}\,$.  The objects of $\,{\C}^\I\,$
are given by $\,\I$-tuples, $\,\hat X=(X_0, X_1,\ldots)\,$, of objects,
$\,X_j\in \Obj({\C})\,$, where $\,j\in \I\,$. The morphisms of
$\,{\C}^\I\,$ are of the form
$\,\hat f=(f_0,f_1,\ldots):\,\hat X\to \hat Y\,$, where
$\,f_j:\,X_j\to Y_j\,$ is a morphism in $\,{\C}\,$ for all
$\,j\in \I\,$.  The category $\,{\C}^\I\,$ naturally inherits a braided
monoidal structure from $\,\C\,$. The unit object in $\,{\C}^\I\,$ is
given by $\,\Hat\1=(\1,0,0,\ldots)\,$, the tensor product is defined by
$\,(\hat X\otimes\hat Y)_n=\bigoplus_{k+l=n}X_k\otimes Y_l\,$, where
$\,n,k,l\in\I\,$, and the tensor product of\/ $\,\I$-graded morphisms
is built analogously.  Besides the natural braided structure, the
category ${\C}^\I$ actually admits a {\em family} of braidings given by
$\,(\hat \Psi_{\hat X,\hat Y}^{(\lambda)})_n=
\bigoplus_{k+l=n}\lambda^{k\,l}\Psi_{X_k,Y_l}\,$,
for any $\,\lambda\in{\rm Aut}_{\C}(\1)\,$, if $\I=\NN$, and for any
$\,\lambda\in{\rm Aut}_{\C}(\1)\,$ with $\lambda^n=1$, if $\I=\Zn$.
\end{definition}

\begin{proposition}\label{3-4f-D11II}\

\begin{enumerate}
\item For a  given object, $\,X\,$, in an additive braided category,
$\,{\C}\,$, the tensor algebra $\{X^{\otimes n}\}$ in ${\C}^{\NN }$
admits two Hopf algebra structures,
$\,T_{\C}(X)=\bigl(\{X^{\otimes n}\},\mu,\eta,\Delta,\varepsilon\bigr)\,$
and
$\,{\overset\circ T}_{\C}(X)=
 \bigl(\{X^{\otimes n}\},\overset\circ\mu,\overset\circ\eta,
 \overset\circ\Delta,\overset\circ\varepsilon\bigr)\,$,
which are given as follows:
\begin{eqnarray}
\mu_{n,m} \cong \id_{X^{\otimes\,n+m}}&:&
X^{\otimes\,n}\otimes X^{\otimes\,m}\to X^{\otimes\,n+m}\,,
\nonumber\\
\Delta_{n,m} \cong \Big[\Atop{n+m}{n}\Big\vert X;\lambda\Big]&:&
 X^{\otimes\,n+m}\to X^{\otimes\,n}\otimes X^{\otimes\,m}\,,
\nonumber\\
{\overset\circ\mu}_{n,m} \cong
\Big[\Atop{n}{n+m}\Big\vert X;\lambda\Big]&:&
X^{\otimes\,n}\otimes X^{\otimes\,m}\to X^{\otimes\,n+m}\,,
\nonumber\\
{\overset\circ\Delta}_{n,m} \cong \id_{X^{\otimes\,n+m}}&:&
X^{\otimes\,n+m}\to X^{\otimes\,n}\otimes X^{\otimes\,m}\,,
\nonumber\\
\eta_n={\overset\circ\eta}_n=
\begin{cases}\id_\1\,,&n=1\,, \\
             0\,,&n\ne 1\,,\end{cases}
\quad&&\quad
\varepsilon_n={\overset\circ\varepsilon}_n=
\begin{cases}\id_\1\,,&n=1\,,   \\
             0\,,&n\ne 1\,,\end{cases}
\nonumber\\
S_n=\overset\circ S_n=
(-1)^n\lambda^{{}^{\binom{n}{2}}}\,(\sigma_n^0)_{\C}(X)
&:& X^{\otimes\,n}\to X^{\otimes\,n}\,,
\qquad\text{where}\quad
\sigma_n^0=\left(\begin{smallmatrix}1 & \ldots & n\\ n & \ldots &1
\end{smallmatrix}\right)\,.
\nonumber
\end{eqnarray}
\item The morphism
$\,\hat A_X:=\bigl([n\vert X]!\bigr)_{n\in\NN}: \,T_{\C}(X)\to
{\overset\circ T}_{\C}(X)\,$ is a Hopf algebra map.
\end{enumerate}
\end{proposition}

Now let $\,\C\,$ be a rigid braided monoidal abelian category with a
biadditive tensor product.  In this case the unit object $\,\1\,$ is
semisimple, see \cite{DelMil},   so that  we can suppose without
restriction of generality  that $\,\End_{\C}(\1)\,$ is a field.
Rigidity implies that the tensor product is exact. In this situation we
have that for any Hopf algebra morphism, $\,f\,$, the intermediate
 object of an epi-mono decomposition, $\,f=\im f\circ\coim f\,$, can be
equipped with a Hopf algebra structure, such that both, $\,\im f\,$ and
$\,\coim f\,$, are Hopf algebra morphisms. This allows us to make the
following definition:

\begin{definition} \label{anti-tens-alg}
We will denote by $\,T^\wedge_{\C}(X)=\{X^{\wedge n}\}\,$
the Hopf algebra determined by the epi-mono decomposition
\begin{equation}
\label{Anti-epi-mono}
\hat A_X=\{T_{\C}(X)\xra{\coim{\hat A_X}} T^\wedge_{\C}(X)
           \xra {\im{\hat A_X}}{\overset\circ T}_{\C}(X)\}\,.
\end{equation}
\end{definition}

Let us suppose that $\,\lambda^n=1\,$, $\,[n\vert X]!=0\,$, and
$\,[n-1\vert X]!\ne0\,$.  In this case also $\,X^{\wedge k}=0\,$ for
any $\,k\ge n\,$, and $\,T^\wedge_{\C}(X)\,$ can be considered as a
Hopf algebra in the rigid braided monoidal abelian category
$\,{\C}^{\Zn}\,$.

\begin{proposition}
Under the previous assumptions the integrals
$\,\intl{T^\wedge_{\C}(X)}\,$, $\,\intr{T^\wedge_{\C}(X)}\,$,
$\,\cointl{T^\wedge_{\C}(X)}\,$  and $\,\cointr{T^\wedge_{\C}(X)}\,$
are all given by  the following expressions:
\begin{equation*}
\Int\bigl(T^\wedge_{\C}(X)\bigr)_k=
 \begin{cases}0\,,&\text{if \  $k\ne n-1$\,,}\\
              X^{\wedge(n-1)}\,,&\text{if \ $k=n-1$}\,,\end{cases}
\qquad\quad
({\textstyle \int})_k=
 \begin{cases}0\,,&\text{if \  $k\ne n-1$\,,}\\
              \id_{X^{\wedge(n-1)}}\,,&\text{if \ $k=n-1$\,.}\end{cases}
\end{equation*}
\end{proposition}
\begin{proof}
We have that
$\,\PPi^\bullet_\bullet\bigl(T^\wedge_{\C}(X)\bigr)_{n-1}=
 \id_{X^{\wedge(n-1)}}\,$,
because
$\,\m^{T^\wedge_{\C}(X)}_{k,\ell}=0\,$ for $\,k+\ell\ge n\,$.
The condition that $\,\End_{\C}\bigl(\Int(T^\wedge_{\C}(X))\bigr)\,$
has no nontrivial idempotents implies that $\,\Int(T^\wedge_{\C}(X))\,$
has only one nonzero component.
\end{proof}

\begin{corollary}
Under the above assumptions $\,X^{\wedge(n-1)}\,$ is an invertible
object in $\,\C\,$.  The multiplication in $\,T_{\C}^\wedge(X)\,$
defines side-invertible pairings,
$\,\mu_{k,n-k-1}:X^{\wedge k}\otimes
X^{\wedge (n-k-1)}\to X^{\wedge (n-1)}\,$.
I.e., it  induces isomorphisms,
$\,\,(X^*)^{\wedge(n-k-1)}\otimes X^{\wedge(n-1)}\simeq X^{\wedge k}\,$,
for  $\,k\in\Zn\,$ and any object $\,X^*\,$ dual to $\,X\,$.
\end{corollary}

The Hopf algebra from Example~\ref{non-triv-int} was obtained as
$\,T^\wedge(k)\,$ for $\,\lambda=q\,$ a primitive $\,n^{\text{th}}\,$
root of unit.  Our results show similarities between the Hopf algebras
 $\,T^\wedge(k)\,$ and $\,T^\wedge(X)\,$ in the most general case.

\begin{remark}
Note that if $\,[n+1\vert X]!=0\,$ and  $\,[m+1\vert X]!=0\,$  then
also $\,[n+m+1\vert X\oplus Y]!=0\,$. This follows from the fact that
the matrix element of $\,[n+m+1\vert X\oplus Y]!\,$, which is a
morphism from $\,X^{\otimes k}\otimes Y^{\otimes\ell}\,$ to
$\,X^{\otimes k}\otimes Y^{\otimes\ell}\,$, equals to
$\,[k\vert X]!\otimes[\ell\vert Y]!\,$ and other matrix elements are
either of such type (up to isomorphism) or zero.
\end{remark}

\begin{remark}
Suppose $\,X\,$ lives in the category $\,\DY{\C}^H_H\,$ of crossed
modules over a Hopf algebra, $\,H\,$, and $\,\d:H\to H\ltimes X\,$
defines a first order differential calculus.  Then
$\,H\ltimes T^\wedge_{\C(X)}\,$ is a Hopf algebra in $\,{\C}^{\Zn}\,$,
and $\,\d\,$ has a unique extension, $\,\d^{\wedge}\,$, on
$\,H\ltimes T^\wedge_{\C(X)}\,$, which turns
$\,H\ltimes T^\wedge_{\C(X)}\,$ into a differential Hopf algebra
(braided De Rham complex) \cite{BesDra:DifCal}.  Integrals on
$\,H\ltimes T^\wedge_{\C(X)}\,$ are calculated by the formulae
(\ref{int-cross-prod}) for integrals on cross products.  They are
nonzero only on one component $\,H\ltimes X^{\wedge(n-1)}\,$ and play
the role of the integration over a volume form.
\end{remark}

\subsection*{Duality for Hopf bimodules}

In this section we shall describe how the previous results can be
applied in order to directly  construct the equivalence functor between
 the categories $\,\sideset{^H_H}{^H_H}\C\,$ and
$\,\sideset{^{(H^\vee)}_{(H^\vee)}}{^{(H^\vee)}_{(H^\vee)}}\C\,$
(see \cite{B:hopf93} for unbraided case).  The existence of such an
equivalence follows from the fact that the category
$\,\sideset{^H_H}{^H_H}\C\,$ is equivalent to the category of crossed
modules $\,\DY{\C}^H_H\,$. The equivalence of the categories
$\,\DY{\C}^H_H\,$ and $\,\DY{\C}^{H^\vee}_{H^\vee}\,$ is described in
\cite{Bes:cross}.  These and the following considerations lead to
alternative proofs of the generalized Radford formula:

\begin{proposition}
\label{Prop-X-over-dual}
There exists an equivalence,
$\,\sideset{^H_H}{^H_H}\C\to
 \sideset{^{(H^\vee)}_{(H^\vee)}}{^{(H^\vee)}_{(H^\vee)}}\C\,$,
of monoidal categories, which converts a Hopf $H$\n-bimodule $X$ into
the Hopf $H^\vee$\n-bimodule $\,X\otimes\Int H^\vee\,$ with actions and
coactions defined via  Figure~\ref{Fig-X-over-dual}.
\end{proposition}

\begin{figure}
\[
\ba{ccccccc}
\ba{l}
\object{H^\vee}\Step\object{X}\Step\object{\Int{H^\vee}}\\
\begin{tangle}\xx\Step\id \\
       \id\Step\hdcd\step\id \\
       \rd\step\S\step\id\step\id \\
       \id\step\hev\step\QQ{\alpha^{-1}_{H^\vee}}\step\id\end{tangle}
\ea
&\quad&
\ba{l}
\step[.6]\object{X}\step[1.5]\object{\Int{H^\vee}}\step[1.6]\object{H^\vee}\\
\begin{tangle}\ld\step\hxx \\
       \hx\step\SS\step\id \\
       \hh\id\step\ev\step\id\end{tangle}
\ea
&\quad&
\ba{l}
\begin{tangle}\hh\coev\step\id\step\id \\
       \id\step\lu\step\id\end{tangle}\\
\object{H^\vee}\step[1.8]\object{X}\step[1.5]\object{\Int{H^\vee}}
\ea
&\quad&
\ba{l}
\begin{tangle}\id\step\id\step\Q{a_{H^\vee}}\step\hcoev \\
       \id\step\id\step\hdcu\step\O{S^2} \\
       \id\step\d\step\hxx \\
       \d\step\hxx\step\id \\
       \step\ru\step\id\step\id\end{tangle}\\
\step[.9]\object{X}\step[1.6]\object{\Int{H^\vee}}\step[1.6]\object{H^\vee}
\ea
\\
\hbox{\scriptsize $\mu^{(X\otimes\Int{H^\vee})}_\ell$}&&
\hbox{\scriptsize $\mu^{(X\otimes\Int{H^\vee})}_r$}&&
\hbox{\scriptsize $\Delta^{(X\otimes\Int{H^\vee})}_\ell$}&&
\hbox{\scriptsize $\Delta^{(X\otimes\Int{H^\vee})}_r$}
\ea
\]
\caption{$X\otimes\Int H^\vee\in
       \Obj(\sideset{^{H^\vee}_{H^\vee}}{^{H^\vee}_{H^\vee}}\C)$}
\label{Fig-X-over-dual}
\end{figure}

Before we can prove Proposition~\ref{Prop-X-over-dual} we need a few
technical results about tensor products on the $\,H$-module categories.
The first observations on the existence and properties of the braided
tensor product actions $\,\mu^{X\otimes Y}\,$ are straightforward to
prove:

\begin{lemma}
\label{Hopf-otimes}\

\begin{itemize}
\item
Let
$\,\left( X, \mu^{X}_{\ell}, \mu^{X}_{r} \right)
             \in\Obj\left( \sideset{_H}{_H}\C \right)\,$ and
$\,\left( Y, \mu^{Y}_{\ell}, \mu^{Y}_{r},
       \Delta^{Y}_{r} \right)
             \in\Obj\left( \sideset{_H}{^H_H}\C \right)\,$.

Then
$\,\left( X\otimes Y, \mu^{X\otimes Y}_{\ell}, \mu^{X\otimes Y}_{r},
        \id_X\otimes\Delta^Y_r \right)
             \in\Obj\left( \sideset{_H}{^H_H}\C \right)\,$,
where the underlying module structures
$\,\mu^{X\otimes Y}_{\ell}\,$ and $\,\mu^{X\otimes Y}_{r}\,$
are the braided tensor product ones:
\begin{multline*}
\mu^{X\otimes Y}_{\ell} :\, H\otimes X\otimes Y\,
\xra{\,\Delta\otimes\id_{X\otimes Y}\,}\, H\otimes H\otimes X\otimes Y \\
\xra{\,\id_H\otimes\Psi_{H,X}\otimes\id_H\,} \,H\otimes X\otimes H\otimes Y
\,\xra{\,\mu^{X}_{\ell}\otimes\mu^{Y}_{\ell}\,} \,X\otimes Y\,.
\end{multline*}
\item
Let
$\,\left( X, \mu^{X}_{\ell}, \mu^{X}_{r},
       \Delta^{X}_{\ell} \right)
             \in\,\Obj\left( \sideset{^H_H}{_H}\C \right)\,$ and
$\,\left( Y, \mu^{Y}_{\ell}, \mu^{Y}_{r},
       \Delta^{Y}_{r} \right)
             \in\,\Obj\left( \sideset{_H}{^H_H}\C \right)\,$.

Then
$\,\left( X\otimes Y, \mu^{X\otimes Y}_{\ell}, \mu^{X\otimes Y}_{r},
       \Delta^X_\ell\otimes\id_Y, \id_X\otimes\Delta^Y_r \right)
             \in\,\Obj\left( \sideset{^H_H}{^H_H}\C \right)\,$.
\end{itemize}
\end{lemma}


Lemma~\ref{Hopf-otimes} further allows us to properly define the
tensor product of bimodules $\,(X,Y)\to X_\bullet\tens_H Y\,$,
which may be viewed as an action of $\,\sideset{^H_H}{^H_H}\C\,$
on $\,\sideset{^{H\pti}_{H\pti}}{^{H\pti}_{H\pti}}\C\,$:

\begin{proposition}
Let $\,X\,$ be a Hopf $\,H$-bimodule and $\,Y\,$ be a Hopf
$\,H\pti$\n-bimodule.  View $\,Y\,$ as a left $\,H$\n-module via the
action $\,{}^\Vee\mu_\ell\,$ from (\ref{Dual-act}).  Then there exists
the tensor $\,X\tens_H Y \in\C\,$ product of $\,H$\n-modules.  There is
a unique Hopf $\,H\pti$\n-bimodule structure $\,X_{\bullet}\tens_H Y\,$
on this object, such that the canonical projection
$\,\lambda: \,X_\bullet\tens Y \to X_\bullet\tens_H Y\,$ is a
homomorphism of $\,H\pti$\n-bimodules. Here $\,X_\bullet\,$ is obtained
via Corollary~\ref{two-fold-Hopf-fun} (the explicit
$\,H\pti$\n-(co)module structures for the underlying object of $\,X\,$
are shown in Figure~\ref{Fig-over-dual} ), and the Hopf bimodule
structure on the tensor product is defined by Lemma~\ref{Hopf-otimes}.
The monoidal category of $\,H$\n-bimodules equipped with the tensor
product $\,\text-\tens_H\text-\,$ of bimodules acts on the category of
$\,H\pti$\n-bimodules via the functor
$\,\sideset{^H_H}{^H_H}\C \times
\sideset{^{H\pti}_{H\pti}}{^{H\pti}_{H\pti}}\C \to
\sideset{^{H\pti}_{H\pti}}{^{H\pti}_{H\pti}}\C\,$:
$\,(X,Y)\to X_\bullet\tens_H Y\,$.
\end{proposition}

\begin{proof}
The $\,H$-module $\,Y\,$ admits in fact the structure $\,Y^\bullet\,$
of a left Hopf $\,H$\n-module given by Figure~\ref{Fig-over-dual}.
According to \cite{BesDra:hopf} there exists a tensor product,
$\,X\tens_H Y\,$, together with the canonical projection
$\,\lambda^H_{X,Y}:\,X\otimes Y\to X\otimes_H Y\,$ given by the
composition
\[
\lambda^H_{X,Y}\,:\;
X\otimes Y\xra{\id_X\otimes\Delta^Y_\ell}X\otimes H\otimes Y
\xra{\mu^X_r\otimes {}_Yp}X\otimes\left( {}_HY \right)\,,
\]
with the following property:

For any $\,f:X\otimes Y\to Z\,$, for which both compositions
\[
X\otimes H\otimes Y
\genfrac{}{}{0pt}{2}
{\xrightarrow[\hphantom{\id_X\otimes\mu^Y_\ell}]{\mu^X_r\otimes\id_Y}}
{\xrightarrow[\id_X\otimes\mu^Y_\ell]{\hphantom{\mu^X_r\otimes\id_Y}}}
X\otimes Y \xra{f} Z
\]
are equal, there is a unique $\,\hat f:\,X\otimes({}_HY)\to Z\,$,
such that
$\,f=\hat f\circ\lambda^H_{X,Y}\,$. It is given by the composition
$\,\hat f:\,X\otimes({}_HY)\xra{\id_X\otimes{}_Yi}X\otimes Y\xra{f}Z\,$.
Moreover, $\lambda$ is a split epimorphism with a splitting as follows:
\[ \theta: X\tens({}_HY) \xra{\id_X\tens {}_Yi} X\tens Y
\xra{\id_X\tens \Delta^Y_\ell} X\tens H\tens Y \xra{\id_X\tens S\tens\id_Y}
X\tens H\tens Y \xra{\mu_r^X\tens\id_Y} X\tens Y\; .
\]
\noindent
One can check that both of the following compositions coincide:
\[ H\pti\tens X\tens H\tens Y\,
\genfrac{}{}{0pt}{2}
{\xra[\hphantom{\id_{H\pti}\tens\mu^X_r\tens\id_Y}]%
{\id_{H\pti\tens X}\otimes\mu^Y_\ell}}
{\xrightarrow[\id_{H\pti}\tens\mu^X_r\tens\id_Y]%
{\hphantom{\id_{H\pti\tens X}\otimes\mu^Y_\ell}}}\,
H\pti\tens X\tens Y\, \xra{\mu_\ell}\, X\tens Y\, \xra\lambda\, X\tens_H Y
\quad .\]
Therefore, there exists a unique structure of a left $\,H\pti$\n-module
on $\,X\tens_H Y\,$, such that $\,\lambda\,$ is a homomorphism of left
$\,H\pti$\n-modules. It must coincide with
\[
H\pti\tens X\tens_H Y \,\xra{\,\id_H\tens\theta\,} \,H\pti\tens X\tens Y
\,\xra{\,\mu_\ell\,}\, X\tens Y \,\xra\lambda\, X\tens_H Y \,\,\,\,.
\]
The right $\,H\pti$\n-module structure is treated analogously.

The existence of a left $\,H\pti$-comodule structure on
$\,X\tens_H Y\,$ follows from the fact that the left coaction of
$\,H\pti\,$ on $\,X\,$ commutes with the right action of $\,H\,$. It
must be given by the formula
\[
X\tens_H Y \,\xra\theta \,X\tens Y \,\xra{\,\Delta_\ell^X\tens\id_Y\,}\,
H\pti\tens X\tens Y \,\xra{\,\id_{H\pti}\tens\lambda\,}\,
H\pti\tens(X\tens_H Y)
\,\,\,\, .
\]
Similarly, the existence of the right $\,H\pti$\n-comodule structure on
$\,X\tens_H Y\,$ follows from the fact that the right coaction of
$\,H\pti\,$ on $\,Y\,$ commutes with the left action of $\,H\,$. Thus
$\,\lambda\,$ is a morphism of left and right modules and comodules
over $\,H\pti\,$. Moreover, it is a split epimorphism so that
$\,X\tens_H Y\,$ is a Hopf $\,H\pti$\n-bimodule.

One can check that the canonical associativity isomorphism
$\,(X\tens_H Y)\tens_H Z \simeq X\tens_H(Y\tens_H Z)\,$ for a pair of
Hopf $\,H$\n-bimodules, $\,X\,$ and $\,Y\,$, and a Hopf
$\,H\pti$\n-bimodule, $\,Z\,$, is an isomorphism of
$\,H\pti$\n-bimodules, where $\,X\tens_H Y\,$ is the standard tensor
product of Hopf bimodules.
\end{proof}

\begin{proof}[Proof of Proposition~\ref{Prop-X-over-dual}]
Applying the above proposition to the regular Hopf bimodule
$\,Y = H\pti\,$ we get a Hopf $\,H\pti$\n-bimodule structure on
$\,X\otimes\Int{H^\vee}\,$.  Indeed, the chosen left $\,H$\n-module
structure of $\,Y=H^\vee\,$ gives
$\,\Pi^\ell_\ell(Y)=\PPi^\ell_r(H^\vee)\,$,
$\,{}_Yp=\intl{H^\vee}\,$, $\,{}_Yi=(\clr)^{-1}\cdot\cointr{H^\vee}\,$,
and
$\,\lambda^H_{X,H^\vee}=
\{X\otimes H^\vee\xra{\id_X\otimes{}_{H^\vee}\overline{\mathcal F}}
X\otimes H\otimes\Int H^\vee\xra{\mu_r\otimes\Int H^\vee}
X\otimes\Int H^\vee \}\,$.
The splitting $\,\theta\,$ simplifies to
$\,(\clr)^{-1}\cdot\id_X\otimes\cointr{H^\vee}:\,
X\otimes\Int{H^\vee} \to X\otimes H^\vee\,$ by the definition of
integrals.  Hence the Hopf $\,H^\vee$-bimodule structure on
$\,X\otimes\Int{H^\vee}\,$ is given by the compositions
{\allowdisplaybreaks
\begin{align}
\mu^{X\otimes\Int{H^\vee}}_\ell\,&:& H^\vee\otimes X\otimes\Int{H^\vee}
\xra{(\clr)^{-1}\cdot\id_{H^\vee\otimes X}\otimes\cointr{H^\vee}}
H^\vee\otimes X\otimes H^\vee\qquad \quad \ \notag \\
&& \hspace*{3cm}
\xra{\mu^{X\otimes H^\vee}_{\ell}}
X\otimes H^\vee
\xra{\lambda^H_{X,H^\vee}}
X\otimes\Int{H^\vee}\,, \notag
\\
\mu^{X\otimes\Int{H^\vee}}_r      \,&:& X\otimes\Int{H^\vee}\otimes H^\vee
\xra{(\clr)^{-1}\cdot\id_X\otimes\cointr{H^\vee}\otimes\id_{H^\vee}}
X\otimes H^\vee\otimes H^\vee \notag \qquad \quad \\
&& \hspace*{3cm}
\xra{\mu^{X\otimes H^\vee}_{r}}
X\otimes H^\vee
\xra{\lambda^H_{X,H^\vee}}
X\otimes\Int{H^\vee}\,, \notag
\\
\Delta^{X\otimes\Int{H^\vee}}_\ell   \,&:& X\otimes\Int{H^\vee}
\xra{(\clr)^{-1}\cdot\id_X\otimes\cointr{H^\vee}}
X\otimes H^\vee
\xra{\tilde\Delta_\ell\otimes\id_{H^\vee}}
H^\vee\otimes X\otimes H^\vee \notag \\
&& \hspace*{3cm}
\xra{\id_{H^\vee}\otimes\lambda^H_{X,H^\vee}}
H^\vee\otimes X\otimes\Int{H^\vee}\,, \notag
\\
\Delta^{X\otimes\Int{H^\vee}}_r\,&:& X\otimes\Int{H^\vee}
\xra{(\clr)^{-1}\cdot\id_X\otimes\cointr{H^\vee}}
X\otimes H^\vee
\xra{\id_X\otimes\Delta_{H^\vee}}
X\otimes H^\vee\otimes H^\vee \notag \\
&& \hspace*{3cm}
\xra{\lambda^H_{X,H^\vee}\otimes\id_{H^\vee}}
X\otimes\Int{H^\vee}\otimes H^\vee\,.
\label{X-Int-Hopf-structure}
\end{align}
}
or in explicit diagrammatic form  in Figure~\ref{Fig-X-over-dual}.
\par
The dual construction $\,X^\bullet\square_H{}^\vee H\,$, where the
cotensor product is used instead of tensor product, gives a Hopf
bimodule structure on $\,X\otimes\Int{}^\vee H\,$ described by the same
diagrams only taken upside-down.  Specializing this construction to the
case of Hopf $\,H^\vee$-bimodules yields a functor,
$\,\sideset{^{H^\vee}_{H^\vee}}{^{H^\vee}_{H^\vee}}\C\to
\sideset{^H_H}{^H_H}\C\,$,
which together with the functor described above defines an equivalence
of categories.  Indeed, let us define an isomorphism
$\,\mu_r\circ(\id_X\otimes a):X\to X\,$ in $\,\C\,$, consider the
source space with the given Hopf bimodule structure, and introduce a
new Hopf bimodule structure, $\,\tilde X\,$, on the target space, which
turns this isomorphism into a Hopf bimodule morphism.  Direct
calculations show that subsequent applications of the two above
functors turn the Hopf $\,H$-bimodule $\,X\,$ into the Hopf
$\,H\,$-bimodule $\,X\otimes\Int H^\vee\otimes\Int H\,$, whose
underlying (co)modules are tensor products of the underlying
(co)modules of $\,\tilde X\,$, with the trivial (co)module
$\,\Int H^\vee\otimes\Int H\,$.
\end{proof}

\begin{remark}
Note that the $\,H^\vee$\n-(co)module structures on
$\,X\otimes\Int H^\vee\,$ presented in Figure~\ref{Fig-X-over-dual} are
the tensor product ones:
\begin{multline*}
\bigl(X\otimes\Int H^\vee\,,\,\,
 \mu^*_{\bullet\,\ell}\otimes\id_{\Int H^\vee}\,,\,\,
 (\mu_{\bullet\,r}\otimes\id_{\Int H^\vee})\circ
                          (\id_X\otimes\Psi_{\Int H^\vee,H^\vee})\,\,,\\
 \Delta_{\bullet\,\ell}\otimes\id_{\Int H^\vee}\,,\,\,
 (\id_X\otimes\Psi_{\Int H^\vee,H^\vee})\circ
               (\Delta^*_{\bullet\,r}\otimes\id_{\Int H^\vee}) \bigr)\,,
\end{multline*}
where the first factor is the underlying right-left Hopf
$H^\vee$\n-module $(X,\mu_{\bullet\,r},\Delta_{\bullet\,\ell})$ of
$X_\bullet$ (see Figure~\ref{Fig-over-dual}) with the modified left
action and right coaction:
\begin{eqnarray*}
\mu^*_{\bullet\,\ell}&:=\{H^\vee\otimes X\xra{\Delta_{H^\vee}\otimes\id_X}
H^\vee\otimes H^\vee\otimes X
\xra{\id_X\otimes\alpha^{-1}_{H^\vee}\otimes\id_X}
H^\vee\otimes X\xra{\mu_{\bullet\,\ell}}X  \}\,,
\\
\Delta^*_{\bullet\,r}&:=\{ X\xra{\Delta_{\bullet\,r}} X\otimes H^\vee
\xra{\id_X\otimes a_{H^\vee}\otimes(\monodromy{\Int H^\vee}{H^\vee})^{-1}}
X\otimes H^\vee\otimes H^\vee \xra{\id_X\otimes\mu_{H^\vee}}
X\otimes H^\vee\}\,,
\end{eqnarray*}
and $\,\Int H^\vee\,$ is considered with trivial (co)module structures.
The Hopf bimodule axioms for $\,X\otimes\Int H^\vee\,$ are equivalent to
the conditions that
$\,(X,\mu^*_{\bullet\,\ell},\mu_{\bullet\,r},\Delta_{\bullet\,\ell})\in
\Obj(\sideset{^{H^\vee}_{H^\vee}}{_{H^\vee}}\C)\,$,\
$\,(X,\mu_{\bullet\,r},\Delta_{\bullet\,\ell},\Delta^*_{\bullet\,r})\in
\Obj(\sideset{^{H^\vee}_{H^\vee}}{^{H^\vee}}\C)\,$,
and the following modified version of the right Hopf module axiom:
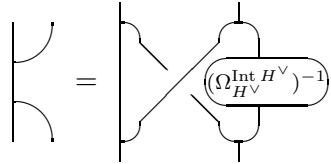
\begin{equation*}
\begin{tangle}\ru \\
       \hh\rd\end{tangle}
\enspace =\enspace
\begin{tangle}\hh\hrd\Step\cd \\
       \id\hstep\xx\step\hbx(0,2){%
\put(0,0){\line(0,1){0.4}}%
\put(0,1.6){\line(0,1){0.4}}%
\put(.2,1){\oval(3.1,1.2)}%
\put(-0.4,0.4){\makebox(1.2,1.2)[cc]%
{\scriptsize $(\monodromy{\Int{H^\vee}}{H^\vee})^{-1}$}}} \\
       \hh\hru\Step\cu\end{tangle}
\end{equation*}
\end{remark}

\begin{proof}[Second proof of Theorem~\ref{thm-Radford-formula}]
Analogs of Proposition~\ref{Prop-X-over-dual} are true if categories of
Hopf modules or two-fold Hopf modules are used instead of category of
Hopf $\,H$\n-bimodules. In particular, one can consider the regular
Hopf bimodule $\,H^\vee\,$ and put $\,X=(H^\vee)^\bullet\,$ or
$\,X={}_\bullet(H^\vee)\,$, considered as an object of
$\,\sideset{}{^H_H}\C\,$ with corresponding structures defined in
Figure~\ref{Fig-over-dual}. The left-right Hopf module axiom for
$\,X\otimes\Int H^\vee\,$ from Proposition~\ref{Prop-X-over-dual} is
equivalent to the left-right Hopf module axiom for
$\,(X,\mu^*_{\bullet\,\ell},\Delta^*_{\bullet\,r})\,$.
For the latter axiom
$\,{\rm L.H.S}={\rm R.H.S.}:
\,H^\vee\otimes H^\vee\to H^\vee\otimes H^\vee\,$
the identity
$\,(\varepsilon\otimes\id)\circ{\rm L.H.S}\circ(\eta\otimes\id)=
 (\varepsilon\otimes\id)\circ{\rm R.H.S.}\circ(\eta\otimes\id)\,$
is the generalized Radford formula for $\,H^\vee\,$.
\end{proof}

\begin{proof}[Third proof of Theorem~\ref{thm-Radford-formula}]
\begin{figure}
\[
\ba{ccc}
\ba{l}
\Step\object{Y} \\
 \begin{tangle}\hh\Step\id \\
        \hh\coev\step\id \\
        \hx\step\id \\
        \SS\step\lu \\
        \hh\id\Step\id\end{tangle} \\
\object{H_\op}\Step\object{Y}
\ea
&\qquad&
\ba{l}
\object{X}\step[3]\object{H\pti} \\
 \begin{tangle}\hh\id\step\coev\step\id \\
        \id\step\hx\step\id \\
        \id\step\SS\step\hcu \\
        \ru\step[1.5]\Ointl{{H\pti}}\end{tangle} \\
\object{X}\step[2.5]\object{\Int H\pti}
\ea
\\
\hbox{\scriptsize a) $Y\in\Obj{{}^{H_\op}_{H_\op}\C}$}&&
\hbox{\scriptsize b) $\bar\lambda^H_{X,H^\vee}$}
\ea
\]
\caption{}
\label{Fig-2-over-dual}
\end{figure}

There is another possibility to extend the left $\,H$-module structure
on
$\,Y\in\sideset{^{H\pti}_{H\pti}}{^{H\pti}_{H\pti}}\C\,$ to a left Hopf
module structure. Specifically, Figure~\ref{Fig-2-over-dual}a)
describes $\,Y\,$ as a left Hopf $\,H_\op$\n-module.  We can apply the
reasonings in the previous proof to this case as well. Now, for
$\,Y = H\pti\,$ we get $\,{}_Yp=\intl{H^\vee}\,$,
$\,{}_Yi=(\cll)^{-1}\cdot\cointl{H^\vee}\,$, and the canonical
projection
$\,\bar\lambda^H_{X,H^\vee}: X\tens H\pti \to X\tens \Int H\pti =
X\tens_H H\pti\,$ is presented on Figure~\ref{Fig-2-over-dual}b).  The
splitting $\,\bar\theta\,$ is equal to
$\,(\cll)^{-1}\cdot\id_X\otimes\cointl{H^\vee}\,$.  The formulae for
the actions and coactions of $\,H\pti\,$ on $\,X\tens \Int H\pti\,$ are
given by expressions similar to \eqref{X-Int-Hopf-structure} with
$\,\lambda\,$ and $\,\theta\,$ replaced by $\,\bar\lambda\,$ and
$\,\bar\theta\,$, respectively.  They are described explicitly by
Figure~\ref{Fig-second-over-dual}.

\begin{figure}
\[
\ba{ccccccc}
\ba{l}
\object{H^\vee}\Step\object{X}\Step\object{\Int{H^\vee}}\\
\begin{tangle}\S\Step\id\step\id \\
       \xx\step\id \\
       \rd\step\id\step\id \\
       \hh\id\step\ev\step\id\end{tangle} \\
\object{X}\step[3]\object{\Int H\pti}
\ea
&\quad&
\ba{l}
\step[.6]\object{X}\step[1.5]\object{\Int{H^\vee}}\step[1.6]\object{H^\vee}\\
\begin{tangle}\step\id\step\id\step\hcd \\
       \ld\step\hxx\step\QQ{\alpha} \\
       \hx\step\SS\step\id \\
       \hh\id\step\ev\step\id\end{tangle} \\
\object{X}\step[3]\object{\Int H\pti}
\ea
&\quad&
\ba{l}
\begin{tangle}\hcoev\step\id\step\id \\
       \id\step\lu\step\id \\
        \hh\id\Step\id\step\id\end{tangle}\\
\object{H^\vee}\step[1.8]\object{X}\step[1.5]\object{\Int{H^\vee}}
\ea
&\quad&
\ba{l}
\object{X}\step[3]\object{\Int H\pti} \\
\begin{tangle}\hh\id\step\coev\step\id \\
       \id\step\hx\step\id \\
       \id\step\O{S^{-2}}\step\hxx\step\Q{a_{H\pti}} \\
       \ru\step\id\step\hcu\end{tangle}\\
\step[0]\object{X}\step[2]\object{\Int{H^\vee}}\step[1.6]\object{H^\vee}
\ea
\\
\hbox{$\bar\mu^{(X\otimes\Int{H^\vee})}_\ell$}&&
\hbox{$\bar\mu^{(X\otimes\Int{H^\vee})}_r$}&&
\hbox{$\bar\Delta^{(X\otimes\Int{H^\vee})}_\ell$}&&
\hbox{$\bar\Delta^{(X\otimes\Int{H^\vee})}_r$}
\ea
\]
\caption{The second Hopf module structure of $X\tens\Int H\pti$}
\label{Fig-second-over-dual}
\end{figure}

Both presentations of $\,X\tens_H H\pti\,$ must agree. Therefore, there
exists an isomorphism,
$\,\phi: X\tens_H H\pti \to X\tens_H H\pti\,$ of Hopf
$H\pti$\n-bimodules, equipped with the first and the second structure,
such that
\[ \bar\lambda\, =\, \left( X\tens H\pti \xra\lambda X\tens_H H\pti
\xra\phi X\tens_H H\pti \right)\quad .\]
Clearly,
\[ \phi\, = \,\left( X\tens\Int H\pti \xra\theta X\tens H\pti
\xra{\bar\lambda} X\tens\Int H\pti \right) \,\,\,.\]
A calculation gives $\,\phi = (\text-.a^{-1}_H) \tens \id\,$. The
statement that $\,\phi\,$ is a homomorphism of $\,H\pti$\n-modules and
left $\,H\pti$\n-comodules adds nothing new to our knowledge. However,
the statement that $\,\phi\,$ is a homomorphism of right
$\,H\pti$\n-comodules implies (and is equivalent to) the generalized
Radford formula. This gives the third proof of
Theorem~\ref{thm-Radford-formula}.
\end{proof}

\begin{remark}
The left dual Hopf $\,H$-bimodule to
$\,X\in\Obj{\sideset{^H_H}{^H_H}\C}\,$ is obtained by application of
the functor $\,{}^\vee(\_)\,$ in $\,\C\,$ to the Hopf $H^\vee$-bimodule
$\,X_\bullet\otimes_HH^\vee\,$.  As a result we obtain an object,
$\,\Int H\otimes{}^\vee X\,$, equipped with tensor product
$H$\n-(co)module structures, where $\,\Int H\,$ is considered with
trivial (co)actions, and the structures of $\,{}^\bullet({}^\vee X)\,$
are modified by means of $\,a_H,\alpha_H,\monodromy{\Int H}{H}\,$.
\end{remark}

\addtolength{\baselineskip}{-1.5mm}


\bigskip

   {\sc\small Yu.B.: Bogolyubov Institute for Theoretical Physics,

         Metrologichna str., 14-b,
         Kiev 143, 252143 Ukraine}

       {\em  E-mail: }{ \tt mmtpitp@gluk.apc.org}
\medskip

{\sc\small T.K.: The Ohio State University, Department of Mathematics,

 Columbus, OH,
 USA }

 {\em E-mail: }{ \tt kerler@math.ohio-state.edu}
\medskip

 {\sc\small V.L.: Department of Mathematics,
  138 Cardwell Hall,

  Kansas State University,
  Manhattan, Kansas 66506-2602,
  U.S.A. }

 {\em E-mail: }{ \tt lub@math.ksu.edu}
 \medskip

 {\sc\small V.T.: Institut de Recherche Math\'ematique Avanc\'ee,
  Universit\'e Louis Pasteur,

  7, rue Ren\'e Descartes,
  67084 Strasbourg,
  France}

  {\em E-mail: }{ \tt turaev@math.u-strasbg.fr}

\end{document}